\documentclass{amsart}
\usepackage{amssymb,amsmath,amsthm,amscd}
\usepackage{leftidx}
\usepackage{graphicx}
\usepackage{color}
\usepackage{microtype}
\usepackage{dsfont}
\usepackage{hyperref}
\usepackage[usenames,dvipsnames,svgnames,table]{xcolor}
\usepackage{tikz}
\usepackage{comment}
\usepackage{caption}
\usepackage{subcaption}
\usepackage{xfrac,faktor}
\usepackage{epstopdf} 
\usepackage{float}
\usepackage{tikz-cd}
\usepackage[all]{xy}
\usepackage{pgfplots}
\pgfplotsset{compat=1.16}

\newtheorem{theorem}{Theorem}[section]
\newtheorem{prop}[theorem]{Proposition}
\newtheorem{lemma}[theorem]{Lemma}

\newtheorem*{theorem*}{Theorem}

\theoremstyle{definition}
\newtheorem{definition}[theorem]{Definition}
\newtheorem{eg}[theorem]{Example}
\newtheorem{rmk}[theorem]{Remark}

\theoremstyle{definition}

\newcommand{\C}{{\mathbb C}}

\newcommand{\natls}{{\mathbb N}}

\newcommand\AAA{{\mathcal A}}
\newcommand\BB{{\mathcal B}}
\newcommand\CC{{\mathcal C}}

\newcommand\DD{{\mathcal D}}
\newcommand\EE{{\mathcal E}}
\newcommand\FF{{\mathcal F}}
\newcommand\GG{{\mathcal G}}
\newcommand\HH{{\mathcal H}}

\newcommand{\G}{{\Gamma}}

\newcommand\LL{{\mathcal L}}
\newcommand\MM{{\mathcal M}}

\newcommand\PP{{\mathcal P}}

\newcommand\RR{{\mathcal R}}
\newcommand\SSS{{\mathcal S}}
\newcommand\TT{{\mathcal T}}

\newcommand\VV{{\mathcal V}}

\newcommand\XX{{\mathcal X}}

\newcommand\TC{{\TT\CC}}

\newcommand\PMF{{\PP\kern-2pt\MM\FF}}
\newcommand\ML{{\MM\LL}}

\newcommand\PML{{\PP\kern-2pt\MM\LL}}

\newcommand\ep{\epsilon}
\newcommand\hhat{\widehat}

 \newcommand\Hyp{{\mathbb H}}

\newcommand\Z{{\mathbb Z}}
\newcommand\R{{\mathbb R}}

\renewcommand\P{{\mathbb P}}

\newcommand{\fsubd}{\mathrel{{\scriptstyle\searrow}\kern-1ex^d\kern0.5ex}}
\newcommand{\bsubd}{\mathrel{{\scriptstyle\swarrow}\kern-1.6ex^d\kern0.8ex}}
\newcommand{\fsubeq}{\mathrel{\raise-.7ex\hbox{$\overset{\searrow}{=}$}}}
\newcommand{\bsubeq}{\mathrel{\raise-.7ex\hbox{$\overset{\swarrow}{=}$}}}

\newcommand{\bbar}{\overline}

\newcommand{\tsh}[1]{\left\{\kern-.9ex\left\{#1\right\}\kern-.9ex\right\}}

\newcommand{\PSL}{{PSL_2 (\mathbb{C})}}
\newcommand{\pslc}{{PSL_2 (\mathbb{C})}}
\newcommand{\pslr}{{PSL_2 (\mathbb{R})}}

\newcommand{\mate}{\bot \!\! \! {\bot}} 

\makeatletter
\DeclareFontFamily{U}{tipa}{}
\DeclareFontShape{U}{tipa}{m}{n}{<->tipa10}{}
\newcommand{\arc@char}{{\usefont{U}{tipa}{m}{n}\symbol{62}}}%

\newcommand{\arc}[1]{\mathpalette\arc@arc{#1}}

\newcommand{\arc@arc}[2]{%
  \sbox0{$\m@th#1#2$}%
  \vbox{
    \hbox{\resizebox{\wd0}{\height}{\arc@char}}
    \nointerlineskip
    \box0
  }%
}
\makeatother

\makeatletter
\providecommand{\bigsqcap}{%
  \mathop{%
    \mathpalette\@updown\bigsqcup
  }%
}
\newcommand*{\@updown}[2]{%
  \rotatebox[origin=c]{180}{$\m@th#1#2$}%
}
\makeatother

\newcommand{\leftrarrows}{\mathrel{\raise.75ex\hbox{\oalign{%
  $\scriptstyle\leftarrow$\cr
  \vrule width0pt height.5ex$\hfil\scriptstyle\relbar$\cr}}}}
\newcommand{\lrightarrows}{\mathrel{\raise.75ex\hbox{\oalign{%
  $\scriptstyle\relbar$\hfil\cr
  $\scriptstyle\vrule width0pt height.5ex\smash\rightarrow$\cr}}}}
\newcommand{\Rrelbar}{\mathrel{\raise.75ex\hbox{\oalign{%
  $\scriptstyle\relbar$\cr
  \vrule width0pt height.5ex$\scriptstyle\relbar$}}}}

\makeatletter
\def\leftrightarrowsfill@{\arrowfill@\leftrarrows\Rrelbar\lrightarrows}
\newcommand{\xleftrightarrows}[2][]{\ext@arrow 3399\leftrightarrowsfill@{#1}{#2}}
\makeatother

\begin{document}
	
\title{Combination Theorems in Groups, Geometry and Dynamics}

\author{Mahan Mj}
\address{School of Mathematics, Tata Institute of Fundamental Research, Mumbai-40005, India}
\email{mahan@math.tifr.res.in, mahan.mj@gmail.com}

\author{Sabyasachi Mukherjee}
\address{School of Mathematics, Tata Institute of Fundamental Research, Mumbai-40005, India}
\email{sabya@math.tifr.res.in, mukherjee.sabya86@gmail.com}
	%
	%

\begin{abstract}
The aim of this chapter is to give a survey  of combination theorems occurring  in hyperbolic geometry, geometric group theory and complex dynamics, with a particular focus on Thurston's contribution and  influence in the field.
\end{abstract}


\date{\today}

\thanks{Both authors were  supported by  the Department of Atomic Energy, Government of India, under project no.12-R\&D-TFR-5.01-0500 as also  by an endowment of the Infosys Foundation.
	MM was also supported in part by   a DST JC Bose Fellowship, Matrics research project grant  MTR/2017/000005, CEFIPRA  project No. 5801-1. SM was supported in part by SERB research project grant  SRG/2020/000018.}

 \maketitle

\tableofcontents

\section{Introduction}\label{history} The aim of this survey is to give an eclectic account of combination theorems in hyperbolic geometry, geometric group theory and complex dynamics. Thurston's contribution and  influence
in the theme is pervasive, and we will only be able to touch upon some of these aspects.  The hope in writing this survey is therefore only to whet the appetite of the reader and provide some references to more detailed articles
and surveys.

Combination theorems have a long history, going back to Klein's paper from 1883 \cite{klein}. A major subsequent development in terms of combination theorems for Kleinian groups is due to Maskit \cite{maskit1,maskit2,maskit3,maskit4,maskit-book}. See Section \ref{klein-maskit} for further details.

A combination theorem of a different, more complex analytic,  flavor was introduced
by Ahlfors and Bers (see Section \ref{sec-su}). The Bers Simultaneous Uniformization Theorem provided a way of combining two abstractly isomorphic Fuchsian surface groups into a single Kleinian surface group. Equivalently,
two discrete faithful representations $\rho_i : \pi_1(S) \to \mathrm{PSL}(2, \R)$ are combined in a dynamically natural way into a single discrete faithful representation $\rho : \pi_1(S) \to \mathrm{PSL}(2, \C)$. The study of Kleinian groups around this time thus took on a rather complex analytic orientation.

 In the 1970's and 80's,
a phase transition occurred in the theory 
with the advent of Thurston, who combined the above two strands into one unified theme, and vastly generalized both. He introduced a  3-dimensional hyperbolic geometry point of view, leading to his proof of hyperbolization of
atoroidal Haken manifolds \cite{thurston-hypstr1,thurston-hypstr2,thurston-hypstr3,otal-book,otal-survey,kapovich-book}. Particular mention must be made of his Double Limit Theorem that may be thought of as a limiting case of the Simultaneous Uniformization Theorem. (See Section \ref{thurston}.)

Thurston's work has had a deep and profound influence 
on hyperbolic geometry ever since, and has provided a template for related
developments in geometric group theory and complex dynamics. In geometric group theory, Bestvina and Feighn \cite{BF} isolated the coarse geometric features of Thurston's combination theorem and proved a highly influential
combination theorem for Gromov-hyperbolic groups \cite{gromov-hypgps},
spawning considerable activity and several generalizations \cite{dahmani-combnthm,alib,mahan-reeves,mahan-sardar,gautero}. In particular, the  main theorem of \cite{BF} was extended to a relatively hyperbolic setup \cite{dahmani-combnthm,alib,mahan-reeves,gautero} and also  to the setup of a coarse-geometric analog of bundles \cite{mahan-sardar}. (See Section \ref{sec-combnthmhyp} for further details.)

In a relatively recent major development leading to a resolution of Thurston's virtual Haken conjecture by Agol and Wise \cite{agol-vhak,wise}, Haglund and Wise \cite{hw-annals} proved
a combination theorem for virtually special cubulable hyperbolic groups \cite{hw-gafa}. For these and related developments, see Section \ref{sec-cube}.

In a closely related theme, Fatou and Julia laid the foundation of the theory of dynamics of rational maps on the Riemann sphere in the first quarter of the twentieth century \cite{fatou-1919,fatou-1920a,fatou-1920b,fatou-1926,julia-1918,julia-1922}. These early developments in the field drove Fatou to observe similarities between limit sets of Kleinian groups and Julia sets of rational maps:
    `L'analogie remarqu{\'e}e entre les ensembles de points limites des groupes Kleineens et ceux qui sont constitu{\'e}s par les fronti{\`e}res des r{\'e}gions de convergence des it{\'e}r{\'e}es d'une fonction rationnelle ne para{\^i}t d'ailleurs pas fortuite et il serait probablement possible d'en faire la synth{\`e}se dans une th{\'e}orie g{\'e}n{\'e}rale des groupes discontinus des substitutions alg{\'e}briques.'  After several decades, this analogy was set on a firm footing by Sullivan with the introduction of quasiconformal techniques in the study of rational dynamics \cite{sullivan-dict}. Shortly afterwards, Sullivan put forward a dictionary between the aforementioned classes of conformal dynamical systems, Thurston proved a topological characterization for an important class of rational maps \cite{DH1} as a philosophical analog of the hyperbolization of atoroidal Haken $3$-manifolds. The theory of polynomial mating, designed by Douady and Hubbard \cite{douady-mating}, extends the notion of a combination theorem from the world of Kleinian groups to that of complex dynamics. This theme too bears the tell-tale stamp of Thurston. In fact, Thurston's topological characterization of rational maps is an invaluable tool in constructing such matings \cite{DH1,shishikura-mating} (see Section \ref{sec-cxdyn}).

The idea of combining Kleinian groups with rational maps was first conceived by Bullett and Penrose in \cite{bullett-penrose}, where they used iterated algebraic correspondences to `mate' the modular group $\mathrm{PSL}(2, \Z)$ with certain quadratic polynomials. More recently, a one complex variable approach was adopted to bind together the actions of Kleinian groups and rational maps in the dynamics of a single map. This perspective can be thought of as a unification of Bers simultaneous uniformization theorem (and in certain cases, Thurston's double limit theorem) with the Douady--Hubbard theory of polynomial mating. The crucial difference between this mating framework and that of Bullett-Penrose is that here one extracts a non-invertible map from a Kleinian group that is \emph{orbit equivalent} to the group on its limit set (i.e., one extracts a semi-group dynamics from the dynamics of a non-commutative group), and `mates' this map with the dynamics of a polynomial. In the anti-holomorphic setting, this is achieved by associating a piecewise circular reflection map, called the \emph{Nielsen map}, to a Kleinian reflection group. The simplest example of this mating phenomenon is given by the \emph{Schwarz reflection map} associated with a simply connected \emph{quadrature domain}; namely, the exterior of a deltoid curve (which is the conformal mating of the anti-polynomial $\overline{z}^2$ and the ideal triangle reflection group). A series of papers 
\cite{LLMM1,LLMM2,LMM2,LMMN} culminated in a comprehensive framework for conformally mating Kleinian reflection groups with anti-holomorphic polynomials (see Section \ref{antiholo_mating_subsec}).

On the holomorphic side, a framework for combining Kleinian groups with polynomial maps was devised in \cite{mj-muk}. The key player in this setting of combination theorems is a class of \emph{piecewise M{\"o}bius} maps, termed \emph{mateable} maps. Such maps are  dynamically and combinatorially compatible with Kleinian groups on the one hand and polynomials on the other. In particular, a mateable map associated to a Kleinian group is orbit equivalent to the group on the limit set. While the simplest examples of mateable maps are given by the classical \emph{Bowen--Series maps} associated with  Fuchsian punctured sphere groups, a new class of examples called \emph{higher Bowen--Series maps} was also described in \cite{mj-muk}. These maps enjoy various close connections with Bowen--Series maps, and are interesting in their own right (for instance, they are responsible for failure of \emph{topological orbit equivalence rigidity} of Fuchsian groups). It turns out that any mateable map can be conformally mated with suitable complex polynomials giving rise to disconnected moduli spaces of matings of punctured spheres with complex polynomials (see Section \ref{holo_mating_subsec}).
\medskip

\noindent \textbf{Acknowledgements.} The authors would like to thank Mikhail Lyubich for helpful conversations and useful suggestions. They also thank Athanase Papadopoulos for a careful reading of the manuscript and suggestions.
 
\section{Klein--Maskit combination  for Kleinian groups}\label{klein-maskit}

A discrete subgroup $\Gamma$ of $\PSL$ is called a \emph{Kleinian group}.
The \emph{limit set} of the Kleinian group $\Gamma$, denoted by
$\Lambda_\Gamma$, is the collection of accumulation points of a $\Gamma$-orbit $\Gamma\cdot z$ for some $z \in \hat{\mathbb{C}}$. $\Lambda_\Gamma$ is independent of $z$. It may be thought of as the locus of chaotic dynamics of $\Gamma$ on $\hat{\mathbb{C}}$, i.e.\ for $\Gamma$ non-elementary and any $z \in \Lambda_\Gamma$, $\Gamma \cdot z$ is dense  in $\Lambda_\Gamma$. We shall identify
the Riemann sphere $\hat{\mathbb{C}}$  with the sphere at infinity $\mathbb{S}^2$ of $\Hyp^3$. The complement of the limit set $\hat{\mathbb{C}} \setminus \Lambda_\Gamma$ is called the domain of discontinuity $\Omega(\Gamma)$ of $\Gamma$. If the Kleinian group $\Gamma$ is torsion-free, it acts freely and properly discontinuously on $\Omega(\Gamma)$ with a Riemann surface quotient. 

\begin{definition}
A set $D$ is called a \emph{partial fundamental domain}  for $\G$, if
\begin{enumerate}
\item  $D\neq \emptyset$,
\item $ D \subset \Omega(\G)$, and
\item $g(D) \cap D = \emptyset$, for all $g\in \G$, $g\neq  1$.
\end{enumerate}

If, further, $\bigcup_{g \in \G} g. D = \Omega (\G)$,
then $D$ is called a \emph{fundamental domain}  for $\G$.
\end{definition}

The story of combination theorems starts with the following theorem of Klein:

\begin{theorem}[Klein combination theorem] \cite{klein}\label{klein}
	Let $\Gamma_1, \Gamma_2$ be Kleinian groups with fundamental domains
	$D_1, D_2$
	respectively.
	Assume that the interior of $D_1$ (resp.\ $D_2$)
	contains the boundary and exterior of $D_2$ (resp.\ $D_1$).Then the group $ \G$ generated by $\Gamma_1, \Gamma_2$ is
	Kleinian,
 and $D=D_1 \cap D_2$ is a fundamental domain for $\G$.
\end{theorem}

In the 1960's, Maskit started working on extending the Klein combination Theorem \ref{klein}
to a more general setup.
Maskit's work on combination theorems for Kleinian groups started with the following.
\begin{theorem}[Klein-Maskit combination theorem for free product with amalgamation] \cite{maskit0}\label{maskit1-amalg}
Let $\Gamma_1, \Gamma_2$ be Kleinian groups with domains of discontinuity $\Omega_1, \Omega_2$ respectively. Let $H=\Gamma_1\cap \Gamma_2$.  Let $D_1, D_2, \Delta$
 be partial fundamental domains for $\G_1, \G_2, H$ respectively. For
i = 1,2, set $E_i = H.D_i$. Denote the interior of $D=E_1 \cap E_2 \cap \Delta$ by $D'$. If
\begin{enumerate}
\item $D'\neq \emptyset$,
\item $E_1 \cup E_2 = \Omega_1 \cup \Omega_2$.
\end{enumerate}
 Then the group $ \G$ generated by $\Gamma_1, \Gamma_2$ is
Kleinian, $D'$ is a partial fundamental domain for $\G$, and $\G=\G_1 \ast_H \G_2$ is the free product with amalgamation of $\G_1, \G_2$ along $H$.  Further, $gD \cap D = \emptyset$, 
for all $g\in \G$, $g\neq  1$.
\end{theorem}

 In \cite{maskit1}, Maskit strengthened the above theorem by determining precisely a fundamental domain for the group. 
\begin{theorem}[Klein-Maskit combination theorem for free product with amalgamation] \cite{maskit1}\label{maskit1-amalg2}
Let $\Gamma_1, \Gamma_2$ be Kleinian groups with domains of discontinuity $\Omega_1, \Omega_2$ respectively. Let $H \subset \Gamma_1\cap \Gamma_2$.
such that $H$ is either cyclic or consists only of the identity. 
Let $D_1, D_2, \Delta$
be  fundamental domains for $\G_1, \G_2, H$ respectively. 
For
$i = 1,2,$ set $E_i = H.D_i$. Denote the interior of $D=E_1 \cap E_2 \cap \Delta$ by $D'$.
 Suppose  $E_1 \cup E_2 = \Omega(H)$ and 
 that $D'\neq \emptyset$. Assume further that
there is a simple closed curve $\gamma$ contained in $int(E_1 \cup E_2) \cup \Lambda_H$
such that $\gamma$  is invariant
under $H$, the closure of $\gamma \cap \Delta$  is contained in 
$int(E_1 \cap E_2)$ and $\gamma$ separates
both $E_1 \setminus E_2$ and $E_2 \setminus E_1$.  Then the group $ \G$ generated by $\Gamma_1, \Gamma_2$ is
Kleinian, $\G=\G_1 \ast_H \G_2$, and
$D$ is a fundamental domain for $G$. 
\end{theorem}

Subsequently, in \cite{maskit2,maskit3}, Maskit upgraded Theorem \ref{maskit1-amalg2} to the following.
We start with two Kleinian groups $ \G_1, \G_2$
 with $H\subset \G_1 \cap \G_2$, where $H\neq \G_1, \G_2$. We are also given a simple closed curve $\gamma$ dividing the Riemann sphere $\hat \C$ into two closed topological discs, $B_1$ and $B_2$, where $B_i$ is precisely
invariant under $H$ in $\G_i$.
 More precisely, $B_i$ is $H-$invariant, and if $g \in  \G_i \setminus H$, then
$g(B_i)\cap B_i = \emptyset$. Then $\G = \langle \G_1, \G_2\rangle$, the group generated by $\G_1$ and $\G_2$, is
also a Kleinian group. What really needs to be proved in all these cases is the discreteness of $\G$.

In all these cases,
Maskit shows  that $\G  = \G_1\ast_H \G_2$, i.e.\ $\G = \langle \G_1, \G_2\rangle$
is equal to the free product with amalgamation of $\G_1, \G_2$ along $H$.
Further, by carefully choosing fundamental domains for $\G_1,
\G_2$ one can ensure that their
intersection will be a fundamental domain for $\G$.
Thus, the basic hypothesis guaranteeing discreteness of $\G$ can be summarized as follows: 
\begin{enumerate}
\item The disks $B_1$ and $B_2$ are both invariant under $H$.
\item The
$(\G_1 \setminus H)-$translates of $B_1$ are disjoint disks in $B_2$.
\item The
$(\G_2 \setminus H)-$translates of $B_2$ are disjoint disks in $B_1$.
\end{enumerate}

 There is also a version of the Klein-Maskit combination theorem for HNN extensions.
 We are given a single group $\G_0$, 
 with two subgroups $H_1$ and $H_2$,
two closed disks $B_1$ and $B_2$ which have disjoint projections to $\Omega(\G_0)/\G_0$,
where 
\begin{enumerate}
\item $H_i$ preserves $B_i$,
\item there exists  a M\"obius transformation $h$ 
mapping
the outside of $B_1$ onto the inside of $B_2$ and conjugating
$H_1$ to $H_2$.
\end{enumerate}
Maskit then shows that $ \G = \langle \G_0, h \rangle$ is a Kleinian group. Also $ \G =\G_0 \ast_H$ is the HNN-extension
of $G_0$ along $H$, where the two inclusions of $H$ map to $H_1, H_2$ and $h$ is the stable letter conjugating one to the other.
Further, by carefully choosing a fundamental domain $D$  for
$\G_0$, one can ensure that $D \setminus (B_1 \cup B_2)$ is a fundamental domain for $\G$. 

Maskit weakens the hypotheses further in \cite{maskit4}, allowing translates of the closed disks $B_1, B_2$ to have common boundary points. However, in \cite{maskit4}
he requires that such points of intersection also
be ordinary points of our original group. 
In \cite{maskit4}, it is also shown that $\G$ is geometrically finite if and only if the original groups are so.
The basic topological tool used in the proof is a Jordan curve $\gamma$ in $\hat{\C}$ and its translates under a Kleinian
group. The standard hypothesis in these papers is
the `almost disjointness' of $\gamma$
from all its translates.
More precisely, if $g(\gamma) \cap \gamma \neq \emptyset$,
then it is required that $g(\gamma)$ is entirely 
contained in the closure of one of the open disks bounded by $\gamma$. Thus, a substantial amount of the technical difficulty in \cite{maskit3,maskit4} comes from controlling the points of intersection $g(\gamma) \cap \gamma$.

To conclude this section, we refer the reader to 
\begin{enumerate}
\item Work of Li-Ohshika-Wang \cite{low1,low2} for generalizations of the Klein-Maskit combination theorems to higher dimensions.
\item Work of Dey, Kapovich and Leeb \cite{dekl} for a combination theorem for Anosov subgroups, a natural class of discrete 
subgroups of higher rank Lie groups that generalizes convex cocompact subgroups of $\pslc$.
\end{enumerate}

\section{Simultaneous uniformization and Quasi-Fuchsian Groups}\label{sec-su} The aim of this section is to give a brief account of Bers simultaneous uniformization theorem.
The reason is 2-fold. First, it provides the context for Thurston's double limit theorem in Section \ref{sec-dlt}.
Second, it is the Kleinian group analog for the Douady--Hubbard mating construction \cite{douady-mating}, and more generally the original motivation for the mating constructions in Section \ref{interbreeding_sec}.

Fix a surface $S$.
The collection of all representations $\rho : \pi_1(S) \to \pslc$ up to conjugacy (in $\pslc$) is called the \emph{character variety} and is represented as $\RR(S)$.
We note in passing that the appropriate quotient by $\pslc$ of the space of  all representations $\rho : \pi_1(S) \to \pslc$ is the GIT quotient. This is needed in order to obtain the structure of a variety on $\RR(S)$.

\subsection{Topologies on space of representations}
\label{sec-tops} 

For future reference, we summarize here a natural collection of topologies on the space of \emph{discrete faithful} $\rho : \pi_1(S) \to \pslc$.
The \emph{algebraic topology} is the topology of pointwise convergence on elements of   $\pi_1(S)$:

\begin{definition}\label{def-altop}
	We shall say that	a sequence of representations $\rho_n : \pi_1(S) \to \pslc$  converges \emph{algebraically} to $\rho_\infty : \pi_1(S) \to \pslc$  if for all $g \in \pi_1(S)$, $\rho_n (g) \to \rho_\infty(g)$ in $\pslc$. 
\end{definition}

The collection of conjugacy classes of discrete faithful representations of  $\pi_1(S)$ into $\pslc$ equipped with the algebraic topology is denoted as $AH(S)$. Thus, $AH(S)\subset \RR(S)$ comes naturally equipped with a complex analytic structure. The space of discrete faithful representations of  $\pi_1(S)$ into $\pslr$ equipped with the algebraic topology is precisely the Teichm\"uller space. Thus, the Teichm\"uller space
sits `diagonally' in $AH(S)$.

For analyzing convergence from a geometric point of view, the natural topology is the \emph{geometric topology}, or equivalently, the Gromov--Hausdorff topology.

\begin{definition}\label{def-geotop}
	Let $\rho_n: \Gamma \rightarrow \pslc$ be a sequence of discrete, faithful
	representations of a finitely generated, torsion-free, nonabelian group $\Gamma$. Thus, $\rho_n(\Gamma)$ is a sequence of closed subsets of $\pslc$.
	If $G \subset \pslc$ is a closed subgroup such that
	$\rho_n(\Gamma)$ converges to $G$ in the Gromov--Hausdorff
	topology, then
	$\rho_n(\Gamma)$ is said to \emph{converge geometrically} to $G$, and  $G$ is called the 
	\emph{geometric limit}.
\end{definition}

\begin{definition}\label{def-strong}
	$G_n$ converges \emph{strongly}
	to  $G$ if $G_n$ converges to
	$G$  both algebraically and geometrically.
\end{definition}

\subsection{Simultaneous uniformization}\label{subsec-su}
\begin{definition}
	Let $\rho: \pi_1(S) \to \pslc$ be a discrete faithful representation such that the limit set of $G= \rho(\pi_1(S))$ is a topological circle in $\mathbb{S}^2$. Then $G$ is said to be \emph{ quasi-Fuchsian}. The collection of conjugacy classes of quasi-Fuchsian representations is denoted as $QF(S)$. 
\end{definition}

Note that $QF(S)$ is contained in $AH(S)$ and hence inherits 
a complex analytic structure.
The domain of discontinuity $\Omega$ of a  quasi-Fuchsian $G$ consists of two open invariant disks $\Omega_1, \Omega_2$. Hence the quotient $\Omega/G$ is the disjoint union $\Omega_1/G \sqcup \Omega_2/G$. Hence we have a map $\tau: QF(S) \to Teich(S) \times Teich(S)$, where $Teich(S)$ denotes the Teichm\"uller space of $S$.
The { \bf Bers simultaneous Uniformization Theorem} asserts:

\begin{theorem}\cite{bers-su,bers-cpam}
	$\tau: QF(S) \to Teich(S) \times Teich(S)$ is a homeomorphism.
	\label{thm-su}
\end{theorem}

\begin{figure}
\centering
\includegraphics[width=0.6\linewidth]{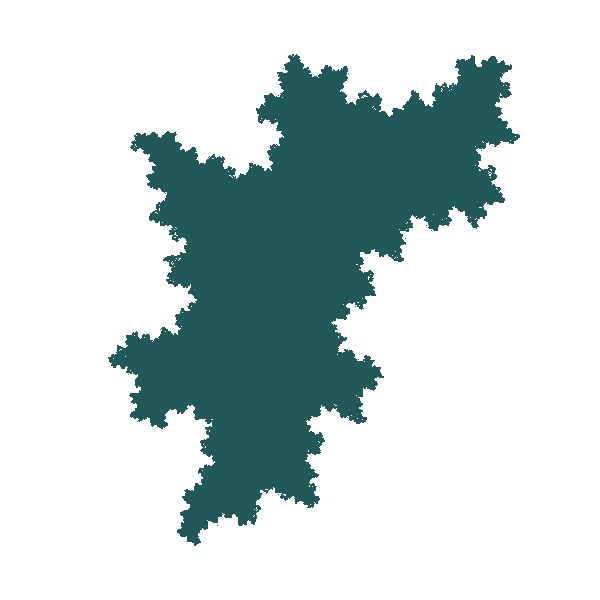}
\caption{Quasi-Fuchsian group limit set}
\label{qf}
\end{figure}

Hence, given any two conformal structures $T_1, T_2$ on a surface, there is a unique discrete quasi-Fuchsian  $G$  whose limit set $\Lambda_G$ is  topologically a circle, and the quotient of whose domain
of discontinuity is $T_1 \sqcup T_2$. See Figure \ref{qf} below \cite{kabaya}, where the inside and the outside of the Jordan curve correspond to $\Omega_1, \Omega_2$.

We refer to \cite{hubbard-book1} for a proof of Theorem \ref{thm-su} and summarize the main ideas here. 
Theorem \ref{thm-su} is essentially complex analytic in nature
and goes back to an understanding of the Beltrami partial differential equation due to Morrey. Let $K_S$ denote the canonical bundle of the Riemann surface $S$ (if $S$ has punctures as a hyperbolic surface, we regard them as marked points in the complex analytic category).

\begin{definition}\label{def-bd}
	A \emph{Beltrami differential} on $S$ is an $L^\infty$ section of 
	$K_S^{-1} \otimes \bbar{K_S}$, where $\bbar{K_S}$ denotes the
	complex conjugate of $K_S$. The space of Beltrami differentials on $S$ will be denoted as $\DD_b(S)$
\end{definition}

The local expression for an element of $\DD_b(S)$ in a complex
analytic chart $U\subset S$ is thus given by
$\mu \frac{d\bbar{z}}{dz}$, where $\mu \in L^\infty(U)$ is called a
\emph{Beltrami coefficient}. 

\begin{definition}\label{def-qc}
	A \emph{quasiconformal map} between  two  Riemann  surfaces $S_1$ and $S_2$
	is  a  homeomorphism $\phi:S_1 \to S_2$
	having locally square-integrable weak partial derivatives  such that $$\mu =\frac{\phi_{\bbar{z}}}{\phi_{{z}}}$$ satisfies 
	$||\mu||_\infty  < 1$. Here, $\mu$ is called the
	\emph{ Beltrami coefficient} of $\phi$.
\end{definition}

The first major ingredient
in the proof of Theorem \ref{thm-su} is the Measurable Riemann
mapping theorem. As pointed out by Hubbard in \cite[p. 149]{hubbard-book1}, the Beltrami coefficient $\mu$ really represents an almost-complex structure on $U$
and the Measurable Riemann
mapping theorem (due to Ahlfors--Bers--Morrey) below     ensures its integrability to
a complex structure.

\begin{theorem}[Measurable Riemann Mapping Theorem]\cite[Theorem 4.6.1]{hubbard-book1}\label{thm-rmt}\\
	\noindent {\bf Existence of quasiconformal maps:}	Let $U \subset \C$ be open. Let $\mu  \in  L^\infty (U)$ satisfying $||\mu||_\infty  < 1$ . Then
	there exists a quasiconformal mapping $f : U\to \C$ solving the
	Beltrami equation $$\frac{\partial f}{\partial \bbar{z}}
	= \mu \frac{\partial f}{\partial {z}}.$$
	\noindent {\bf Uniqueness of quasiconformal maps:}
	If $g$ is another quasiconformal solution to the
	Beltrami equation above, then
	there exists a univalent analytic function $\phi : f (U)\to \C$ such that
	$g = \phi\circ f$.
\end{theorem}

The rest of this brief account of Theorem \ref{thm-su} follows
\cite{gupta-teich} which captures the relevant
conformal geometry.
Recall that we have fixed  a base Riemann surface
$S$.  Let $\G < \pslr$ be a (base)
Fuchsian group uniformizing $S$.
Let $S'$ be an arbitrary point in the Teichm\"uller 
space of the underlying topological surface.
Theorem \ref{thm-su} then associates to $S'$ a quasiconformal map
$\Phi: \hat{\C} \to \hat{\C}$ fixing
the three points $0,1,\infty$,  and conjugating the action of $\G$ to that of a
Kleinian group $\G(S,S')$, such that
\begin{enumerate}
\item $\Phi$ is conformal on the lower half-plane.
\item $\G(S,S')$ leaves  invariant  the  images  of the lower
and upper  half-planes.
\item The quotient of the lower half-plane by $\G(S,S')$ is
$S$.
\item The quotient of the upper half-plane by $\G(S,S')$ is
$S'$.
\end{enumerate}

To prove the existence of a $\Phi$ as above, we first note that the Teichm\"uller space
can  be identified (via Theorem \ref{thm-rmt}) with Beltrami differentials on $S$ with norm bounded by one.
Let $\mu$ be the Beltrami differential on $S$ corresponding to $S'$.
Next. lift
 $\mu$ to the upper half plane.
 Extend to a Beltrami coefficient $\mu_0$ on $\hat{\C}$ by defining it to be identically zero on
 the lower half plane.  The map
 $\Phi$ above is then given by
 the normalized solution to the Beltrami equation with the Beltrami coefficient $\mu_0$.  Invariance of $\mu_0$ under 
 $\Gamma$ ensures the existence of an isomorphism
 $\rho: \G \to \G(S,S')$ to the desired Kleinian group $\G(S,S')$ such that
 $$\Phi \circ g = \rho(g)\circ \Phi, $$ for all $g \in \G$.

\subsection{Geodesic Laminations} We turn now to the hyperbolic geometry 
of quasi-Fuchsian groups.

\begin{definition} A geodesic lamination on a hyperbolic surface is a foliation of a closed subset with geodesics. 
\end{definition} 

A geodesic lamination on a surface may further be equipped with a transverse (positive) measure to obtain a \emph{measured lamination}. The space $\ML (S)$ of measured (geodesic) laminations on $S$  then has the structure of a positive cone in a vector space, i.e.\ for every $\lambda \in \ML (S)$ and $t \in \R_+$, $t\lambda \in \ML (S)$. It can be projectivized to obtain the space of projectivized measured  laminations $\PML(S)$. It was shown by Thurston \cite{FLP} that

\begin{theorem}\label{thm-pmlcptfs}
$\PML(S)$ is homeomorphic to a sphere and can be adjoined to $Teich(S)$ compactifying the latter to a closed ball.
\end{theorem} 

\begin{definition} \cite{thurstonnotes}[Definition 8.8.1] A \emph{pleated surface} in a hyperbolic three-manifold $N$ is
	a complete hyperbolic surface $S$ of finite area, together with an isometric map $f :
	S \to N$ such that every $x \in S$ is in the interior of some geodesic segment which
	is mapped by $f$ to a straight line segment. Also, $f$ must take every cusp of $S$ to a
	cusp of $N$
	
	The pleating locus of the pleated surface $f: S \to M$  is the set $\gamma \subset S$ consisting of those points in the pleated surface which are in the
	interior of unique line segments mapped to line segments.
\end{definition}

\begin{prop} \cite{thurstonnotes}[Proposition 8.8.2] The pleating locus $\gamma$
	is a geodesic lamination on $S$. The map $f$ is totally geodesic in
	the complement of 
	$\gamma$. \end{prop}

The geometry of quasi-Fuchsian groups and their relationship with geodesic laminations arises out of the geometry of the convex core that we now describe.

\begin{definition}\label{def-ch}
Let $\G$ be an infinite Kleinian group and let $\Lambda
\subset \hat{\C}$ denote its limit set. The \emph{convex hull} of $\Lambda$ is the smallest non-empty closed convex subset of $\Hyp^3$ whose set of accumulation points in 
$\hat{\C}$ equals $\Lambda$. We denote the convex hull of $\Lambda$ by $CH(\Lambda)$.

The convex hull $CH(\Lambda)$ of a Kleinian group $\G$ is invariant under $\G$. The quotient $CH(\Lambda)/\G \subset \Hyp^3/\G$ is called the \emph{convex core} of the hyperbolic
3-manifold $M=\Hyp^3/\G$.
\end{definition}
For a quasi-Fuchsian group $\G=\rho(\pi_1(S))$, the convex
core is homeomorphic to a product $S \times [a,b]$ (the Fuchsian case corresponds to $a=b$). The hyperbolic distance
between $S \times \{a\}$ and $S \times \{b\}$ is a measure 
of the geometric complexity of $\G$.
In \cite{thurstonnotes}[Ch. 8],
Thurston further shows:

\begin{prop}\label{prop-ps}
Let $M$ be a complete hyperbolic 3-manifold corresponding to a quasi-Fuchsian group $\G$, and let $CC(M)$ denote its convex core. Then each component of the convex core boundary $\partial CC(M)$ is a  pleated surface. 
\end{prop}

\section{Thurston's combination theorem for Haken manifolds}\label{thurston} The material in this section provides the core inspiration for most combination theorems that came subsequently.

\begin{definition}\cite{hempel-book}
	A properly embedded surface $(F, \partial F) \subset 
	(M, \partial M)$ in a 3-manifold $M$ with boundary $\partial M$ (possibly empty) is said to be \emph{incompressible} if 
the inclusion map $i: (F, \partial F) \subset 
	(M, \partial M)$ induces an injective homomorphism of fundamental groups $i_\ast: \pi_1(F) \to \pi_1(M)$.
	Further, we require that for every boundary component
	$\gamma$ of $\partial F$, $i_\ast: \pi_1(\gamma) \to \pi_1(\partial M)$ is injective. (The second condition 
	is automatic when $F$ is not a disk.)
	
	An embedded incompressible  surface $(F, \partial F) \subset 
	(M, \partial M)$ is said to be \emph{boundary parallel}
	if $F$ can be isotoped into $\partial M$ keeping
	$\partial F \subset \partial M$ fixed.

A compact 3-manifold $M$ (possibly with boundary $\partial M$) is said to be \emph{Haken} if 
\begin{enumerate}
\item $\pi_2(M)=0$.
\item There exists an embedded incompressible  surface 
$(F, \partial F)$ that is not boundary parallel.
\end{enumerate}

$M$ is said to be atoroidal if $\pi_1(M)$ contains no $\Z\oplus \Z$  subgroups.  $M$ is said to be acylindrical if
any embedded incompressible annulus in $	(M, \partial M)$
is boundary parallel.
\end{definition}

We summarize Thurston's celebrated
hyperbolization theorem now and then give a brief account of the ingredients that go into the proof.

\begin{theorem}\label{thm-monster}\cite{thurston-hypstr1,thurston-hypstr2,thurston-hypstr3}
Let $M$ be a compact atoroidal Haken 3-manifold. Then $M$ is hyperbolic.
\end{theorem}

There is a version of Theorem \ref{thm-monster} for 3-manifolds with torus boundary components also. But, in the interests of exposition, we shall largely focus on the compact atoroidal case. The proof of Theorem \ref{thm-monster} breaks into two principal pieces:

\begin{enumerate}
\item $M$ is compact atoroidal Haken and does not fiber over
$\mathbb{S}^1$. This case will be described in Section \ref{sec-monster}.
\item $M$ fibers over the circle with fiber $F$.
This case will be described in Section \ref{sec-dlt}.
\end{enumerate}

\subsection{Non-fibered Haken 3-manifolds}\label{sec-monster}
There are a number of detailed expositions for the 
compact atoroidal Haken non-fibered case and we point out 
\cite{otal-survey,kapovich-book,ctm-itn} in particular.

It is a fundamental fact of 3-manifold topology \cite[Chapter 13]{hempel-book} that any Haken manifold admits a Haken hierarchy. Cutting $(M,\partial M)$ open along 
$(F,\partial F)$ gives us a new (possibly disconnected) atoroidal 3-manifold with boundary. The cut open manifold
is automatically Haken, and we can proceed inductively.
At the last stage, we are left with a finite collection of balls, and these are clearly hyperbolic.

Thus, in order to prove Theorem \ref{thm-monster} in the non-fibered case, an essential step is the following:

\begin{theorem}\cite{otal-survey}\label{thm-acyl}
Let $M_1$ be an acylindrical atoroidal
 3-manifold with non-empty incompressible boundary $\partial M_1$ whose interior admits a complete hyperbolic metric. Let $\tau : \partial M_1 \to \partial M_1$ be
an orientation-reversing involution. Then the interior of
$M=M_1/\tau$ admits a complete hyperbolic metric.
\end{theorem}

To prove Theorem \ref{thm-acyl}, a first tool is the following generalization of Theorem \ref{thm-su}:

\begin{theorem}\label{thm-teichbdy}
Let $M_1$ be a compact 3-manifold with boundary such that 
\begin{enumerate}
	\item 
The interior of $M_1$ admits a complete hyperbolic metric.
\item No component of $\partial M_1$ is homeomorphic to a torus or a sphere.
\end{enumerate}
 Then the space of complete hyperbolic metrics on $M_1$ is given by $Teich(\partial M_1)$.
\end{theorem}

Let $\partial M_1 = \sqcup_i \Sigma_i$, where each $\Sigma_i$
is a surface of genus greater than one. Then 
$Teich(\partial M_1) = \prod_i Teich(\Sigma_i)$. Fix a complete hyperbolic structure on $M_1$ (the existence of such a structure is guaranteed by the hypothesis of Theorem \ref{thm-teichbdy}). This is equivalent to a discrete faithful representation $\rho: \pi_1(M_1) \to \pslc$.
Let $\G = \rho(\pi_1(M_1))$. Then
each $\Sigma_i \subset \partial M_1$
gives (via inclusion) 
 a conjugacy class of quasi-Fuchsian subgroups of $\G$.
Thus the involution $\tau$ of Theorem \ref{thm-acyl}
induces a map 
$$\sigma:Teich(\partial M_1) \to Teich(\partial M_1).$$ The map $\sigma$ is called the \emph{skinning map}. The existence
of a complete hyperbolic structure on $M=M_1/\tau$ is equivalent
to the existence of
 a fixed point of the skinning map $\sigma$, as such a fixed point ensures an  isometric gluing. Thurston's
 fixed point
theorem can now be stated as the following:
If
$M$
is atoroidal
then
$\sigma$
has
a
fixed point.
It now follows from the Klein-Maskit combination theorem (Section \ref{klein-maskit})
that $M$
admits a complete hyperbolic structure. 
The acylindricity hypothesis of Theorem
\ref{thm-acyl}
guarantees
that $M$
is atoroidal, completing an outline of the proof of
Theorem \ref{thm-acyl}.

An effective method of proving the existence of
a fixed point of the skinning map $\sigma$ was taken by McMullen in \cite{ctm-itn}. As mentioned in \cite{otal-survey}
Hubbard had observed that the
analytical formula for the coderivative of the
skinning map relates it to  the Theta operator in Teichm\"uller theory.
McMullen studies the fixed point problem via
$||D\sigma||$, the norm of the
derivative  of  the  skinning  map. He reproves Theorem \ref{thm-acyl}
by showing that  if $M_1$
is acylindrical, then there exists $c<1$ such that
$||D\sigma||<c$  guaranteeing  a solution to the gluing problem. 

Both Thurston's fixed-point theorem and McMullen's estimates
in \cite{ctm-itn} are easiest to state when $M_1$ is acylindrical. However, both approaches can be refined to conclude hyperbolicity of $M$ as long as $M_1$ is not of the
form $S \times I$ and $\tau$ glues $S \times \{0\}$ to 
$S \times \{1\}$. The excluded case is that of 3-manifolds 
fibering over the circle and involves a completely different approach that we describe now.

\subsection{The double limit theorem}\label{sec-dlt}
We shall follow \cite{otal-book} to give an outline of the steps involved in the hyperbolization of 3-manifolds fibering over the circle.  Recall (Theorem \ref{thm-pmlcptfs}) that the space of projectivized measured laminations $\PML(S)$
compactifies $Teich(S)$. Thurston's double limit theorem may be thought of as an extension of the simultaneous uniformization Theorem \ref{thm-su} to the case where the pair
$(\tau_1, \tau_2)$ of Riemann surfaces in $Teich(S) \times Teich(S)$ is replaced by a pair $(\ell_1, \ell_2) \in
\bbar{Teich(S)} \times \bbar{Teich(S)}$, where 
$\bbar{Teich(S)} = Teich(S) \cup \PML(S)$ denotes the Thurston
compactification of $Teich(S)$ as in Theorem \ref{thm-pmlcptfs}.

Dual to any measured lamination $\ell \in \ML(S)$ there is an
action of $\pi_1(S)$ on an $\R-$tree. An $\R-$tree is a
geodesic metric space such that any two distinct points are joined by a unique arc isometric to an interval in $\R$.
We refer to \cite{bestvina-rtree} for an expository account
of group actions on $\R-$trees and convergence of $\G-$spaces, 
and mention only the following theorem. Fix a group $\G$.
A triple $(X,o,\rho)$ is called a based $\G-$space if $o\in X$ is a
base-point, and $\G$ acts on $X$ via a homomorphism $\rho:\G \to Isom(X)$ from $\G$ to the isometry group $Isom(X) $ of $X$.

\begin{theorem} \label{cgnce}\cite{bestvina-rtree}[Theorem 3.3]
	Let $(X_i,o_i,\rho_i)$ be a convergent sequence of based $\G$-spaces such that
	\begin{enumerate}
		\item Each $X_i$ is $\delta$ hyperbolic, for some
		 $\delta\geq 0$.
		\item there exists $g\in \G$ such that the sequence
		$d_i=d_{X_i}(o_i,\rho_i(g)(o_i))$ is unbounded.
	\end{enumerate}
	Then there is a based $\R$-tree $(T,o)$ and an isometric action
	$\rho:\G\to Isom(T)$ such that $(X_i,o_i,\rho_i)\to (T,o,\rho)$.
	
	Further, the (pseudo)metric on 
	the $\R-$tree $T$ 
	is obtained as the
	limit of pseudo-metrics $\frac{d_{(X_i,o_i,\rho_i)}}{d_i}$.
\end{theorem}

Finally, we shall need the following theorem of Skora \cite{skora} on the structure of groups admitting small actions on $\R-$trees.

\begin{theorem}\label{skora}\cite{skora}
Let $S$ be a finite area hyperbolic surface. Suppose $\pi_1(S)$ acts non-trivially on an $\R-$tree $T$, such that
for every cusp $\P$ of $S$,
$\pi_1(\P)$ fixes a point in $T$.
 Then the stabilizer of each non-degenerate arc of $T$ contains no free subgroup of rank 2 if and only if the action is dual to an element of $\ML(S)$.
\end{theorem}

An action of $\pi_1(S)$ on an $\R-$tree $T$ such that 
the stabilizer of each non-degenerate arc of $T$ contains no free subgroup of rank 2 is called a \emph{small action}.
Morgan and Shalen \cite{ms1,ms2,ms3} constructed a compactification of the  variety
$\RR(S)$ by the   space of projectivized
length functions arising from small actions of $\pi_1(S)$ on $\R-$trees.
Skora's theorem \ref{skora} allows us to replace $\PML(S)$
in Theorem \ref{thm-pmlcptfs} by  such
length functions.

With this background in place we return to an outline of Thurston's double limit theorem \cite{thurston-hypstr2} following Otal \cite{otal-book}. Let $(\tau_i^+, \tau_i^-)
\in Teich(S) \times Teich(S)$ 
be
a
sequence
of points
converging to  
$(\ell_+, \ell_-) \in
\bbar{Teich(S)} \times \bbar{Teich(S)}$. By Theorem \ref{thm-su}, we can identify $Teich(S) \times Teich(S)$ with
$QF(S)$ and hence assume that $(\tau_i^+, \tau_i^-)
\in QF(S)$. For convenience of exposition, we assume that
$\ell_+, \ell_-$ are both in $\PML(S)$ (a similar statement holds if only one of $\ell_+, \ell_-$ lies in $\PML(S)$).
Assume further that $\ell_+, \ell_-$ \emph{fill} $S$, i.e.\
 each component of
 $S\setminus
 (\ell_+\cup\ell_-)$
 is either  simply connected or else is topologically a punctured disk. Let $\rho_i: \pi(S) \to \pslc$ be the quasi-Fuchsian representation corresponding to $(\tau_i^+, \tau_i^-)
 \in QF(S)$ and let $\G_i = \rho_i(\pi(S))$.
 Thurston's double limit theorem now says:
 
 \begin{theorem}\cite{thurston-hypstr2}\label{thm-dlt}
 Under the above assumptions, there exists a Kleinian group
 $\G$ such that $\G_i$ converges to $\G$ in $AH(S)$.
 \end{theorem}

We sketch Otal's proof following \cite{otal-book} and argue by contradiction. If $\G_i$ diverges, then Theorem \ref{cgnce}
shows
that
there
is
a
limiting small action of
$\pi_1(S)$
on an $\R-$tree $T$.  By Theorem \ref{skora} such a small  action is dual
to a
measured
lamination
$\ell$
on
$S$.

 It is then shown in \cite{otal-book} that any measured lamination that intersects $\ell$ essentially
is realizable in $T$.  Hence at least one of $\ell_+$ and $\ell_-$ must be realizable in $T$, since the two together fill $S$. Without loss of generality, suppose $\ell_+$ is realizable in $T$.
This allows us to approximate $\ell_+$ by simple closed curves
$\sigma$ on $S$ and 
 estimate the translation length
 $l_i(\sigma)$ of $\sigma$ in $\Hyp^3/\G_i$.
 The estimate thus obtained   contradicts a classical estimate of $l_i(\sigma)$ due to Ahlfors obtained in terms of
 the length of the geodesic realization of $\sigma$ in
 $\tau_i^+$ and $\tau_i^-$. This final contradiction proves
 Theorem \ref{thm-dlt}.
 
 Finally, to hyperbolize an atoroidal 3-manifold fibering over the circle with monodromy $\phi$, one picks a base Riemann surface $\tau$, and sets $\tau_i^+ = \phi^i(\tau)$
 and $\tau_i^- = \phi^{-i}(\tau)$. Then $\ell_+, \ell_-$
 turn out to be the stable and unstable laminations of $\phi$.
 The 3-manifold $M$ obtained from the double limit theorem is easily seen to be invariant under $\phi$, and hence $M$ admits a quotient which is the required hyperbolic 3-manifold.

\section{Combination theorems in geometric group theory: hyperbolic groups}\label{sec-combnthmhyp}
The fundamental combination theorem in the context of
hyperbolic groups in the sense of Gromov \cite{gromov-hypgps}
 is due to Bestvina and Feighn \cite{BF}. The theorem was motivated by Thurston's combination Theorem \ref{thm-monster}.
 In the context of geometric group theory, free products with amalgamation and HNN extensions can be treated on a common footing by passing to the universal cover and looking at
 the resulting Bass--Serre tree of spaces \cite{scott-wall}.
 Thus, while the main combination theorem of \cite{BF} provides only a weaker conclusion than Theorem \ref{thm-monster} inasmuch as it establishes Gromov-hyperbolicity, the context is considerably more general and works for trees of spaces.
 It turns out that the sufficient 
 condition in \cite{BF} is also necessary and this converse direction was established
by Gersten \cite{gersten}, Bowditch
\cite{bowditch-ct} and others. The paper \cite{BF} spawned
considerable activity in geometric group theory and have been giving rise
to  a number of combination theorems
\cite{dahmani-combnthm,alib,mahan-reeves,gautero-conv,gautero-weidmann,mahan-sardar,gautero,martin-osajda}  
right up to the time of writing this article. A forthcoming book of Kapovich and Sardar \cite{kap-sar} furnishes a definitive account and rather general versions
of the material in Sections \ref{sec-trees} and \ref{sec-mbdl}.

\subsection{Trees of  spaces}\label{sec-trees}
The framework of \cite{BF} is that of a  tree of spaces. We follow
the exposition in \cite{mahan-tight} to define the relevant notions.
\begin{definition} \label{def-tree}
	\cite{BF} Let $(X,d)$ be a geodesic
	metric space.
	 Let $T$ be a simplicial tree.  Let $\VV(T)$
	  and $\EE(T)$ denote the vertex set and  edge set of $T$
	  respectively. Then $P: X \rightarrow T$ is said to be a  tree of geodesic
	metric spaces satisfying
	the \emph{quasi-isometrically embedded condition} (or simply, the \emph{qi condition}) if  there exists
	a map $P : X \rightarrow T$, and constants
	$K \geq 1, \epsilon \geq 0$ satisfying the
	following: 
	\begin{enumerate}
		\item  $\forall v\in{\VV(T)}$,
		$X_v = P^{-1}(v) \subset X$ equipped with the induced path metric $d_{v}$ is
		a geodesic metric space $X_v$. Also, the
		inclusion maps ${i_v}:{X_v}\rightarrow{X}$
		are uniformly proper, i.e.\ $\forall M > 0$, $v\in{T}$ and $x, y\in{X_v}$,
		there exists $N > 0$ such that $d({i_v}(x),{i_v}(y)) \leq M$ implies
		${d_{v}}(x,y) \leq N$.
		\item Let $e=[v_1,v_2] \in \EE(T)$ with initial and final vertices $v_1$ and
		$v_2$ respectively (we assume that all edges have length 1).
		Let $X_e$ be the pre-image under $P$ of the mid-point of  $e$.
		There exist continuous maps ${f_e}:{X_e}{\times}[v_1,v_2]\rightarrow{X}$, such that
		$f_e{|}_{{X_e}{\times}(v_1,v_2)}$ is an isometry onto the pre-image of the
		interior of $e$ equipped with the path metric $d_e$. \\
		Further, we demand that $f_e$ is fiber-preserving,
		i.e. projection to the second co-ordinate in ${X_e}{\times}[v_1,v_2]$ corresponds via $f_e$
		to projection to the tree $P: X \rightarrow T$.
		\item  ${f_e}|_{{X_e}{\times}\{{v_1}\}}$ and
		${f_e}|_{{X_e}{\times}\{{v_2}\}}$ are $(K,{\epsilon})$-quasi-isometric
		embeddings into $X_{v_1}$ and $X_{v_2}$ respectively.
We shall often use the shorthand $f_{e,v_1}$ and $f_{e,v_2}$ for	${f_e}|_{{X_e}{\times}\{{v_1}\}}$ and
		${f_e}|_{{X_e}{\times}\{{v_2}\}}$ 
		 respectively. 
	\end{enumerate}
	We shall refer to
	$K, \epsilon$ as the constants or parameters of the  qi-embedding condition.
\end{definition}

If there exists $\delta > 0$ such that the vertex and edge spaces $X_v, X_e$ above are all
$\delta$-hyperbolic metric spaces for all vertices $v$ and edges $e$ of $T$, then
 $P: X \rightarrow T$ 
will be called \emph{ a tree of
hyperbolic metric
spaces}.

\begin{definition}\label{definitionofhallway}\cite{BF}
	A continuous map $f : [-k,k]{\times}{I} \rightarrow 
	X$ is called a \emph{hallway} of length $2k$ if it satisfies the following:
	
	\begin{enumerate}
		\item $f^{-1} ({\cup}{X_v} : v \in T) = \{-k,  \cdots , k \}{\times}
		I$
		\item $f$ is transverse, relative to condition (1) to $\cup_e X_e$.
		\item for all $i \in \{-k,  \cdots , k \}$, $f$ maps $i{\times}I$ to a geodesic in  $X_v$ for some vertex
		space $X_v$.
	\end{enumerate}
\end{definition}

\begin{definition}\label{def-rhothin}\cite{BF} A hallway $f : [-k,k]{\times}{I} \to 
	X$ is said to be \emph{$\rho$-thin} if 
	\[d({f(i,t)},{f({i+1},t)}) \leq \rho\] for all $i, t$.
	
	A hallway $f : [-m,m]{\times}{I} \rightarrow 
	X$ is called \emph{$\lambda$-hyperbolic}  if 
	\[\lambda l(f(\{ 0 \} \times I)) \leq \, {\rm max} \ \{ l(f(\{ -m \} \times I)),
	l(f(\{ m \} \times I)) \}.\]

The \emph{girth} of the hallway	is defined to be the quantity \[{\rm min_i} \, \{ l(f(\{ i \} \times I))\}.\]

	A hallway is \emph{essential} if the edge path in $T$ 
	resulting from projecting the hallway under $P\circ f$
	onto $T$ does not backtrack (and is therefore a geodesic segment in
	the tree $T$).
\end{definition}

\begin{definition}[Hallways flare condition \cite{BF}:]\label{def-flare}
	The tree of spaces, $X$, is said to satisfy the \emph{hallways flare}
	condition if 
	there exist $\lambda > 1$ and $m \geq 1$ such that
	the following holds:\\
	$\forall \rho > 0$ there exists $H(=H(\rho ))$ such that  any
	$\rho$-thin essential hallway of length $2m$ and girth at least $H$ is
	$\lambda$-hyperbolic. 
	
	The constants $\lambda, m$ are referred to as the constants or parameters of the hallways flare condition. If, further, the constant $\rho$ is fixed,
	then $H$ will also be called a constant or parameter of the hallways flare condition.
\end{definition}

With these notions in place, we can state the main geometric combination theorem of \cite{BF}:

\begin{theorem}\label{thm-bf}
Let $P: X \to T$ be a tree of hyperbolic spaces satisfying the qi-embedded condition (as in Definition \ref{def-tree}).
Further, suppose that the hallways flare condition (as in Definition \ref{def-flare}) is satisfied. Then $X$ is hyperbolic. 
\end{theorem}

The proof of Theorem \ref{thm-bf}
in \cite{BF} proceeds by establishing a linear isoperimetric inequality ensuring hyperbolicity. We shall indicate a different proof scheme below in the special case that the edge-to-vertex inclusion maps are uniform quasi-isometries rather than qi-embeddings. The forthcoming
book \cite{kap-sar} provides a new proof as well.

\subsection{Metric Bundles}\label{sec-mbdl}
The notion of a metric bundle \cite{mahan-sardar} adapts the idea of a fiber bundle to a coarse geometric context. We shall describe below
the main combination theorem of \cite{mahan-sardar} which is
an analog of Theorem \ref{thm-bf} in this context.

\begin{definition}\label{def-mbdl}
	Let $(\XX,d_X)$ and $(\BB, d_B)$ be geodesic metric spaces. Let $c, K\geq 1$ be  constants and 
	$h:{\mathbb R}^+ \rightarrow {\mathbb R}^+$  a function.
	$P: \XX\to \BB$ is called an \emph{$(h,c,K)-$ metric bundle} if
	\begin{enumerate}
		\item $P$ is 1-Lipschitz.
		\item For each  $z\in \BB$, $F_z=P^{-1}(z)$ is a geodesic metric space
		with respect to the path metric $d_z$ induced from $(\XX,d_X)$. We refer to $F_z$ as the \emph{fiber} over $z$.
		
		We further demand that the inclusion maps
		$i_z: (F_z,d_z) \rightarrow \XX$ are uniformly metrically proper as measured with respect to $h$, i.e.\ for  all $z\in \BB$ and $u, v\in F_z$,  $d_X(i_z(u),i_z(v))\leq N$ implies that $d_z(u,v)\leq f(N)$. 
		\item For $z_1,z_2\in \BB$ with $d_B(z_1,z_2)\leq 1$, let $\gamma$ be
		a geodesic in $\BB$ joining them. 
		Then for any $z\in \gamma$ and $x\in F_z$,  there is a path in $p^{-1}(\gamma)$
		of length at most $c$ joining $x$ to both $F_{z_1}$ and $F_{z_2}$.
		\item For $z_1,z_2\in \BB$ with $d_B(z_1,z_2)\leq 1$ and $\gamma \subset \BB$
		a geodesic  joining them,
		let $\phi: X_{z_1}\rightarrow X_{z_2}$, be any map such that
		for all $ x_1\in X_{z_1}$ there is a path of length at most $c$ in $P^{-1}(\gamma)$
		joining $x_1$ to $\phi(x_1)$. Then $\phi$ is a $K-$quasi-isometry.
	\end{enumerate}
	If in addition, there exists $\delta'$ such that each $X_z$ is $\delta'-$hyperbolic, then $P: \XX\to \BB$ is called an $(h,c,K)-$  metric bundle
	of $\delta'-$hyperbolic spaces (or simply a metric bundle
	of hyperbolic spaces if the constants are implicit).
\end{definition}
It is pointed out in \cite{mahan-sardar} that  condition (4) follows from the previous three (with suitable $K$); but it is more convenient to have it as part of our definition. 

A closely related notion of a \emph{metric graph bundle} often turns out to be more useful:

\begin{definition}\label{defn-mgbdl}\cite[Definition 1.2]{mahan-sardar}
	Suppose $\mathcal X$ and $\BB$ are metric graphs and $f:\mathbb N \rightarrow \mathbb N$ is a proper function.
	We say that $\XX$ is an \emph{$f$-metric graph bundle} over $\BB$ if there exists a surjective simplicial 
	map $\pi:\XX\rightarrow \BB$  such that the following hold.
	\begin{enumerate}
		\item For all $b\in V(\BB)$, $\FF_b:=\pi^{-1}(b)$ is a connected subgraph of $\XX$. Moreover, the inclusion maps
		$\FF_b\rightarrow \XX$, $b\in V(\BB)$ are uniformly metrically proper as measured by $f$.
		\item For all adjacent vertices $b_1,b_2\in V(\BB)$, any $x_1\in V(\FF_{b_1})$ is connected 
		by an edge to some $x_2\in V(\FF_{b_2})$.
	\end{enumerate}
\end{definition}

For all $b\in V(B)$, $F_b$ is called the \emph{ fiber} over $b$ and its path metric is denoted by $d_b$. It is pointed out in \cite{mahan-sardar} that any metric bundle is quasi-isometric to
a metric graph bundle, where the quasi-isometry coarsely preserves fibers and restricts to a quasi-isometry of fibers.
Condition (2) of Definition \ref{defn-mgbdl} immediately shows that if $\pi:\XX\rightarrow \BB$ is a metric graph bundle then for any points $v,w\in V(\BB)$ we have that 
	$Hd(F_v, F_w)<\infty$, where $Hd$ denotes the Hausdorff distance.

\begin{eg}\label{eg-es}
Let \[1 \to N \to G \to Q \to 1\] be an exact sequence of finitely generated groups.
Choose a finite generating set of $N$ and extend it 
to a finite generating set
of $G$. The image of the finite generating set of $G$ in $Q$ under the quotient map
is then a generating set of $Q$.  This gives a natural 
simplicial map $P: \Gamma_G \to \Gamma_Q$ between the respective Cayley graphs. This is the prototypical example of a metric graph bundle.
 The fibers are all copies of $\Gamma_N$.
\end{eg}

\begin{definition}
	Suppose $\XX$ is an $f$-{\em metric graph bundle} over $\BB$.
	Given $k\geq 1$ and a connected subgraph $\AAA\subset \BB$, a {\em $k$-qi section} over $\AAA$ is
	a map $s:\AAA\to \XX$ such that $s$ is a $k$-qi embedding and $\pi\circ s$ is the identity
	map on $\AAA$.
\end{definition}

For any hyperbolic metric space $F$ with more than two points in its Gromov boundary $\partial F$, there is a coarsely well-defined \emph{barycenter map} \[\phi : \partial^3 F \rightarrow F\] 
mapping an unordered triple $(a,b,c)$ of distinct points in $\partial F$ to a centroid of the ideal triangle spanned by $(a,b,c)$. We shall say that the
barycenter map $\phi : \partial^3 F \rightarrow F$ is  $N-$coarsely surjective if $F$ is contained in the $N$-neighborhood of the image of $\phi$. 
A \emph{$K-$qi-section} $\sigma: \BB\to \XX$ is a $K-$qi-embedding from $\BB$ to $\XX$ such that $P \circ \sigma $ is the identity map.
The following  guarantees the existence of qi-sections for metric bundles:

\begin{prop}\label{qi-section}\cite[Section 2.1]{mahan-sardar} Given $\delta,N, c, K\geq 0$ and 
	proper $f:{\mathbb{N}} \rightarrow {\mathbb{N}} $, there exists $K_0$ such that the following holds.\\
	Let $P : \XX \rightarrow \BB$ be an $(f,c, K)$-metric bundle of $\delta-$hyperbolic spaces such that all barycenter maps $\phi_b : \partial^3 F_b \rightarrow F_b$ are  $N-$coarsely surjective, 
	Then  through each point of $X$, there exists a $K_0$-qi section.\\
	
	A similar statement holds for metric graph bundles.
\end{prop}

The following gives the analog of Definition \ref{def-flare} in the context of metric bundles and metric graph bundles:

\begin{definition}\label{defn-flare}
	Let $P:\XX\rightarrow \BB$ be a metric bundle or a metric graph bundle.
 $P:\XX\rightarrow \BB$ is said to satisfy a \emph{flaring condition} if $\forall k \geq 1$, there exist
	$\lambda_k>1$ and  $n_k,M_k\in \mathbb N$ such that
	the following holds:\\
	Let $\gamma:[-n_k,n_k]\rightarrow \BB$ be a geodesic and let
	$\tilde{\gamma_1}$ and $\tilde{\gamma_2}$ be two
	$k$-qi sections of $\gamma$ in $X$.
	If $d_{\gamma(0)}(\tilde{\gamma_1}(0),\tilde{\gamma_2}(0))\geq M_k$,
	then 
	\[\lambda_k.d_{\gamma(0)}(\tilde{\gamma_1}(0),\tilde{\gamma_2}(0))\leq \max\{d_{\gamma(n_k)}(\tilde{\gamma_1}(n_k),\tilde{\gamma_2}(n_k)),d_{\gamma(-n_k)}(\tilde{\gamma_1}(-n_k),\tilde{\gamma_2}(-n_k))\}.
	\]
\end{definition}

The following Theorem is the analog of Theorem \ref{thm-bf} in the context of metric (graph) bundles.
\begin{theorem}\label{effectiveBF} Suppose that 
	$P:\XX\to \BB$ is a metric bundle
	or a metric  graph bundle  such that all fibers $F_z$ are uniformly hyperbolic, and
	 the barycenter maps are uniformly coarsely surjective. Equivalently, by Proposition \ref{qi-section}, there exists $\rho \geq 1$ such that
		for every $x \in \XX$, there exists a $\rho-$qi section $s: \BB\to \XX$ passing through $x$, i.e.\ $s\circ P (x) =x$.

	Then if $\XX$  satisfies the 
	qi-embedded condition  and the   flaring
	condition (as in Definition \ref{defn-flare})  corresponding to $\rho-$qi sections,
	then $\XX$ is  hyperbolic. \\
	
	Conversely,  if 
	$X$ is hyperbolic,
	then as a metric bundle or metric graph bundle, $\XX$ satisfies the flaring condition. 
\end{theorem}

\subsubsection{Ladders} A tool that has turned out to be considerably useful
in the context of both trees of spaces and metric bundles is the notion of a ladder. In particular, for our proof of Theorem \ref{effectiveBF} (sketched in Section \ref{pfofflare}), we shall use it.
The notion is related to, but different from that of a hallway.
Ladders were introduced in \cite{mitra-trees} in the context of trees of spaces and in \cite{mitra-ct}
in the context of groups. Instead of going through the construction in detail,  we extract the relevant features from the ladder construction of \cite{mitra-trees,mitra-ct}. The following is a restatement of \cite[Theorem 3.6]{mitra-trees} reformulated to emphasize the connection with hallways.

\begin{theorem}\label{qcladder} Given $\delta \geq 0, K\geq 1, \ep\geq 0$ there exists $D$ such that the following holds.\\ We consider one of the two following situations:
	\begin{enumerate}
		\item  $P:\XX \to \TT$ is a tree of $\delta-$hyperbolic spaces  as in Definition \ref{def-tree} with parameters $K, \ep$. Let $F_v$ be a vertex space,
		\item $P: \XX \to \BB$ is a metric bundle or metric graph bundle, and $F_v$
		is a fiber. 
	\end{enumerate}  
In both cases, the intrinsic metric on $F_v$ is denoted by $d_v$.
Then for every geodesic segment $\mu \subset (F_v, d_v)$ there exists a $D-$qi-embedded subset $\LL_\mu$ of $X$ such that the following holds.
	\begin{enumerate}
		\item $F_v \cap \LL_\mu = \mu$,
		\item (a) For $P:\XX\to \TT$ a tree of hyperbolic metric spaces and every $w \in \TT$, $F_w \cap \LL_\mu$ is either empty or a geodesic $\mu_w$ in $(F_w, d_w)$. Further, there exists a subtree $\TT_1 \subset \TT$ such that the collection of vertices $w \in \TT$ satisfying $ F_w \cap \LL_\mu \neq \emptyset$ equals the vertex set of $T_1$.\\
		(b) For $P: \XX \to \BB$ a metric bundle or metric graph bundle, 
		 $F_w \cap \LL_\mu$ is a geodesic $\mu_w$ in $(F_w, d_w)$. 
		\item There exists $ \rho_0\geq 1$ such that  through every $z \in \LL_\mu$, there exists a $ \rho_0-$qi-section $\sigma_z$ of $[v, P(z)]$ contained in $\LL_\mu$ satisfying $$ \sigma_z (P(z)) = z, \quad
		\sigma_z (v) \in \mu.$$
	\end{enumerate}
	Further, there exists a $D-$coarse Lipschitz retraction $\Pi_\mu:\XX \to \LL_\mu$, i.e.\
	\begin{enumerate}
		\item $d(\Pi_\mu(x), \Pi_\mu(y)) \leq D d(x, y) + D, \, \forall \, x, y \in X$,
		\item $\Pi_\mu(x) = x, \, \forall \, x \in \LL_\mu$.
	\end{enumerate}
\end{theorem}

The qi-embedded set $\LL_\mu$ is called a \emph{ladder} in \cite{mitra-ct,mitra-trees}. 
Theorem \ref{qcladder} shows in particular that there is a $(2D,2D)-$  quasigeodesic of $(X,d_X)$ joining the end-points of $\mu$ and lying on 
$\LL_\mu$. 

\begin{rmk}
Note that in Theorem \ref{qcladder}, we have \emph{not} assumed that $X$ is hyperbolic: no assumptions on the global geometry of $X$ are necessary here.
\end{rmk}

\subsubsection{Idea behind the proof of Theorem \ref{effectiveBF}}\label{pfofflare} We focus on the metric graph bundle case for convenience.
Theorem \ref{qcladder} guarantees that for any pair of points $x, y$  in a metric graph bundle
 $\XX$, there exist 
 \begin{enumerate}
 \item Qi-sections $\Sigma_x, \Sigma_y$ through $x, y$.
 \item A ladder $\LL(x,y)$ bounded by $\Sigma_x, \Sigma_y$. In fact, in this case (as shown in \cite{mitra-ct,mahan-sardar}), $\LL(x,y) \cap F_b$ is equal to a geodesic in $F_b$ joining
 $\Sigma_x(b), \Sigma_y(b)$ (here we are abusing notation slightly by identifying the qi-sections $\Sigma_x,
 \Sigma_y$ with their images).
 \end{enumerate}

Thus, for every $x, y\in \XX$ there are preferred quasigeodesics in $\XX$ contained in $\LL(x,y)$. We have not used the flaring condition so far. The flaring condition guarantees
hyperbolicity of $\LL(x,y)$. We shall return to this shortly. Hyperbolicity of $\LL(x,y)$
ensures (by the Morse Lemma) that all quasigeodesics in $\LL(x,y)$ joining $x, y$ are in a bounded neighborhood of each other. This gives a family of paths in $\XX$, one for every pair
$x,y$.
We then use a path-families argument following Hamenst\"adt
\cite{hamen} and a criterion due to Bowditch to conclude that $\XX$ is hyperbolic.

We return to the proof of hyperbolicity of $\LL(x,y)$. We note that $\LL(x,y)$ is a bundle
over $\BB$ where the fibers are intervals. The flaring condition is inherited by $\LL(x,y)$
with slightly worse constants.
Thus, we are reduced to proving Theorem \ref{effectiveBF} in the special case that fibers are intervals. To do this, we decompose the ladder $\LL(x,y)$ using qi-sections contained in
$\LL(x,y)$ into a finite number of ladders `stacked one on top of another'. Thus, there 
exist disjoint sections $\Sigma_x=\Sigma_0, \Sigma_1, \cdots , \Sigma_n=\Sigma_y$ and ladders
$\LL_i$ bounded by $\Sigma_{i-1}, \Sigma_i$ such that distinct $\LL_i$'s have disjoint interiors.
The ubiquity of qi-sections allows us to ensure that each of these smaller ladders has bounded
girth (in the spirit of Definition \ref{def-rhothin}), i.e. $\Sigma_{i-1}, \Sigma_i$ are at a bounded distance from each other along some fiber $F_b$ and flare away from each other as one goes to infinity in $B$. A further path families argument following \cite{hamen} allows us
to prove that $\LL(x,y)$ is hyperbolic.

A word about the proof sketch above. Note that we use only the 1-dimensional property
of quasigeodesics flaring and path families to prove the combination theorem in this case,
as opposed to the more `2-dimensional' area argument of \cite{BF}. This has been considerably refined in \cite{kap-sar} to give a new path-families proof of Theorem \ref{thm-bf}.

\subsection{Relatively hyperbolic combination theorems}\label{sec-effrhcomb} 
We refer the reader to \cite{farb-relhyp,gromov-hypgps,bowditch-relhyp} for the basics of relative hyperbolicity.
Theorem \ref{thm-bf} was generalized to the context of trees of relatively hyperbolic spaces in two different ways:

\begin{enumerate}
	\item Using an acylindricity hypothesis in \cite{dahmani-combnthm} and \cite{alib}. This is
	in the spirit of Theorem \ref{thm-acyl}.
\item Using the flaring condition in \cite{mahan-reeves,gautero}. This in the spirit of Theorem \ref{thm-dlt}.
\end{enumerate}

\subsubsection{Relatively hyperbolic combination theorem using acylindricity}
Let $G$ be hyperbolic relative to a finite collection $\PP = \{P_1, \cdots , P_k\}$ of parabolic subgroups. Let $\partial_h G$ denote the Bowditch boundary of $G$.
 Let $H \subset G$ be a relatively quasiconvex subgroup \cite{hruska-relqc}. We shall give 
 Dahmani's version of the combination theorem \cite{dahmani-combnthm}
  below. Let $\Lambda_H \subset \partial_hG$ denote the  limit set of $H$.
  A relatively quasiconvex subgroups $H$ is \emph{full relatively quasi-convex} if  it  is  quasi-convex  and  if,  for  any  infinite  sequence $g_n \in G$ in distinct left cosets of $H$, the intersection  $\cap_n g_n(\Lambda_H)$ is empty.

 \begin{lemma}\cite[Lemma 1.7]{dahmani-combnthm}
 Let
 $G$ be hyperbolic relative to a finite collection $\PP = \{P_1, \cdots , P_k\}$ of parabolic subgroups. Let H be a full relatively quasi-convex subgroup.  Let $P$ be a conjugate of one of the $P_i$'s. Then $P\cap H$ is either  finite, or of  finite index in $P$.
 \end{lemma}

\begin{definition}\label{def-acyltree}\cite{sela-acyl}
 The action of a  group $G$ on a tree $T$ is $k-$acylindrical for some $k\in \natls$ if the stabilizer of any geodesic of length $k$ in $T$ is finite.
 The action of a  group $G$ on a tree $T$ is acylindrical if it is
 $k-$acylindrical for some $k\in \natls$

 A finite graph of groups is said to be acylindrical, if the action on the associated Bass--Serre tree is acylindrical.
\end{definition}

Then Dahmani's combination theorem states:

\begin{theorem}\cite{dahmani-combnthm}
Let $G$ be the fundamental group of an acylindrical  finite graph of relatively hyperbolic groups, whose edge groups are full quasi-convex subgroups of the adjacent  vertex  groups.  Let $\GG$ be the  family  of   images  of  the  maximal parabolic subgroups of the vertex groups, and their conjugates in $G$.  Then $G$ is strongly  hyperbolic relative to $\GG$.
\end{theorem}

The  approach in \cite{dahmani-combnthm} is quite different from that of \cite{BF}.
From the Bowditch boundaries of the vertex and edge groups, a metrizable compact space $Z$
is constructed in such a way that $G$ naturally acts on $Z$. 
 It is then shown that this action is a convergence action. Finally, it is shown that the action is geometrically finite, forcing $G$ to have a relatively hyperbolic structure.
 
 The Bass--Serre tree of $G$ has vertex groups $G_v$ and edge groups $G_e$. Hence, associated
 to the 
 Bass--Serre tree $T$ there is a natural tree ($T$) of compact spaces given by $\partial_h G_v$
 and $\partial_h G_e$.
 The set $Z$ is built \cite[Section 2]{dahmani-combnthm} from these copies of  $\partial_h G_v$
 and $\partial_h G_e$. Suppose $e=[v_1,v_2]$ is an edge of $T$.
 For all such edges $e$, glue together  $\partial_h G_{v_1}$
 and $\partial_h G_{v_2}$ along the limit set $\partial_h G_{e}$. The relevant identification space is thus obtained from the set $\sqcup_{v \in V(T)}\partial_h G_{v}$
 by identifying pairs of points according to the images of $\partial_h G_{e}$. Finally, 
 the base tree $T$ encodes (infinite) directions that are `transverse' to all the vertex spaces. The set $Z$ is then obtained from topologizing $\partial T \cup \sqcup_{v \in V(T)}\partial_h G_{v}/\ \sim$, where $\sim$ is the equivalence relation given by edge spaces.
 
 Alibegovic \cite{alib} proves a similar combination theorem for relatively hyperbolic groups 
 following the original strategy of Bestvina and Feighn in Theorem \ref{thm-bf} using
 the linear isoperimetric inequality characterization of hyperbolicity.

\subsubsection{Relatively hyperbolic combination theorem using flaring} We next define
a tree of relatively hyperbolic spaces in general.
\begin{definition}\label{def-treerh}\cite{mahan-reeves} A tree $P: \XX \rightarrow \TT$ of geodesic
	metric spaces is said to be a  tree of relatively hyperbolic metric spaces
	if in addition to the conditions of Definition \ref{def-tree}
	\begin{enumerate}
		\item[(4)] each  vertex space $X_v$ is strongly hyperbolic relative to a collection of subsets $\HH_v$ and  each 
		edge space $X_e$ is strongly hyperbolic relative to a collection of subsets $\HH_e$.  The  sets
		$H_{v,\alpha}\in \HH_{v}$ or $H_{e,\alpha}\in \HH_{e}$ are referred to as \emph{horosphere-like sets}.
		\item[(5)]  the maps $f_{e,v_i}$ above, for $i = 1, 2$, are \emph{
			strictly type-preserving}. That is, for $i =
		1, 2$ and for
		any $H_{v_i,\alpha}\in \HH_{v_i}$, $f_{e,v_i}^{-1}(H_{v_i,\alpha})$,
		is
		either empty or is equal to some $H_{e,\beta}\in \HH_{e}$.
		 Further, for all
		$H_{e,\beta}\in \HH_{e}$, there exists $v$ and
		$H_{v,\alpha}$, such that $f_{e,v} ( H_{e,\beta}) \subset H_{v,\alpha} $.
		\item[(6)]  There exists $\delta > 0$ such that each $\EE (X_v, \HH_v )$ is $\delta$-hyperbolic
		(here, $\EE (X_v, \HH_v )$ denotes the electric space obtained  from $X_v$ by electrifying all the horosphere-like sets in $\HH_v$).
		\item[(7)]  The induced maps of the coned-off edge spaces into the
		coned-off vertex spaces $\hhat{f_{e,v_i}} : \EE ({X_e}, \HH_e ) \rightarrow
		\EE ({X_{v_i}}, \HH_{v_i})$ ($i = 1, 2$) are uniform quasi-isometries. This is called the
		\emph{qi-preserving electrification condition}
	\end{enumerate}
\end{definition}

We state conditions (4) and (6)  in conjunction by saying that $X_v$ is \emph{strongly $\delta-$hyperbolic} relative to $\HH_v$.

We explain condition (7) briefly.
Given the tree of spaces $P: \XX \to \TT$ with vertex spaces $X_v$ and edge spaces $X_e$ there exists a naturally associated  tree whose vertex spaces are the electrified spaces
$\EE (X_v, {\HH}_v)$ and edge spaces are the electrified spaces
$\EE (X_e, {\HH}_e)$ obtained by electrifying the respective horosphere like sets.
Condition (4) of the above definition ensures that we have natural inclusion maps of edge spaces
$\EE (X_e, {\HH}_e)$ into adjacent vertex spaces $\EE (X_v, {\HH}_v)$.
The resulting tree of coned-off spaces $P: \TC (X) \rightarrow \TT$
is referred to simply as the \emph{ induced
	tree of coned-off spaces}. 
The \emph{cone locus} of $\TC (X)$ is
the  forest given by the following:
\begin{enumerate}
\item the vertex set $\VV(\TC (X))$ consists of the
cone-points $c_{v,\alpha}$ in the vertex spaces $X_v$ resulting from the electrification operation of the horosphere-like sets $H_{v,\alpha} \in \HH_v$.
\item the edge set $\EE(\TC (X))$ consists of the
cone-points $c_{e,\alpha}$ in the edge set $X_e$ resulting from the electrification operation of the horosphere-like sets $H_{e,\alpha} \in \HH_e$.
\end{enumerate}

Each connected component of the
cone-locus is a \emph{maximal cone-subtree}. The collection
of maximal cone-subtrees is denoted by $\CC\TT$ and elements
of $\CC\TT$ are denoted as $CT_\alpha$. Note that each maximal
cone-subtree $CT_\alpha$ naturally gives rise to a tree $CT_\alpha$ of
horosphere-like subsets depending on which cone-points arise as
vertices and edges of $CT_\alpha$. The metric space that $CT_\alpha$
gives rise to is denoted as $C_\alpha$. We refer to 
any such $C_\alpha$ as
a \emph{maximal cone-subtree of horosphere-like spaces}. The induced tree of horosphere-like sets is denoted by $$g_\alpha : C_\alpha
\to CT_\alpha.$$
The collection of these maps will be denoted as $\GG$.
The collection
of the maximal cone-subtree of horosphere-like spaces $C_\alpha$ is denoted as $\CC$. 
Note  that each $CT_\alpha$ thus appears both as a subset of $\TC (X)$ as well as the underlying tree of $C_\alpha$.

\begin{definition}[Cone-bounded hallways strictly flare
		condition:]
		
		An essential  hallway of length $2k$ is \emph{cone-bounded} if
		$f(i \times {\partial I})$
		lies in the cone-locus for $i = \{ -k, \cdots , k\}$.
		
	The tree of spaces, $X$, is said to satisfy the \emph{ cone-bounded hallways flare}
	condition if there are numbers $\lambda > 1$ and $k \geq 1$ such that
	any
	cone-bounded hallway of length $2k$  is
	$\lambda$-hyperbolic, where $\lambda, k$ are called the constants or parameters of the strict flare condition.
\end{definition}

We now state the combination theorem for relative hyperbolicity using the flaring condition.

\begin{theorem}\label{effectiverelBF}\cite{mahan-reeves,gautero}
	Let $P:\XX\to \TT$ be a  tree of 
	uniformly relatively hyperbolic spaces in the sense of Definition \ref{def-treerh}
	satisfying the
	qi-embedded condition, such that the resulting tree of coned-off spaces satisfies
	\begin{enumerate}
	\item the  hallways flare
	condition,
	\item the cone-bounded hallways  flare condition.
	\end{enumerate} Then $X$ is   hyperbolic relative to 
the maximal cone-subtrees of horosphere-like spaces. 
\end{theorem}

\subsection{Effective quasiconvexity and flaring}\label{sec-effecqc}  We now state a couple of theorems along the lines of Theorem \ref{thm-bf} and \ref{effectiveBF}  ensuring quasiconvexity of a subspace of a vertex space. The first, due to Ilya Kapovich \cite{kapovich-combin} is in the setup of an acylindrical graph of groups:

\begin{theorem}\label{ilya} Let $\GG$ be a
	finite acylindrical graph of groups
where all vertex and edge groups are hyperbolic and edge-to-vertex inclusions
are  quasi-isometric embeddings.
Let $\TT \subset \GG$ be a maximal subtree. Let $G$ denote the group  corresponding to the
tree $\TT$. (By Theorem \ref{thm-bf}, $G$ is hyperbolic.)
Then each vertex group $G_v$ of $G$ is quasiconvex in $G$.
\end{theorem}

Next, we shall consider in a unified way the two following situations:
\begin{enumerate}
	\item $P: \XX \to \TT$ is a tree  of hyperbolic metric spaces satisfying the 
	qi-embedded condition with constants $K, \ep$ and the  hallways flare
	condition with parameters $\lambda_0, m_0$. Further, if $\rho_0$ is given we shall assume an additional constant $ H_0$ as a lower bound  for girths of $\rho_0-$thin hallways. 
	\item $P: \XX \to \BB$ is a metric bundle or a metric graph bundle satisfying the flaring condition with constants as in Definitions \ref{def-mbdl}, \ref{defn-mgbdl}, \ref{defn-flare}.
\end{enumerate}

Also $(X_v,d_v)$ will  be a vertex space of $\XX$ (in the tree of spaces case) or $P^{-1}(v)$ equipped with the induced metric in the metric (graph) bundle case.

\begin{definition}\label{def-flareinall} Let $P: \XX\to \TT$ be a tree of hyperbolic spaces.  Let $Y\subset (X_v,d_v)$  be a $C-$quasiconvex subset of $(X_v,d_v)$.
	We say that $Y$ \emph{flares in all directions with parameter $K$} if for any geodesic segment $[a,b] \subset (X_v,d_v)$ with $a, b \in Y$ and any $\rho-$thin hallway $f:[0,k] \times I \to \XX$ satisfying 
	\begin{enumerate}
		\item $\rho \leq \rho_0$,
		\item $f(\{0\} \times I) = [a,b]$,
		\item $l([a,b]) \geq K$,
		\item $k \geq K$,
	\end{enumerate}
	the length of $f(\{k\} \times I)$ satisfies $$l(f(\{k\} \times I)) \geq \lambda l( [a,b]).$$
	
	Similarly, let $P: \XX\to \BB$ be a metric bundle or metric graph bundle with hyperbolic fiber. Let $Y \subset X_v$ be quasiconvex. Further, assume that there is a $\rho-$qi section through every $x \in \XX$ (cf.\ the second hypothesis of Theorem \ref{effectiveBF}).
	
		We say that $Y$  \emph{flares in all directions with parameter $K \geq 0, D \geq 1, \lambda > 1$ } if
	the following holds:\\
	Let $\gamma:[0,D]\rightarrow \BB$ be a geodesic such that $\gamma(0)=v$ and let
	$\tilde{\gamma_1}$ and $\tilde{\gamma_2}$ be two
	$\rho$-qi lifts (sections) of $\gamma$ in $\XX$.
	If $d_{v}(\tilde{\gamma_1}(0),\tilde{\gamma_2}(0))\geq K$,
	then we have

	 $$\lambda.d_{v}(\tilde{\gamma_1}(0),\tilde{\gamma_2}(0))\leq  d_{\gamma(D)}(\tilde{\gamma_1}(D),\tilde{\gamma_2}(D)).$$

\end{definition}

We can now state  a proposition guaranteeing quasiconvexity of subsets of vertex spaces.

\begin{prop}\label{effectiveqc}\cite{mahan-tight} Given $K, C$, there exists $C_0$ such that the following holds.\\
	Let $P: \XX \to \TT$  and $X_v$ be as in Theorem \ref{qcladder} above. If $Y$ is a $C-$quasiconvex subset of $(X_v,d_v)$ and flares in all directions with parameter $K$, then $Y$ is  $C_0-$quasiconvex in $(X,d_X)$.
	
	Conversely, given $C_0$, there exist $K, C$ such that the following holds.\\
	For $P: \XX \to \TT$  and $X_v$  as above, if $Y \subset X_v$ is   $C_0-$quasiconvex in $(X,d_X)$, then it is a $C-$quasiconvex subset in $(X_v,d_v)$ and flares in all directions with parameter $K$.
\end{prop}

A similar statement holds for metric (graph) bundles. 

\begin{prop}\label{effectiveqcmbdl} Given $K,D, \lambda,  C$, there exists $C_0$ such that the following holds.\\
	Let $P: \XX \to \BB$ be a metric (graph) bundle  and $X_v$ be as in Definition \ref{def-flareinall}. If $Y$ is a $C-$quasiconvex subset of $(X_v,d_v)$ and flares in all directions with parameters $K, D, \lambda$, then $Y$ is  $C_0-$quasiconvex in $(\XX,d_X)$.
	
	Conversely, given $C_0$, there exist $K,D, \lambda, C$ such that the following holds.\\
	For $P: \XX \to \BB$ a metric (graph) bundle  and $X_v$  as above, if $Y \subset X_v$ is   $C_0-$quasiconvex in $(X,d_X)$, then it is a $C-$quasiconvex subset in $(X_v,d_v)$ and flares in all directions with parameters $K, D, \lambda$.
\end{prop}

\section{Combination theorems in geometric group theory: cubulations}\label{sec-cube} We turn now to the remarkable work
during the last decade on special cube complexes. We refer to \cite{hw-gafa}
for the basics of special cube complexes. Let $\GG$ denote a finite graph, $RAAG(\GG)$ the right-angled Artin group associated to $\GG$, and $\SSS(\GG)$ its Salvetti complex.
A cube complex $\CC$ is said to be \emph{special} if 
there exists a combinatorial local isometry from $\CC$ to $\SSS(\GG)$ for some finite graph $\GG$. By Agol's resolution of Wise's conjecture in \cite{agol-vhak} (see Theorem \ref{thm-agol} below), hyperbolic groups that are virtually special are precisely those that act geometrically on a CAT($0$) cube complex. We give a brief account of some of the combination theorems that have been proved around this theme.

In \cite{hsuw-ajm}, Hsu and Wise proved the precursor of all virtually special combination theorems by showing that if a hyperbolic group $G$ splits as a finite graph of finitely generated free groups with cyclic edge groups, then $G$ is virtually special. In \cite{hsuw-in}, they later generalized this to
amalgamated products of free groups over a finitely generated malnormal subgroup.
A landmark combination theorem due to Haglund and Wise concerns the combination of  hyperbolic virtually special cubulable groups along malnormal quasiconvex subgroups:

\begin{theorem} \cite{hw-annals}\label{wise-combo}
Let $A, B, M $
be compact virtually special cube complexes. Suppose that $G_A=\pi_1(A)$, $G_B= \pi_1(B)$, and
$G_M= \pi_1(M)$ are  hyperbolic. Let $M\stackrel{i_A}\longrightarrow A$, and  $M\stackrel{i_B}\longrightarrow B$ be local isometries of cube complexes such that  $i_{A\ast} (G_M)$ and $i_{B\ast} (G_M)$ are quasiconvex and malnormal in $G_A$ and $G_B$
respectively. Let $X=A\cup_MB$ be the cube complex obtained by gluing $A$ and $B$ together along $M$ using $M \times [0,1]$. Then $X$ is virtually special. 
\end{theorem} 

Theorem \ref{wise-combo} generalizes earlier work of Wise \cite{wise02} where he showed that
any 2-complex built by amalgamating (in terms of fundamental group) two finite graphs along a malnormal immersed  graph  is  virtually  special. Theorem \ref{wise-combo} is also  a crucial ingredient in Wise's proof of the virtual specialness of hyperbolic groups admitting a quasiconvex
hierarchy. This is very much in the spirit  of the Haken hierarchy for Haken 3-manifolds
and Thurston's hyperbolization of such manifolds (cf.\ Section \ref{sec-monster}).

\begin{theorem}\cite{wise}\label{thm-hier}
Let $G$ be a hyperbolic group admitting a quasiconvex hierarchy.
 Then $G$ is the fundamental group of a compact non-positively curved cube complex that is virtually special.
\end{theorem} 

We list below some of the important consequences of Theorem \ref{thm-hier}.
The following resolved a conjecture of Baumslag:
\begin{theorem}\cite{wise}
Every one-relator group with torsion is virtually special.
\end{theorem}

In the context of hyperbolic $3$-manifolds, Wise showed the following.
\begin{theorem}\cite{wise}
Compact hyperbolic Haken manifolds are  virtually special.
\end{theorem}

Using work of Kahn and Markovic \cite{kahnm}, Bergeron and Wise \cite{bergeronw} proved
that all hyperbolic 3-manifolds can be cubulated, i.e.\ they act geometrically on CAT(0) cube complexes. This led Wise to conjecture that hyperbolic groups that act geometrically on CAT(0) cube complexes are virtually special. The following celebrated theorem of Agol resolved this conjecture affirmatively:

\begin{theorem} \cite{agol-vhak}\label{thm-agol}
Hyperbolic groups acting geometrically on CAT(0) cube complexes are virtually special.
\end{theorem}

A flurry of activity ensued in trying to show that several naturally defined hyperbolic groups are, in fact, cubulable. In \cite{hw1,hw2}, Hagen and Wise proved that hyperbolic groups $G$ admitting an exact sequence of the form $$1 \to F_n \to G \to \Z \to 1$$ are cubulable.
(Here $F_n$ denotes the free group on $n$ generators.) Hence, by Agol's Theorem \ref{thm-agol}, such groups $G$ are virtually special. In a different direction, Manning,
the first author and Sageev \cite{mms} showed that there exist cubulable hyperbolic groups 
$G$ admitting an exact sequence of the form $$1 \to \pi_1(S) \to G \to F_n \to 1,$$ where $S$ is a closed surface of genus greater than one. Again, by Agol's Theorem \ref{thm-agol}, such groups $G$ are virtually special. 

Finally, we mention work of Przytycki and Wise \cite{przwise}, who proved the virtual specialness of fundamental groups  of $3$-manifolds whose JSJ decomposition has both a hyperbolic as well as a Seifert-fibered piece. $3$-manifolds admitting such a JSJ decomposition are called \emph{mixed}. As a consequence of their result, the authors show that mixed manifolds virtually fiber.

The proof in \cite{przwise} proceeds by first showing that there are enough codimension one surface subgroups to ensure cubulability. This is established  by combining  surfaces coming  from the graph manifold pieces with those coming  from hyperbolic pieces. Once cubulability has been established, the malnormal special quotient theorem  \cite{wise} (see also \cite{agm-msqt}) is used to establish specialness of 
the cube complex thus built. 

In this section, we have given only a cursory treatment of a topic that, starting with \cite{sageev-thesis} has undergone tremendous development over the last two decades. We refer the reader to the books \cite{afw} and \cite{wise-raags} for a more comprehensive treatment.

\section{Holomorphic dynamics and polynomial mating}\label{sec-cxdyn}

Before entering the theme of combination theorems in holomorphic dynamics, we say a few words on the history of the subject, and sketch briefly some of the philosophical parallels between holomorphic dynamics and Kleinian groups and some of the developments inspired by this synergy. These will also serve as a motivation for combination theorems involving complex polynomials and Kleinian groups that we will discuss later in the section.

\subsection{Historical comments}\label{hist_comm_subsec}
The study of dynamics of rational maps on the Riemann sphere started with groundbreaking work of Fatou and Julia \cite{fatou-1919,fatou-1920a,fatou-1920b,fatou-1926,julia-1918,julia-1922} in the 1920s. The subject remained dormant for several decades barring a handful of important contributions, most notably by Siegel \cite{siegel} and Brolin \cite{brolin}. Around the 1970s, the availability of  computers allowed Feigenbaum and Mandelbrot to perform numerical experiments on finer structures of dynamical and parameter spaces of real/complex-analytic maps. Their pioneering discoveries infused fresh blood into the field, and gave rise to problems and conjectures that played pivotal roles in the development of the modern theory of holomorphic dynamics.

A revolutionary contribution came from Sullivan, who introduced quasiconformal methods into the study of rational dynamics to prove nonexistence of wandering domains in the Fatou set for rational maps \cite{sullivan-dict}. The seminal work of Douady and Hubbard on the dynamics of quadratic polynomials and the structure of the Mandelbrot set \cite{orsay1,orsay2} turned out to be equally fundamental in that the techniques devised by them were robust enough to be applied to the study of a wide variety of holomorphic dynamical systems.

Sullivan proposed a dictionary between Kleinian groups and rational dynamics that was motivated by various common features shared by them \cite[p. 405]{sullivan-dict}. In addition to the apparent similarities between the topological structures of the \emph{limit set} (respectively, the \emph{domain of discontinuity}) of a Kleinian group and the \emph{Julia set} (respectively, the \emph{Fatou set}) of a rational map, there are deeper similarities between the techniques employed in proving various statements in these two parallel worlds. In fact, in the same paper, Sullivan gave a new proof of Ahlfors' finiteness theorem which closely parallels the proof of the `no wandering Fatou component' theorem for rational maps.

Around the same time, Thurston proved a topological characterization for an important class of rational maps \cite{DH1}. This result, which is a philosophical analog of the hyperbolization of atoroidal Haken $3$-manifolds, has given rise to a wealth of rich and beautiful results that we will not be touching upon in this survey (see \cite[\S 9]{rees-survey} and the references therein). 

We should emphasize that the aforementioned dictionary is not an automatic method for translating results in one setting to the other, but rather an inspiration for results and proof techniques. We now list a few prominent pieces of work motivated by this dictionary. Sullivan and McMullen introduced Teichm{\"u}ller spaces of conformal dynamical systems in the spirit of Teichm{\"u}ller spaces of Riemann surfaces in \cite{sul-mc}.  Bullett and Penrose constructed matings of holomorphic quadratic polynomials and the modular group as {\it holomorphic correspondences} \cite{bullett-penrose}. In \cite{lyubich-minsky}, Lyubich and Minsky constructed ``an explicit object that plays for a rational map the role played by the hyperbolic $3$-orbifold quotient of a Kleinian group''.McMullen established conceptual connections between renormalization ideas used in holomorphic dynamics and the study of $3$-manifolds fibering over the circle \cite{ctm-renorm}. Pilgrim proved a canonical decomposition theorem for Thurston maps as an analog of the torus decomposition theorem for $3$-manifolds \cite{pilgrim}. Another noteworthy development in the framework of the above dictionary is the recent work of Luo \cite{Luo19a,Luo21a,Luo21b}, where results in rational dynamics were proved using techniques that are closely related to Thurston's work on $3$-manifolds.

We refer the reader to \cite{milnor-book,CG1} for a basic introduction to rational dynamics, to \cite{lyubich-book} for a comprehensive account on the dynamics of quadratic polynomials and the Mandelbrot set, to the recent survey article by Rees on major advances in the field \cite{rees-survey}, and a survey by DeZotti on connections between holomorphic dynamics and other branches of mathematics \cite{dezotti-survey}.

\subsection{Mating of polynomials}\label{poly_mating_subsec}

The operation of polynomial mating, which was introduced by Douady and Hubbard in \cite{douady-mating}, constructs a rational map on $\widehat{\C}$ by combining the actions of two complex polynomials. Since the first appearance of the notion, several closely related definitions and perspectives have been put forward. In this survey, we will follow the route adopted in \cite{petersen-mayer-mating} (see \cite{roesch-intro} for the original formulation and some historical comments, and \cite{milnor-mating-example} for a lucid account of the mating construction along with a detailed worked out example).

To define the operation  of polynomial mating formally, we need to introduce some terminology. The \emph{Fatou set} of a rational map $R$, denoted by $\mathcal{F}(R)$, is the largest open subset of $\widehat{\C}$ on which the sequence of iterates $\{R^{\circ n}\}_{n\geq 0}$ forms a normal family. Its complement is called the \emph{Julia set}, and is denoted by $\mathcal{J}(R)$. For a complex polynomial $P$, the \emph{filled Julia set} (i.e., the set of points with bounded forward orbits) and the \emph{basin of attraction of infinity} (i.e., the complement of the filled Julia set) are denoted by $\mathcal{K}(P)$ and $\mathcal{B}_\infty(P)$, respectively. We refer the reader to \cite{milnor-book} for basic topological and dynamical properties of these sets.

A rational map $R$ is called \emph{postcritically finite} if each of its critical points has a finite forward orbit. $R$ is called \emph{hyperbolic} if each of its critical points converges to an attracting cycle under forward iteration.

If $P$ is a monic, centered polynomial of degree $d$ with a connected Julia set, then there exists a conformal map $\phi_P: \widehat{\C}\setminus\overline{\mathbb{D}}\rightarrow \mathcal{B}_\infty(P)$ that conjugates $z^d$ to $P$, and satisfies $\phi_P'(\infty)=1$ \cite[Theorem~9.1, Theorem~9.5]{milnor-book}. We will call $\phi_P$ the \emph{B\"ottcher coordinate} for $P$. Furthermore, if $\partial \mathcal{K}(P)=\mathcal{J}(P)$ is locally connected, then $\phi_P$ extends to a semiconjugacy between $z^{d}\vert_{\mathbb{S}^1}$ and $P\vert_{\mathcal{J}(P)}$. In this case, the map $\phi_P:\mathbb{S}^1\to\mathcal{J}(P)$ is called the \emph{Carath{\'e}odory loop/semi-conjugacy} for $\mathcal{J}(P)$.

Now let $P_1, P_2$ be two monic polynomials of the same degree $d\geq 2$ with connected and  locally connected filled Julia sets. We consider the disjoint union $\mathcal{K}(P_1)\sqcup\mathcal{K}(P_2)$ and the map 
\begin{center}
$P_1\sqcup P_2: \mathcal{K}(P_1)\sqcup \mathcal{K}(P_2)\to \mathcal{K}(P_1)\sqcup \mathcal{K}(P_2),$\\
 $P_1\sqcup P_2\vert_{\mathcal{K}(P_1)}=P_1,\quad  P_1\sqcup P_2\vert_{\mathcal{K}(P_2)}=P_2.$
\end{center}
Let $\sim$ be the equivalence relation on $\mathcal{K}(P_1)\sqcup\mathcal{K}(P_2)$ generated by $\phi_{P_1}(z)\sim \phi_{P_2}(\overline{z})$, for all $z\in\mathbb{S}^1$. It is easy to check that $\sim$ is $P_1\sqcup P_2-$invariant, and hence it descends to a continuous map $P_1\mate P_2$ to the quotient $\mathcal{K}(P_1)\mate\mathcal{K}(P_2):=\left(\mathcal{K}(P_1)\sqcup\mathcal{K}(P_2)\right)/\sim$ (see \cite[\S 4.1]{petersen-mayer-mating} for details). The map $P_1\mate P_2$ is called the \emph{topological mating} of the polynomials $P_1, P_2$. Moreover, if $\mathcal{K}(P_1)\mate\mathcal{K}(P_2)$ is homeomorphic to a $2$-sphere, we say that the topological mating is \emph{Moore-unobstructed}. We refer the reader to \cite[Theorem~2.12]{petersen-mayer-mating} for the statement of Moore's theorem, which provides a general sufficient condition for the quotient of $\mathbb{S}^2$ under an equivalence relation to be a topological $2-$sphere, and to \cite[Proposition~4.12]{petersen-mayer-mating} for a useful application of Moore's theorem giving a sufficient condition for the topological mating of $P_1, P_2$ (as above) to be Moore-unobstructed (note that the conditions of Moore's theorem are not necessary, see \cite[Example~13.18]{bonkmeyer}). By \cite[Proposition~4.3]{petersen-mayer-mating}, if the topological mating of $P_1, P_2$ is not Moore obstructed (i.e., if $\mathcal{K}(P_1)\mate\mathcal{K}(P_2)\cong \mathbb{S}^2$), then $P_1\mate P_2$ is topologically conjugate to an orientation-preserving branched covering of $\mathbb{S}^2$.
The following definition relates the topological mating to rational maps of $\widehat{\C}$. We refer the reader to \cite{DH1} for the notion of \emph{Thurston equivalence} appearing below.

\begin{definition}\cite[Definition~4.4]{petersen-mayer-mating}
Let the topological mating of $P_1\mate P_2$ be Moore-unobstructed, and $h:\mathcal{K}(P_1)\mate\mathcal{K}(P_2)\to\mathbb{S}^2$ be a homeomorphism.
\begin{enumerate}
\item The polynomials $P_1, P_2$ are called \emph{combinatorially mateable} if they are postcritically finite and if the branched covering $h\circ P_1\mate P_2\circ h^{-1}:\mathbb{S}^2\to\mathbb{S}^2$ is Thurston equivalent to a rational map $R$.

\item The polynomials $P_1, P_2$ are called \emph{conformally/geometrically mateable} if the homeomorphism $h$ can be so chosen that
$R=h\circ P_1\mate P_2\circ h^{-1}:\mathbb{S}^2\to\mathbb{S}^2$ is a rational map and $h$ is conformal on the interior of $\mathcal{K}(P_1)\mate\mathcal{K}(P_2)$.
\end{enumerate}

Conversely, a rational map $R$ is said to be combinatorially (respectively, conformally) a mating if there exist polynomials $P_1, P_2$ satisfying the corresponding property above with $R=h\circ P_1\mate P_2\circ h^{-1}$.
\end{definition}

The following equivalent definition of conformal mating is often useful in practice (see \cite[\S 4.7]{petersen-mayer-mating} for other definitions). In fact, this definition can be easily adapted for the other frameworks of combination theorems that we will discuss in this section (compare Definition~\ref{mating_def_2}).

\begin{definition}\label{mating_def_1}\cite[Definition~4.14]{petersen-mayer-mating}
A rational map $R:\widehat{\C}\to\widehat{\C}$ of degree $d\geq 2$ is said to be the \emph{conformal mating} of two degree $d$ monic, centered, polynomials $P_1$ and $P_2$ with connected and locally connected filled Julia sets if and only if there exist continuous maps
\[\psi_{1}: \mathcal{K}(P_1) \rightarrow \widehat{\mathbb{C}} \textrm{ and } \psi_{2}: \mathcal{K}(P_2) \rightarrow \widehat{\C},  \] conformal on $\textrm{int}~{\mathcal{K}(P_1)}$, $\textrm{int}~{\mathcal{K}(P_2)}$, respectively, such that 

\begin{enumerate}\upshape
\item $\psi_1(\mathcal{K}(P_1))\bigcup\psi_2(\mathcal{K}(P_2))=\widehat{\C}$,

\item $\psi_i\circ P_i=R\circ\psi_i$,\ \textrm{for}\ $i\in\{1,2\}$, and

\item $\psi_1(z)=\psi_2(w)$ if and only if $z\sim w$, where $\sim$ is the equivalence relation defined above.
\end{enumerate}
\end{definition}

With the above notions of mating in place, we can now mention the first major results on mateability of complex polynomials. In fact, these provided the first main application of Thurston's theorem on the topological characterization for rational maps. The following theorem completely answers the question of conformal mateability of postcritically finite quadratic polynomials (see \cite{orsay1,orsay2} for a detailed study of the Mandelbrot set, or \cite{milnor-mandelbrot} for a quick introduction).

\begin{theorem}\cite{rees-quad-rat,tan-mating,shishikura-mating}
Let $P_1(z) = z^2+c_1$ and $P_2(z) = z^2+c_2$ be two postcritically finite quadratic polynomials. Then $P_1$ and $P_2$ are conformally
mateable if and only if $c_1$ and $c_2$ do not belong to conjugate limbs of the Mandelbrot set.
\end{theorem}

Among other early works on matings of postcritically finite polynomials, we ought to mention the work of Shishikura and Lei which highlighted additional complexities that are absent in the quadratic setting, but arise for cubic rational maps \cite{shishikura-tan-mating}. 

In \cite{tan-newton-mating}, Lei described the dynamics of postcritically finite cubic Newton maps (these maps, which are obtained by plugging in complex polynomials in Newton's classical root-finding method, form an important and well-studied class of rational maps), and exhibited in the process the fact that a large subclass of such maps are matings. (See also the more recent work \cite{roesch-aspenberg} for a description of certain postcritically infinite cubic Newton maps as matings.)

The next theorem, due to Yampolsky and Zakeri, was the first existence result for conformal matings of polynomials that are not `close cousins' of postcritically finite ones. A quadratic polynomial $P$ is said to have a \emph{bounded type Siegel fixed point} if it has a fixed point $z_0$ with $P'(z_0)=e^{2\pi i\theta}$ such that the continued fraction expansion of $\theta\in\mathbb{R}/\mathbb{Z}$ has uniformly bounded partial fractions.

\begin{theorem}\cite{saeed-yampolski}
Suppose $P_1, P_2$ are quadratic polynomials which are not anti-holomorphically conjugate and each of which has a bounded type Siegel fixed point. Then $P_1$ and $P_2$ are conformally mateable.
\end{theorem}

The question of unmating a rational map; i.e., deciding whether a given rational map appears as the mating of two polynomials (and if so, whether such a decomposition is unique) has also been studied by various authors. For a general combinatorial characterization of hyperbolic, postcritically finite rational maps arising as matings, see \cite[Theorem 4.2]{meyer-unmating}. The situation is a bit more subtle for postcritically finite, non-hyperbolic rational maps, as discussed in the same paper. However, the next theorem gives a positive answer to the unmating question for a class of postcritically finite, non-hyperbolic rational maps:

\begin{theorem}\cite[Theorem~1.1]{meyer-09}
Let $R:\widehat{\C}\to\widehat{\C}$ be a postcritically finite rational map such that its Julia set is the whole sphere. Then every sufficiently high iterate $R^{\circ n}$ of $R$ arises as a mating (i.e., is topologically conjugate to the topological mating of two polynomials).
\end{theorem}

We refer the reader to the excellent survey article \cite{meyer-unmating} for more on this topic.

Matings of geometrically finite polynomials (i.e., a polynomial whose postcritical set intersects the Julia set in a finite set, or equivalently, if every critical point is either preperiodic, or attracted to an attracting or parabolic cycle) were studied by Ha{\"i}ssinsky and Lei using techniques of David homeomorphisms. They showed that two geometrically finite polynomials $P_1$ and $P_2$ with connected Julia sets and parabolic periodic points are mateable if and only if the postcritically finite polynomials $\mathcal{T}(P_1), \mathcal{T}(P_2)$ canonically associated to $P_1, P_2$ (such that $\mathcal{T}(P_i)$ and $P_i$ have topologically conjugate Julia set dynamics, $i=1,2$) are mateable \cite[Theorem~D]{haissinsky-tan} (cf. \cite[Theorem~5.2]{LMMN}).
\vspace{4mm}

\noindent\textbf{Mating of anti-holomorphic polynomials.}\ 
Let us now mention a class of anti-holomorphic polynomials (anti-polynomials for short) for which a complete solution to the conformal mating problem is known. These are the so-called \emph{critically fixed} anti-polynomials; i.e., anti-polynomials that fix all of their critical points. The proof of the following theorem crucially uses \cite[Theorem~3.2]{pilgrim-tan}, which in many situations, facilitates the application of Thurston's topological characterization of rational maps.

\begin{theorem}\cite[Theorem~1.3]{LLM1}\label{mating_crit_fixed_thm}
Let $P_1$ and $P_2$ be two (marked) anti-polynomials of equal degree $d\geq 2$, where $P_1$ is critically fixed and $P_2$ is postcritically finite, hyperbolic.
Then there is an anti-rational map $R$ that is the conformal mating of $P_1$ and $P_2$ if and only if there is no Moore obstruction.
\end{theorem}

In the opposite direction, the question of unmating critically fixed anti-rational maps was also settled in \cite[Theorem~1.2]{LLM1}, and examples of \emph{shared matings} were demonstrated (cf. \cite{rees-shared-mating}).

\begin{rmk}
1) Combined with \cite[Lemma~4.22]{LLM1}, Theorem~\ref{mating_crit_fixed_thm} yields an effective procedure to decide conformal mateability of a critically fixed anti-polynomial $P_1$ and a postcritically finite, hyperbolic anti-polynomial $P_2$. This is particularly useful in applying Theorem~\ref{antiholo_mating_thm} below to concrete examples.

2) It is worth mentioning that the above mating (respectively, unmating) results for critically fixed anti-polynomials (respectively, anti-rational maps) serve as a precise philosophical counterpart of the double limit theorem for (geometrically finite) Kleinian reflection groups in the complex dynamics world (see the discussion before \cite[Theorem~1.3]{LLM1} and \cite[\S 4.3]{LLM1}).
\end{rmk}

To conclude, we list a few relevant works that we did not touch upon in this survey: \cite{mashanova-timorin, epstein-sharland, cheritat, cui-peng-tan, buff-epstein-koch, sharland-cubic-1}.
A good part of the mating theory discussed above carries over to the setting of Thurston maps (i.e., postcritically finite, orientation-preserving branched coverings of $\mathbb{S}^2$), for which we encourage the reader to consult \cite{bonkmeyer, dudko-bartholdi-4}. Several beautiful visual illustrations of polynomial matings can be found in \cite{arnaud-movies}.
For a list of open questions on polynomial matings, we refer the reader to \cite{mating-questions}.

\section{Combining rational maps and Kleinian groups}\label{interbreeding_sec}

In this Section, we will expound recently developed frameworks for combining polynomials (respectively, anti-polynomials) with Kleinian (respectively, reflection) groups.

\subsection{Mating anti-polynomials with reflection groups}\label{antiholo_mating_subsec}

The story of mating anti-polynomials with Kleinian reflection groups began with the study of a new class of anti-holomorphic dynamical systems given by \emph{Schwarz reflection maps} associated with \emph{quadrature domains}. We will recall the definitions of these objects, and sketch the simplest examples of the mating phenomenon in Subsection~\ref{schwarz_subsec}. To put these examples in a general framework, we will introduce in Subsection~\ref{necklace_subsec} a class of Kleinian reflection groups (called \emph{necklace reflection groups}), that are central to the mating construction. Further, we will associate a map (called the \emph{Nielsen map}) to each necklace reflection group that is \emph{orbit equivalent} to the group. In Subsection~\ref{conformal_mating_subsec}, we formalize the notion of conformal mating of a necklace group and an anti-polynomial. Subsection~\ref{mating_examples_subsec} summarizes some of the main results of \cite{LLMM1,LLMM2,LMM2}, where various explicit examples of Schwarz reflection maps were shown to be conformal matings of necklace groups and anti-polynomials. Finally in Subsection~\ref{gen_thm_subsec}, we state a general combination theorem for necklace groups and anti-polynomials proved in \cite{LMMN}.

\subsubsection{Schwarz reflection maps and motivating examples}\label{schwarz_subsec}

By definition, a domain $\Omega\subsetneq\widehat{\C}$ satisfying $\infty \not\in \partial\Omega$ and $\Omega = \textrm{int}~{\overline{\Omega}}$ is a {\it quadrature domain} if there exists a continuous function $\sigma : \overline{\Omega} \to \widehat{\C}\ $ such that $\sigma$ is anti-meromorphic in $\Omega$ and $\sigma(z) = z$ on the boundary $\partial\Omega$. Such a function $\sigma$ is unique (if it exists), and is called the {\it Schwarz reflection map} associated with $\Omega$.
(See \cite{AS76}, \cite{lee-makarov} and the references therein.)

It is well known that except for a finite number of {\it singular} points (cusps and double points), the boundary of a quadrature domain consists of finitely many disjoint real analytic curves \cite{sakai-acta}. Every non-singular boundary point has a neighborhood where the local reflection in $\partial\Omega$ is well-defined. The (global) Schwarz reflection $\sigma$ is an anti-holomorphic continuation of all such local reflections.

Round disks on the Riemann sphere are the simplest examples of quadrature domains. Their Schwarz reflections are just the usual circle reflections. Further examples can be constructed using univalent polynomials or rational functions. In fact, simply connected quadrature domains admit a simple characterization.

\begin{prop}\cite[Theorem~1]{AS76}\label{simp_conn_quad}
	A simply connected domain $\Omega\subsetneq\widehat{\C}$ with $\infty\notin\partial\Omega$ and $\textrm{int}~{\overline{\Omega}}=\Omega$ is a quadrature domain if and only if the Riemann uniformization $f:\mathbb{D}\to\Omega$ extends to a rational map on $\widehat{\C}$. The Schwarz reflection map $\sigma$ of $\Omega$ is given by $f\circ(1/\overline{z})\circ(f\vert_{\mathbb{D}})^{-1}$. 
	
	In this case, if the degree of the rational map $f$ is $d$, then $\sigma:\sigma^{-1}(\Omega)\to\Omega$ is a (branched) covering of degree $(d-1)$, and $\sigma:\sigma^{-1}(\textrm{int}~{\Omega^c})\to\textrm{int}~{\Omega}^c$ is a (branched) covering of degree $d$.
\end{prop}

\[ \begin{tikzcd}
	\overline{\mathbb{D}} \arrow{r}{f} \arrow[swap]{d}{1/\overline{z}} & \overline{\Omega} \arrow{d}{\sigma} \\
	\widehat{\C}\setminus\mathbb{D} \arrow[swap]{r}{f}& \widehat{\C}.
\end{tikzcd}
\]

In \cite{lee-makarov}, questions on equilibrium states of certain $2$-dimensional Coulomb gas models were answered using iteration of Schwarz reflection maps associated with quadrature domains. It transpired from their work that these maps give rise to dynamical systems that are interesting in their own right. The general situation is as follows. Given a disjoint collection of quadrature domains, we call the complement of their union a \emph{droplet}. Removing the double points and cusps from the boundary of a droplet yields the \emph{desingularized droplet} or the \emph{fundamental tile}. One can then look at a partially defined anti-holomorphic dynamical system $\sigma$ that acts on (the closure of) each quadrature domain as its Schwarz reflection map. Under this dynamical system, the Riemann sphere $\widehat{\C}$ admits a dynamically invariant partition. The first one is an open set called the {\it escaping/tiling set}, it is the set of all points that eventually escape to the fundamental tile (on the interior of which $\sigma$ is not defined). The second invariant set is the {\it non-escaping} set, the complement of the tiling set or equivalently, the set of all points on which $\sigma$ can be iterated forever. When the tiling set contains no critical points of $\sigma$, it is often the case that the dynamics of $\sigma$ on its non-escaping set resembles that of an anti-polynomial on its filled Julia set, while the $\sigma-$action on the tiling set exhibits features of reflection groups.

\begin{figure}[ht!]
	\begin{center}
		\includegraphics[scale=0.32]{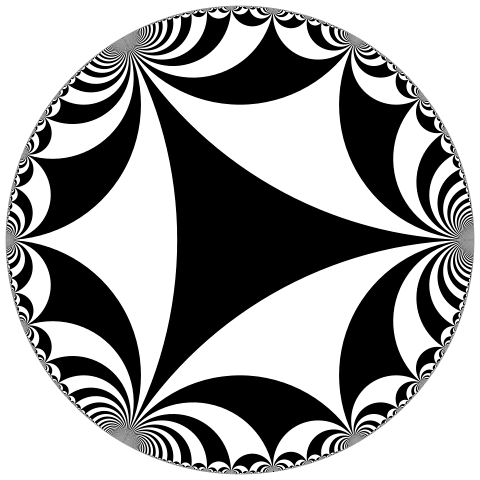}\hspace{4mm}\ \includegraphics[scale=0.35]{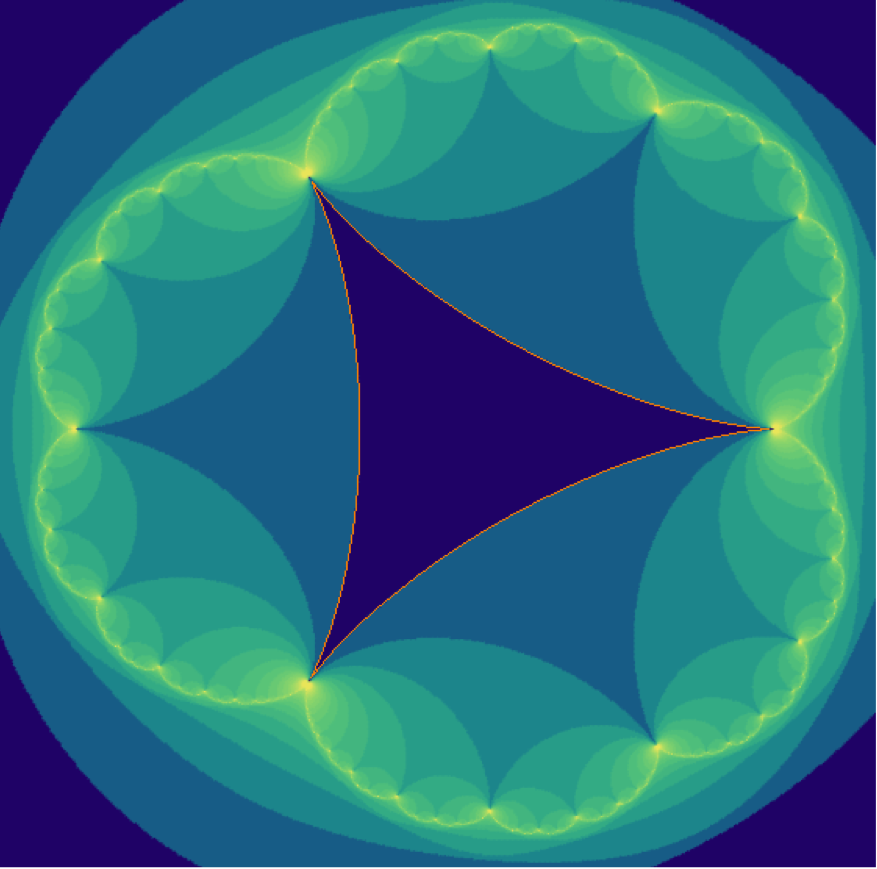}
	\end{center}
	\caption{Left: The tessellation of $\mathbb{D}$ for the ideal triangle reflection group. Right: The dynamical plane of the Schwarz reflection map associated with the quadrature domain $\Omega_0=f_0(\mathbb{D})$ (the exterior of a deltoid curve), where $f_0(z)=1/z+z^2/2$. The dynamics on the exterior of the bright green fractal curve is conformally conjugate to $\overline{z}^2\vert_{\mathbb{D}}$, while the dynamics on the interior is conformally equivalent to the Nielsen map of $\mathbf{\Gamma}_3$.}
	\label{tessellation_pic_1}
\end{figure}

This is precisely the case for the Schwarz reflection map of the exterior of a deltoid curve: this map is conformally conjugate to the anti-polynomial $\overline{z}^2$ on its non-escaping set, and conformally conjugate to a suitable piecewise circular reflection map associated with the ideal triangle reflection group on its tiling set. In this sense, this map is a conformal mating of $\overline{z}^2$ and the ideal triangle reflection group \cite[\S 5]{LLMM1} (see Theorem~\ref{deltoid_mating} and Figure~\ref{tessellation_pic_1}). 

\begin{figure}[ht!]
\centering
\includegraphics[scale=0.25]{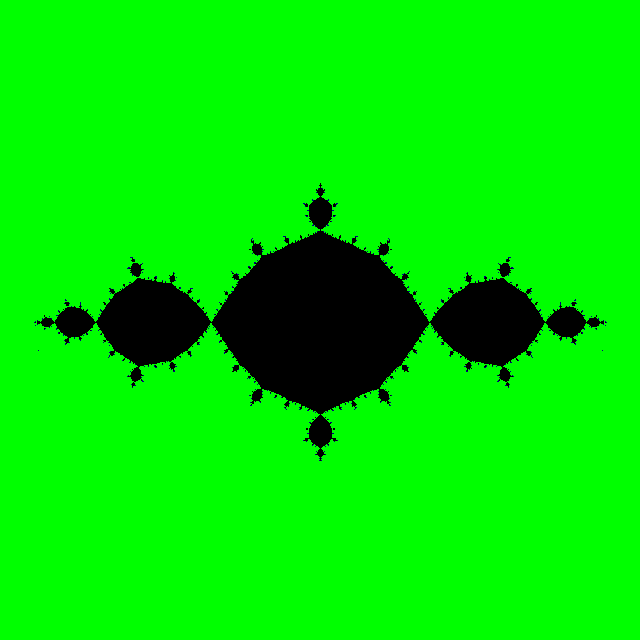}\hspace{2mm} \includegraphics[scale=0.184]{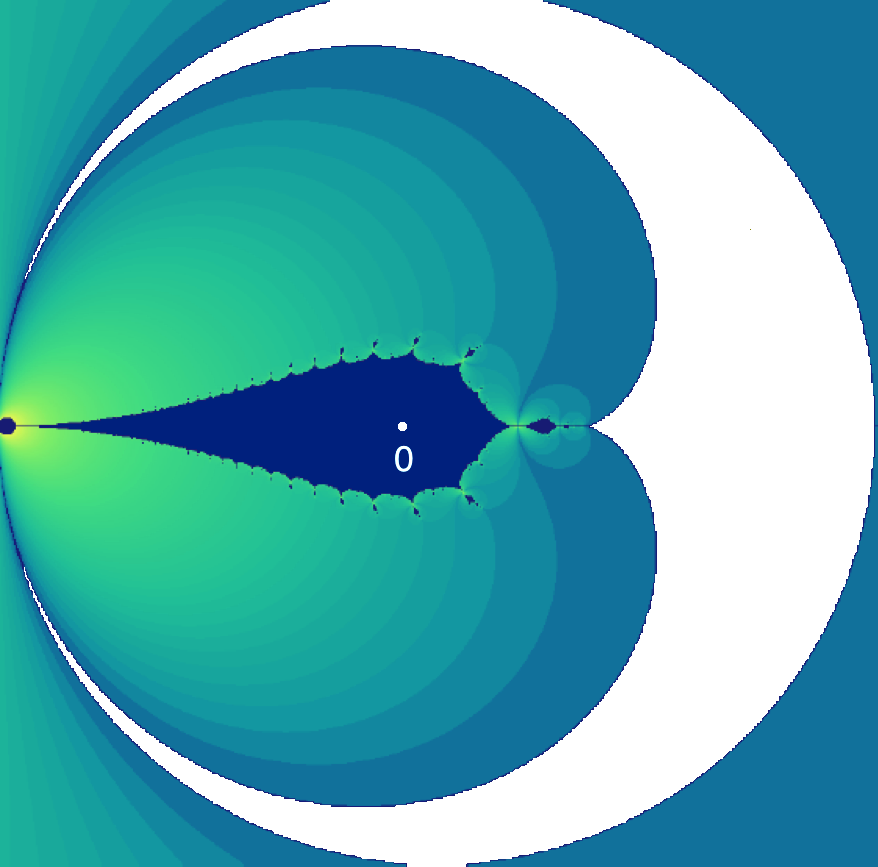}
\caption{Under the bijection $\chi$ of Theorem~\ref{c_and_c_mating}, the  postcritically finite quadratic anti-polynomial $\overline{z}^2-1$ corresponds to $F_a$ with $a=0$. Left: The filled Julia set of $\overline{z}^2-1$. Right:  The part of the non-escaping set of $F_0$ inside the cardioid (in dark blue) with the critical point $0$ marked. Both maps have a critical cycle of period $2$.}
\label{schwarz_basilica_fig}
\end{figure}

The above example was extended in \cite{LLMM2} by studying the Schwarz reflection maps associated with a fixed cardioid and a family of circumscribing circles. Such Schwarz reflection maps were shown to be conformal matings of generic quadratic anti-holomorphic polynomials with the ideal triangle reflection group (see Theorem~\ref{c_and_c_mating} for the precise statement and Figure~\ref{schwarz_basilica_fig} for a specific example). 

While the above examples produce matings of a rigid group with quadratic anti-polynomials, conformal matings of a large class of Kleinian reflection groups (called \emph{necklace groups}) with the anti-polynomial $\overline{z}^d$ were constructed in \cite{LMM2}, and these matings were realized as Schwarz reflection maps arising from a natural space of `univalent rational maps'.

\subsubsection{Necklace reflection groups}\label{necklace_subsec}

A circle packing is a connected collection of oriented circles in $\C$ with disjoint interiors (where the interior is determined by the orientation). Up to a M{\"o}bius map, we can always assume that no circle of the circle packing contains $\infty$ in its interior; i.e., the interior $\mathrm{int}~C$ of each circle $C$ in the circle packing can be assumed to be the bounded complementary component of $C$. Combinatorially, a circle packing can be described by its \emph{contact graph}, where we associate a vertex to each circle, and connect two vertices by an edge if and only if the two associated circles touch. By the Koebe-Andreev-Thurston circle packing theorem \cite[Corollary~13.6.2]{thurstonnotes}, every connected, simple, planar graph is the contact graph of some circle packing.

\begin{definition}\label{necklace_group} 
A \emph{necklace reflection group} is a group generated by reflections in the circles of a finite circle packing whose contact graph is \emph{$2$-connected} and \emph{outerplanar}; i.e., the contact graph remains connected if any vertex is deleted, and has a face containing all the vertices on its boundary. 
\end{definition}

Note that since a necklace reflection group is a discrete subgroup of the group of all M{\"o}bius and anti-M{\"o}bius automorphisms of $\widehat{\C}$, definitions of limit set and domain of discontinuity can be easily extended to necklace reflection groups. By \cite[Proposition~3.4]{LLM1}, the limit set of a necklace reflection group is connected. Moreover, for a necklace reflection group $\Gamma$ generated by reflections in the circles $C_1$, $\cdots$, $C_{d+1}$, the set
$$
\mathcal{F}_\Gamma:=\widehat{\C}\setminus\left(\bigcup_{i=1}^{d+1}\textrm{int}~{C_i}\bigcup_{j\neq k} (C_j \cap C_k)  \right)
$$
is a fundamental domain for the $\Gamma-$action on $\Omega(\Gamma)$ \cite[Proposition~7]{LMM2}.

\begin{figure}[ht!]
\begin{tikzpicture}[thick]
   \node at (2.5,3.6) [circle,fill=black,inner sep=2pt] {};
   \node at (2.5,0.4) [circle,fill=black,inner sep=2pt] {};
      \node at (2.5,2) [circle,fill=black,inner sep=2pt] {};
    \node at (0.6,2) [circle,fill=black,inner sep=2pt] {};
     \node at (4.4,2) [circle,fill=black,inner sep=2pt] {};

    \draw[-,line width=1pt] (2.5,3.6)->(2.5,0.4);
     \draw[-,line width=1pt] (2.5,3.6)->(0.6,2);
   \draw[-,line width=1pt] (2.5,3.6)->(4.4,2);
    \draw[-,line width=1pt] (2.5,0.4)->(0.6,2);
   \draw[-,line width=1pt] (2.5,0.4)->(4.4,2);

   \node at (9,3.6) [circle,fill=black,inner sep=2pt] {};
   \node at (7,3.6) [circle,fill=black,inner sep=2pt] {};
    \node at (6,2) [circle,fill=black,inner sep=2pt] {};
      \node at (10,2) [circle,fill=black,inner sep=2pt] {};
         \node at (8,0.4) [circle,fill=black,inner sep=2pt] {};
         
   \draw[-,line width=1pt] (10,2)->(9,3.6);
  \draw[-,line width=1pt] (9,3.6)->(7,3.6);
  \draw[-,line width=1pt] (7,3.6)->(6,2);
  \draw[-,line width=1pt] (6,2)->(8,0.4);
    \draw[-,line width=1pt] (8,0.4)->(10,2);
  \draw[-,line width=1pt] (8,0.4)->(7,3.6);
 \draw[-,line width=1pt] (8,0.4)->(9,3.6);
 
 \node at (-0.2,1) {};
     \end{tikzpicture}
   \caption{Left: A $2$-connected graph that is not outerplanar. Right: A $2$-connected, outerplanar graph. A circle packing realizing this graph and the limit set of the associated necklace reflection group are shown in Figure~\ref{sigma_mating_fig}(Right).}
   \label{outerplanar_fig}
   \end{figure}
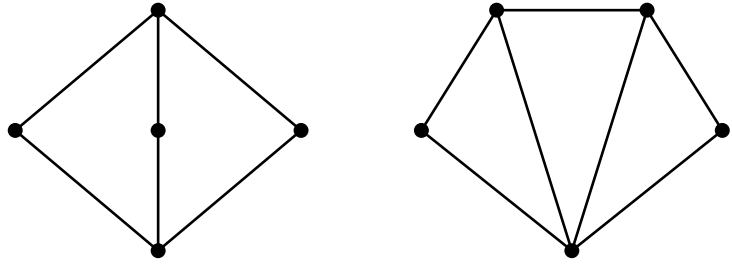

The domain of discontinuity of a necklace group has a simply connected, invariant component. Thus, such groups play the role of \emph{Bers slice closure} Kleinian groups in the world of reflection groups.

To a necklace reflection group $\Gamma$, one can associate a piecewise anti-M{\"o}bius reflection map $\rho_\Gamma$ that plays an important role in the mating construction.

\begin{definition}\cite[Definition~14]{LMM2}, \cite[Definition~6.6]{LMMN}\label{reflection_map} Let $\Gamma$ be a necklace reflection group generated by reflections $\{r_i\}_{i=1}^{d+1}$ in circles $\{C_i\}_{i=1}^{d+1}$. We define the associated \emph{Nielsen map} $\rho_\Gamma$ by: 
$$
\rho_{\Gamma} : \bigcup_{i=1}^{d+1} \overline{\textrm{int}~{C_i}} \rightarrow \widehat{\mathbb{C}},\quad z\longmapsto r_i(z)\ \textrm{  if }\ z \in \overline{\textrm{int}~{C_i}}. 
$$
\end{definition}

The next proposition underscores the intimate dynamical connection between a necklace group $\Gamma$ and its Nielsen map $\rho_\Gamma$.

\begin{prop}\cite[Proposition~16]{LMM2}\label{orbit_equiv_prop}
	Let $\Gamma$ be a necklace reflection group. The map $\rho_\Gamma$ is orbit equivalent to $\Gamma$ on $\widehat{\C}$; i.e., for any two points $z,w\in\widehat{\mathbb{C}}$, there exists $g\in\Gamma$ with $g(z)=w$ if and only if there exist non-negative integers $n_1, n_2$ such that $\rho_\Gamma^{\circ n_1}(z)=\rho_\Gamma^{\circ n_2}(w)$.
\end{prop}

The simplest examples of necklace reflection groups are regular ideal polygon reflection groups.

\begin{definition}\label{ideal_group} Consider the Euclidean circles $\mathbf{C}_1,\cdots, \mathbf{C}_{d+1}$ where $\mathbf{C}_j$ intersects $\mathbb{S}^1$ at right angles at the roots of unity $\exp{(\frac{2\pi i\cdot(j-1)}{d+1})}$, $\exp{(\frac{2\pi i\cdot j}{d+1})}$. (By \cite[Part~II, Chapter~5, Theorem~1.2]{VS93}, the group generated by reflections in these circles is discrete.) We denote this group by $\mathbf{\Gamma}_{d+1}$. 
\end{definition}

Note that the Nielsen map $\rho_{{\tiny{\mathbf{\Gamma}}}_{d+1}}$ of the regular ideal polygon reflection group $\mathbf{\Gamma}_{d+1}$ restricts to an expansive degree $d$ orientation-reversing covering of $\mathbb{S}^1$. By \cite{coven-reddy}, there exists a homeomorphism $\mathcal{E}_d$ of the circle that conjugates $\rho_{\mathbf{\Gamma}_{d+1}}$ to $\overline{z}^d$. The conjugacy $\mathcal{E}_d$ serves as a connecting link between reflection groups and quadratic anti-polynomials.

\subsubsection{Conformal mating of anti-polynomials and necklace groups}\label{conformal_mating_subsec} The precise meaning of conformal matings of the Nielsen map of a necklace group and an anti-polynomial is given below. The definition is an adaptation of the classical definition of conformal matings of two polynomials. 

Let $\Gamma$ be a necklace group generated by reflections in circles $C_1$, $\cdots$, $C_{d+1}$. The unbounded component of the domain of discontinuity $\Omega(\Gamma)$ is $\Gamma-$invariant \cite[Proposition~15]{LMM2}, and we denote it by $\Omega_\infty(\Gamma)$. We also set $\mathcal{K}(\Gamma):=\mathbb{C}\setminus\Omega_\infty(\Gamma)$. According to \cite[Proposition~22]{LMM2}, the restriction of $\rho_\Gamma$ to $\Omega_\infty(\Gamma)$ is conformally conjugate to the $\rho_{\mathbf{\Gamma}_{d+1}}-$action on $\widehat{\C}\setminus\overline{\mathbb{D}}$, and (the inverse of) this conformal conjugacy continuously extends to yield a semiconjugacy $\phi_\Gamma: \mathbb{S}^1 \rightarrow \Lambda(\Gamma)=\partial\mathcal{K}(\Gamma)$ between $\rho_{\mathbf{\Gamma}_{d+1}}\vert_{\mathbb{S}^1}$ and $\rho_\Gamma\vert_{\Lambda(\Gamma)}$ such that $\phi_\Gamma(1)$ is the point of tangential intersection of $C_1$ and $C_{d+1}$. Recall also that $\mathcal{E}_d: \mathbb{S}^1 \rightarrow \mathbb{S}^1$ is a topological conjugacy between $\rho_{\mathbf{\Gamma}_{d+1}}\vert_{\mathbb{S}^1}$ and $z\mapsto\overline{z}^{d}|_{\mathbb{S}^1}$. 

Let $P$ be a monic, centered, anti-polynomial of degree $d$ such that $\mathcal{J}(P)$ is connected and locally connected. Denote by $\phi_P: \mathbb{D}^*\rightarrow \mathcal{B}_\infty(P)$ the B\"ottcher coordinate for $P$ such that $\phi_P'(\infty)=1$. We note that since $\partial \mathcal{K}(P)=\mathcal{J}(P)$ is locally connected by assumption, it follows that $\phi_P$ extends to a semiconjugacy between $z\mapsto\overline{z}^{d}\vert_{\mathbb{S}^1}$ and $P\vert_{\mathcal{J}(P)}$. 

The equivalence relation below specifies a gluing of $\mathcal{K}(\Gamma)$ with $\mathcal{K}(P)$ along their boundaries. The presence of the topological conjugacy $\mathcal{E}_d$ in the definition of the equivalence relation ensures that the maps $\rho_\Gamma$ and $P$ fit together to produce a continuous map on the resulting topological $2$-sphere (when there is no Moore obstruction).

\begin{definition}\label{conf_mating_equiv_reltn} 
We define the equivalence relation $\sim$ on $\mathcal{K}(\Gamma) \sqcup \mathcal{K}(P)$ generated by $\phi_\Gamma(t)\sim\phi_P(\overline{\mathcal{E}_d(t)})$ for all $t\in\mathbb{S}^1$.
\end{definition}

The following definition essentially says that an anti-holomorphic map $F$ (defined on a subset of the Riemann sphere) is a conformal mating of $\Gamma$ and $P$ if there are continuous semi-conjugacies from $\mathcal{K}(\Gamma), \mathcal{K}(P)$ into the dynamical plane of $F$ (conformal on the interiors) such that the images fill up the whole sphere and intersect only along their boundaries as prescribed by the equivalence relation $\sim$ (compare Definition~\ref{mating_def_1}).

\begin{definition}\cite[Definition~10.16]{LMMN}\label{mating_def_2} 
Let $\Gamma$ be a necklace group as above, and let $P$ be a monic, centered anti-polynomial such that $\mathcal{J}(P)$ is connected and locally connected. Further, let $\Omega\subsetneq\widehat{\C}$ be an open set, and $F:\overline{\Omega}\to\widehat{\C}$ be a continuous map that is anti-meromorphic on $\Omega$. We say that $F$ is a \emph{conformal mating} of $\Gamma$ with $P$ if there exist continuous maps \[ \psi_P: \mathcal{K}(P) \rightarrow \widehat{\mathbb{C}} \textrm{ and } \psi_\Gamma: \mathcal{K}(\Gamma) \rightarrow \widehat{\C},  \] conformal on $\textrm{int}~{\mathcal{K}(P)}$, $\textrm{int}~{\mathcal{K}(\Gamma)}$, respectively, such that 

\begin{enumerate}\upshape
\item\label{whole_sphere} $\psi_P(\mathcal{K}(P))\bigcup\psi_\Gamma(\mathcal{K}(\Gamma))=\widehat{\C}$,

\item\label{domain_omega} $\Omega=\widehat{\C}\setminus\psi_\Gamma(\overline{\mathcal{F}_\Gamma})$,

\item\label{poly_semiconj} $\psi_P\circ P=F\circ\psi_P$ on $\mathcal{K}(P)$, 

\item\label{nielsen_semiconj} $\psi_\Gamma\circ \rho_\Gamma=F\circ\psi_\Gamma$ on $\mathcal{K}(\Gamma)\setminus \textrm{int}~{\mathcal{F}_\Gamma}$, and

\item\label{identifications} $\psi_\Gamma(z)=\psi_P(w)$ if and only if $z\sim w$ where $\sim$ is as in Definition~\ref{conf_mating_equiv_reltn}. 
\end{enumerate}
\end{definition}

\begin{rmk}
For the purposes of mating necklace groups with anti-polynomials, it is important to work with labeled circle packings, or equivalently, to regard the space of necklace groups as a space of representations of the ideal polygon reflection group $\mathbf{\Gamma}_{d+1}$. While we have suppressed this abstraction for ease of exposition, we refer the reader to \cite[\S 2.2]{LMM2} or \cite[\S 10.1]{LMMN}, where necklace groups are organized in Bers slices of $\mathbf{\Gamma}_{d+1}$. Although this point of view may seem like an artificial complication at a first glance, the language of representations turns out to be an unavoidable technicality in the mating theory. Roughly speaking, different representations give rise to different ways of gluing the limit set of a necklace group with the Julia set of an anti-polynomial, and the choice of gluing determines whether or not a conformal mating exists (compare \cite[Remark~10.21]{LMMN}). 
\end{rmk}

\subsubsection{Examples of the mating phenomenon}\label{mating_examples_subsec}
By studying the dynamics and parameter spaces of specific families of Schwarz reflection maps, one can often recognize such maps as matings of anti-polynomials and necklace reflection groups. This strategy was successfully implemented in \cite{LLMM1,LLMM2,LMM2}. We collect some results from these papers in this subsection.\\

\noindent\textbf{Example 1: The deltoid reflection.} We will start with the simplest instance of the mating phenomenon; namely, the conformal mating of the anti-polynomial $\overline{z}^2$ and the ideal triangle reflection group $\mathbf{\Gamma}_3$.

\begin{theorem}\cite[Theorem~1.1]{LLMM1}\label{deltoid_mating} 
 The map $f_0(z)=1/z+z^2/2$ is injective on $\overline{\mathbb{D}}$, and hence $\Omega_0:=f_0(\mathbb{D})$ is a simply connected quadrature domain. The associated Schwarz reflection map $\sigma_0$ is the unique conformal mating of $\overline{z}^2$ and $\mathbf{\Gamma}_3$.
\end{theorem}

\begin{rmk}
A \emph{welding homeomorphism} is a homeomorphism of the circle that arises as the composition of a conformal map from the unit disk onto the interior region of a Jordan curve with a conformal map from the exterior of this Jordan curve onto the exterior of the unit disk. A complex-analytic corollary of Theorem~\ref{deltoid_mating} is that the circle homeomorphism $\mathcal{E}_2$ is a welding homeomorphism. That the same is true for each $\mathcal{E}_d$ ($d\geq 2$) follows from a straightforward higher degree generalization of Theorem~\ref{deltoid_mating} worked out in \cite[Appendix~B]{LLMM3} (also compare Theorem~\ref{sigma_mating_thm} below). We refer the reader to \cite[Theorem~5.1]{LMMN} for a general conformal welding result for circle homeomorphisms conjugating suitable covering maps of the circle.
\end{rmk}

\noindent\textbf{Example 2: The circle and cardioid family.} To describe Schwarz reflection maps that are conformal matings of other quadratic anti-polynomials with the ideal triangle reflection group, we need to recall the Circle and Cardioid family which was introduced in \cite[\S 6]{LLMM1}. We consider the fixed cardioid
$$
\heartsuit:=\left\{w=z/2-z^2/4:~|z|<1\right\},
$$ 
and for each complex number $a\in\C\setminus\left(-\infty,-1/12\right)$, let $B(a,r_a)$ be the smallest open disk containing $\heartsuit$ centered at $a$ (in other words, $\{w: \vert w-a\vert=r_a\}$ is a circumcircle of the cardioid; see \cite[Figure~2]{LLMM2}). Let $\Omega_a := \heartsuit\ \cup\ \overline{B}(a,r_a)^c$ (where $\overline{B}(a,r_a)$ is the closed disk $\{w: \vert w-a\vert\leq r_a\}$), and $T_a:=\Omega_a^c$. We now define a piecewise Schwarz reflection dynamical system $F_a:\overline{\Omega}_a\to\widehat{\C}$ as, 
$$
w \mapsto \left\{\begin{array}{ll}
                    \sigma(w) & \mbox{if}\ w\in\overline{\heartsuit}, \\
                    \sigma_a(w) & \mbox{if}\ w\in B(a,r_a)^c, 
                                          \end{array}\right. 
$$
where $\sigma$ is the Schwarz reflection of $\heartsuit$, and $\sigma_a$ is reflection with respect to the circle $\partial B(a,r_a)$. The family
$$
\mathcal{S}:=\left\{F_a:\overline{\Omega}_a\to\widehat{\C}:a\in\C\setminus (-\infty,-1/12)\right\}
$$ 
is referred to as the C\&C family.

For any $a\in\C\setminus(-\infty,-1/12)$, $\partial T_a$ has two singular points; namely, the double point $\alpha_a$ where $\partial B(a,r_a)$ touches $\partial\heartsuit$, and the cusp point $\frac{1}{4}$. Both of them are fixed points of $F_a$. The \emph{fundamental tile} of $F_a$ is defined as $T_a^0:= T_a\setminus\{\alpha_a,\frac{1}{4}\}$. A parameter $a\in\C\setminus(-\infty,-1/12)$ (equivalently, the corresponding map $F_a\in\mathcal{S}$) is said to be \emph{postcritically finite} if the unique (simple) critical point $0$ of $F_a$ has a finite forward orbit that does not meet $T_a^0$. The following mating description for postcritically finite maps in $\mathcal{S}$ was given in \cite{LLMM2}.

\begin{theorem}\cite[\S 8]{LLMM1}, \cite[Theorems~1.1, 1.2]{LLMM2}\label{c_and_c_mating}
There exists a bijection $\chi$ between postcritically finite maps in $\mathcal{S}$ and (the M{\"o}bius conjugacy classes of) postcritically finite quadratic anti-polynomials $\overline{z}^2+c$ (excluding $\overline{z}^2$) such that the postcritically finite map $F_a\in\mathcal{S}$ is a conformal mating of the ideal triangle reflection group $\mathbf{\Gamma}_3$ and the quadratic anti-polynomial $\overline{z}^2+\chi(a)$.
\end{theorem}

\begin{rmk}
For conformal matings of $\mathbf{\Gamma}_3$ with more general quadratic anti-polynomials, see \cite[Theorem~1.1]{LLMM2}. Further, a combinatorial model of the \emph{connectedness locus} of the family $\mathcal{S}$ is given in \cite[Theorem~1.4]{LLMM2} in terms of the \emph{Tricorn}, which is the connectedness locus of quadratic anti-polynomials (for a quick introduction to the Tricorn, see \cite[\S 2]{LLMM2}, and for its detailed topological properties, see \cite{hub-sch,MNS,inou-muk-1,inou-muk-2}).
\end{rmk}

\noindent\textbf{Example 3: The space $\Sigma_d^*$.} The family of `univalent rational maps'
$$
\Sigma_d^*:= \left\{ g(z)= z+\frac{a_1}{z} + \cdots +\frac{a_d}{z^d} : a_d=-\frac{1}{d}\textrm{ and } g\vert_{\widehat{\C}\setminus\overline{\mathbb{D}}} \textrm{ is conformal}\right\}
$$
was introduced in \cite{lee-makarov-sharpness} and studied extensively in \cite{LMM1,LMM2} in terms of the associated Schwarz reflection maps.

\begin{rmk}
The family $\Sigma_d^*$ is closely related to the classically studied space $\Sigma$ of suitably normalized schlicht functions on $\widehat{\C}\setminus\overline{\mathbb{D}}$, see \cite[\S 4.7, \S 9.6]{duren}. 
\end{rmk}

Combining the pinching deformation theory for $\Sigma_d^*$ (developed in \cite{LMM1}) with tools from holomorphic dynamics, it was proved in  \cite{LMM2} that:

\begin{figure}[ht!]
\centering
\includegraphics[scale=0.23]{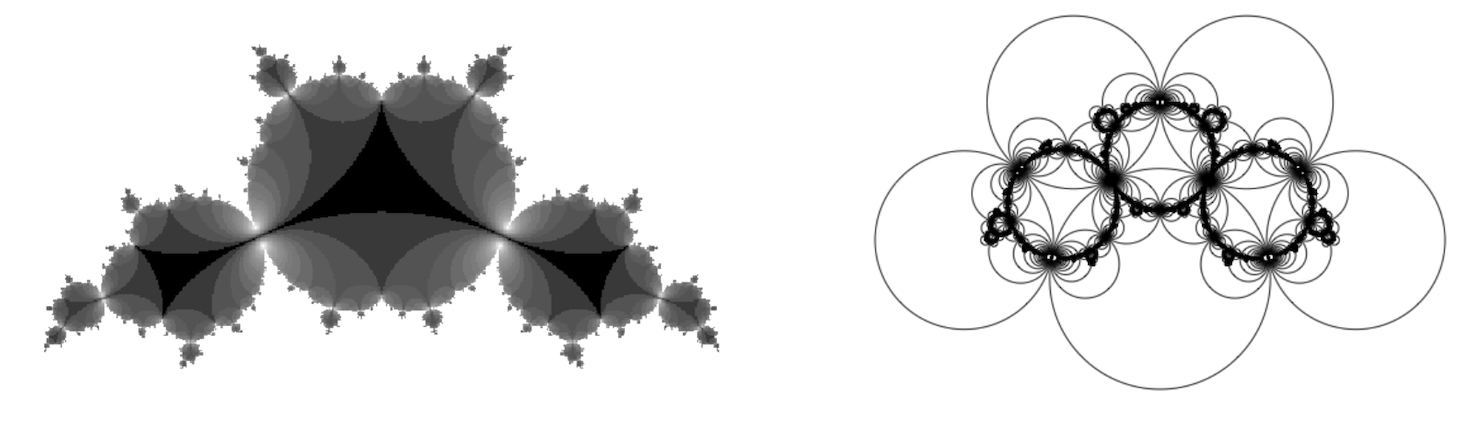}
\caption{Left: The dynamical plane of the Schwarz reflection map associated with some $f\in\partial\Sigma_4^*$. Right: The limit set of the corresponding necklace reflection group $\Gamma_f$.}
\label{sigma_mating_fig}
\end{figure}

\begin{theorem}\cite[Theorem~A]{LMM2}\label{sigma_mating_thm}
There is a bijection $f\mapsto\Gamma_f$ between $\Sigma_d^*$ and the space of necklace reflection groups of rank $d+1$ (up to a natural equivalence) such that the Schwarz reflection map associated with $f\in\Sigma_d^*$ is a conformal mating of the anti-polynomial $\overline{z}^d$ with the corresponding necklace group $\Gamma_f$.
\end{theorem}

\subsubsection{The general theorem}\label{gen_thm_subsec}

We conclude our discussion of combinations of necklace reflection groups and anti-polynomials with a general existence theorem:

\begin{theorem}\cite[Lemma~10.17, Theorem~10.20]{LMMN}\label{antiholo_mating_thm}
	Let $P$ be a monic, postcritically finite, hyperbolic anti-polynomial of degree $d$, and let $\Gamma$ be a necklace group. Then, $P$ and $\Gamma$ are conformally mateable if and only if $\mathcal{K}(P)\sqcup\mathcal{K}(\Gamma)/\sim$ is homeomorphic to $\mathbb{S}^2$ (where $\sim$ is the equivalence relation from Definition~\ref{conf_mating_equiv_reltn}). 
	
	Moreover, if $F:\overline{\Omega}\to\widehat{\C}$ is a conformal mating of $\Gamma$ and $P$, then each component of $\Omega$ is a simply connected quadrature domain, and $F$ is the piecewise defined Schwarz reflection map associated with these quadrature domains.
\end{theorem}

  \begin{figure}[ht!]
\begin{tikzpicture}
\node[anchor=south west,inner sep=0] at (0,-4.5) {\includegraphics[width=0.42\textwidth]{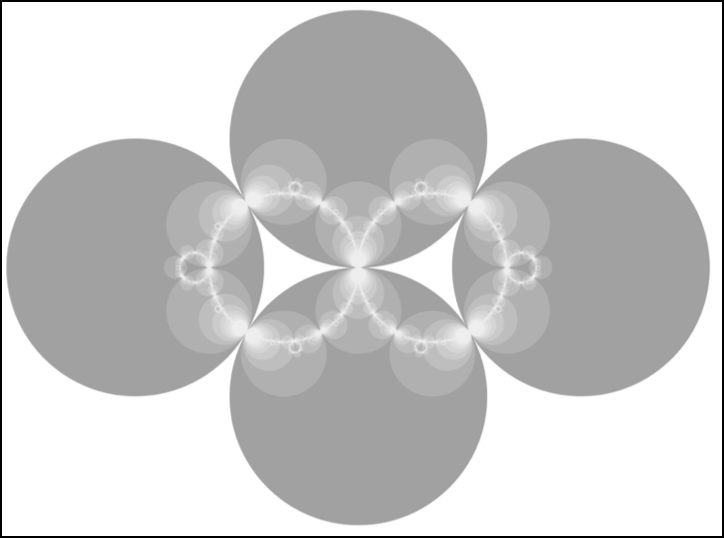}};
\node[anchor=south west,inner sep=0] at (5.8,-4.5) {\includegraphics[width=0.54\textwidth]{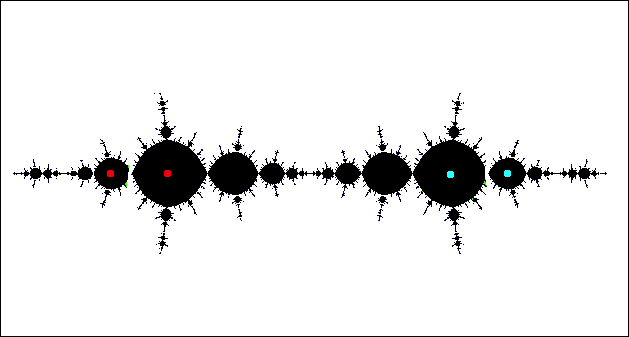}};
 \node[anchor=south west,inner sep=0] at (0,0) {\includegraphics[width=0.96\textwidth]{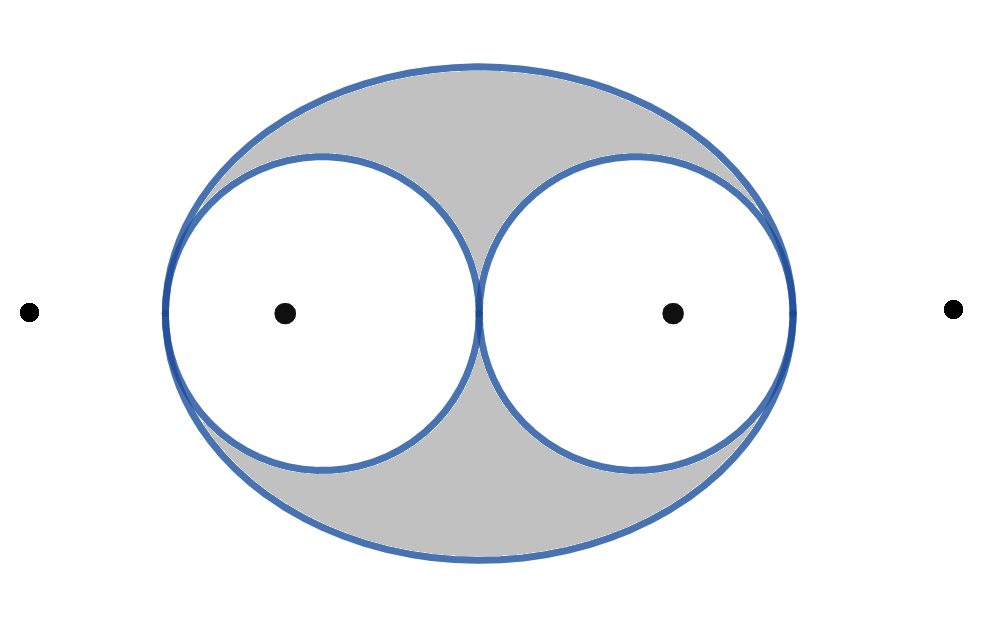}};
\draw [->] (8.48,4) to [out=35,in=145] (11.6,4);
\draw [->] (11.55,3.7) to [out=215,in=325] (8.48,3.66);
\draw [->] (0.5,4) to [out=35,in=145] (3.5,4);
\draw [->] (3.42,3.6) to [out=215,in=325] (0.6,3.7);
\draw[-,line width=0.5pt] (0,0)->(12.2,0);
\draw[-,line width=0.5pt] (12.2,0)->(12.2,7.64);
\draw[-,line width=0.5pt] (12.2,7.64)->(0,7.64);
\draw[-,line width=0.5pt] (0,7.64)->(0,0);
\node at (4.9,-1.5) {$C_4$};
\node at (3.8,-1) {$C_1$};
\node at (0.4,-1.4) {$C_2$};
\node at (3.8,-4) {$C_3$};
\end{tikzpicture}
\caption{Bottom left: The circles $C_{i}$ generate a necklace reflection group $\Gamma$. Bottom right: The dynamical plane of $P(z)=\overline{z}^3-\frac{3i}{\sqrt{2}}\overline{z}$; each critical point of which forms a $2$-cycle. Top: The conformal mating of $P$ and the necklace group $\Gamma$ is given by the piecewise Schwarz reflection map associated with the disjoint union of three quadrature domains: the exterior of an ellipse, and two round disks contained in the interior of the ellipse. Each of the two critical points of $F$ forms a $2$-cycle. (See \cite[\S 11.2]{LMMN} for proofs of these statements.)}
\label{ellipse_disk}
\end{figure}
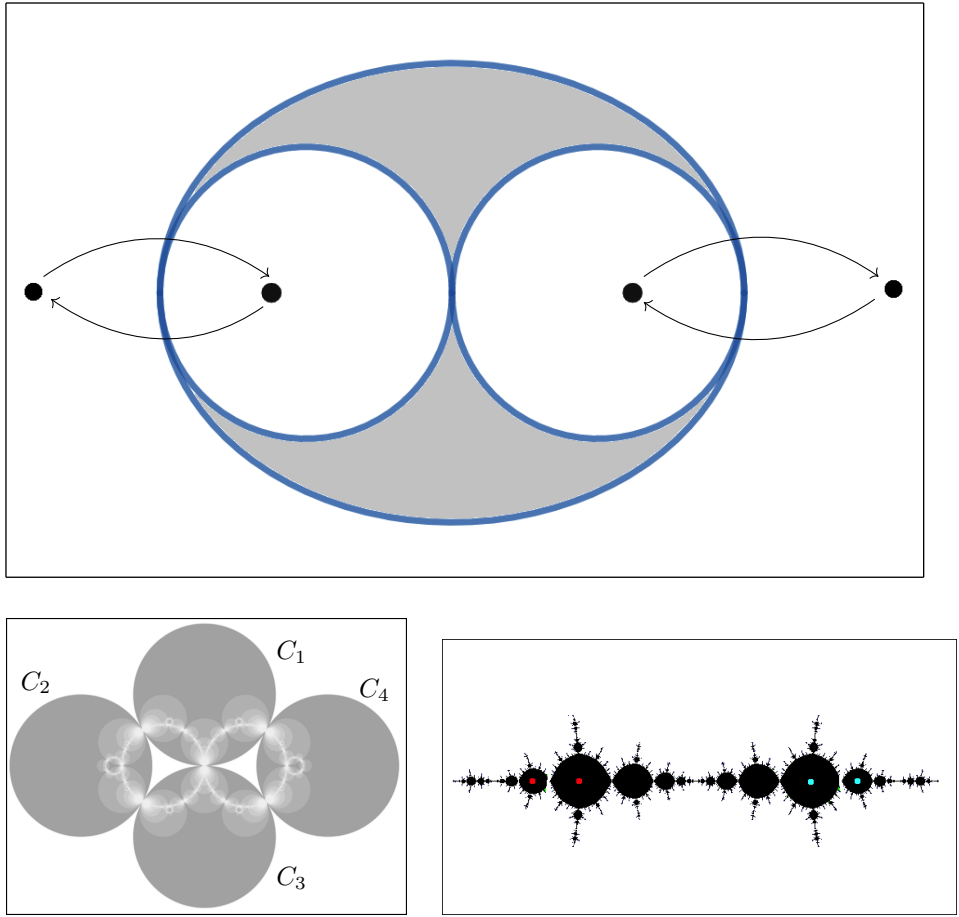

The hard part of the above theorem is to show that if $\mathcal{K}(P)\sqcup\mathcal{K}(\Gamma)/\sim$ is homeomorphic to $\mathbb{S}^2$, then a conformal mating of $P$ and $\Gamma$ exists. In fact, the condition that $\mathcal{K}(P)\sqcup\mathcal{K}(\Gamma)/\sim\ \cong\mathbb{S}^2$ guarantees the existence of a topological mating on a $2$-sphere, but promoting the topological mating to an anti-holomorphic map lies at the heart of the difficulty. This goal is achieved in two steps. One first uses Thurston's topological characterization theorem to construct a hyperbolic anti-rational map $R$ that is a conformal mating of $P$ and another postcritically finite (in fact, critically fixed), hyperbolic anti-polynomial $P_\Gamma$ such that the Julia dynamics of $P_\Gamma$ is topologically conjugate to the limit set dynamics of the Nielsen map $\rho_\Gamma$. The existence of such an anti-polynomial $P_\Gamma$ follows from \cite{LMM2} or \cite{LLM1}, while conformal mateability of $P$ and $P_\Gamma$ follows from the general mateability criterion given in Theorem~\ref{mating_crit_fixed_thm} (in fact, the condition $\mathcal{K}(P)\sqcup\mathcal{K}(\Gamma)/\sim\ \cong\mathbb{S}^2$ is equivalent to saying that the topological mating of $P$ and $P_\Gamma$ is Moore-unobstructed, so Theorem~\ref{mating_crit_fixed_thm} can be applied to produce $R$). Finally, to turn $R$ into a conformal mating of $P$ and $\Gamma$, one needs to glue Nielsen maps of ideal polygon reflection groups in suitable invariant Fatou components of $R$. The fact that all fixed points of $R$ on its Julia set are hyperbolic while those of a Nielsen map are parabolic prohibits the use of purely quasiconformal tools to carry out this task. This problem is tackled by employing surgery techniques involving David homeomorphisms: generalizations of quasiconformal homeomorphisms.

\begin{rmk}
Although Theorem~\ref{antiholo_mating_thm} guarantees the existence of conformal matings of suitable anti-polynomials and necklace reflection groups, in general, it may be hard to find explicit Schwarz reflection maps realizing such conformal matings. However, in certain low complexity situations, the second statement of the theorem (that the conformal matings are piecewise Schwarz reflection maps associated with simply connected quadrature domains) allows one to use Proposition~\ref{simp_conn_quad} and the desired dynamical properties to explicitly characterize the conformal matings (see Figure~\ref{ellipse_disk} for an illustration, and \cite[\S 11]{LMMN} for various worked out examples). 
\end{rmk}

\subsection{Mating polynomials with Kleinian groups}\label{holo_mating_subsec}

This subsection is a summary of \cite{mj-muk}, where a new setup for combination theorems of complex polynomials and Kleinian surface groups was designed using the notion of orbit equivalence.

\subsubsection{The Fuchsian case} A foundational problem that arises in trying to make sense of what it means to combine a polynomial $P$ with a Kleinian group $\Gamma$ is that on one side of the picture we have the semigroup $\langle P \rangle$ generated by $P$, while on the other side we have a non-commutative group
$\Gamma$ generated by more than one element. To formulate a precise notion of mateability between Fuchsian groups and complex polynomials (with Jordan curve Julia sets), one needs to address this inherent discord between these two objects, and this leads to the notion of \emph{mateable} circle maps: \emph{single} maps $A$  that capture essential  dynamical and combinatorial features of Fuchsian groups acting on ${\mathbb S}^1$. Further, $A$ should also be dynamically compatible with polynomial maps.
Before giving a precise definition of mateable maps, let us outline the underlying motivation: the following features are required of a mateable map $A:\mathbb{S}^1\to\mathbb{S}^1$.
\begin{enumerate}	
	\item $A$ must be dynamically compatible with a Fuchsian group $\Gamma$. This leads to \begin{enumerate}
	\item orbit equivalence between $A$ and $\Gamma$.
	\item $A$ has to be piecewise Fuchsian.
	\end{enumerate}  

        \item $A$ must be dynamically compatible with complex polynomials. Hence we demand the existence of a topological conjugacy between $A$ and the polynomial $z^d\vert_{\mathbb{S}^1}$ (where $d\geq 2$ is the degree of $A$),  

	\item $A$ must be combinatorially compatible with $z^d$ leading to a Markov condition, 
	
	\item $A$ must be conformally compatible with $z^d$ requiring absence of asymmetrically hyperbolic periodic break-points of $A$ (this is a weaker version of the $C^1$-condition).
\end{enumerate}

To fulfill the above requirements, we make the following definition. We denote the group of conformal automorphisms of the unit disk $\mathbb{D}$ by $\mathrm{Aut}({\mathbb{D}})$.

\begin{definition}\cite[Definitions~2.7,~2.16]{mj-muk}\label{pwfm_def}
\begin{enumerate}
\item A map $A:\mathbb{S}^1\to\mathbb{S}^1$ is called \emph{piecewise M{\"o}bius} if there exist $k\in\mathbb{N}$, closed arcs $I_j\subset\mathbb{S}^1$, and $g_j\in\mathrm{Aut}({\mathbb{D}})$, $j\in\{1,\cdots, k\}$,  such that
	\begin{enumerate}
		\item $\displaystyle\mathbb{S}^1=\bigcup_{j=1}^k I_j,$
		
		\item $\textrm{int}~{I_m}\cap\textrm{int}~{I_n}=\emptyset$ for $m\neq n$, and
		
		\item $A\vert_{I_j}=g_j$.
	\end{enumerate}
	
A piecewise M{\"o}bius map is called \emph{piecewise Fuchsian} if $g_1$, $\cdots$, $g_k$ generate a Fuchsian group, which we denote by $\Gamma_A$.
	
\item A map $A:\mathbb{S}^1\to\mathbb{S}^1$ is called \emph{piecewise Fuchsian Markov} if it is a piecewise Fuchsian expansive covering map (of degree at least two) such that the pieces (intervals of definition) of $A$ form a Markov partition for $A:\mathbb{S}^1\to\mathbb{S}^1$.

\item A piecewise Fuchsian Markov map $A$ is said to be \emph{mateable} if $A$ is orbit equivalent to the Fuchsian group $\Gamma_A$ generated by its pieces, and none of the periodic break-points of $A$ is asymmetrically hyperbolic.
\end{enumerate} 
\end{definition}

\noindent We refer the reader to \cite[\S 2]{mj-muk} for the definition of the term `symmetrically hyperbolic', and for a detailed discussion on the necessity of each of the requirements in the definition of a mateable map. We also note that the expansivity condition above ensures that a mateable map is topologically conjugate to the polynomial $z^d$ (for some $d\geq 2$).

\begin{rmk}
In the anti-holomorphic setting, the role of mateable maps was played by Nielsen maps of necklace reflection groups (see Subsection~\ref{necklace_subsec}). 
\end{rmk}

The simplest example of a mateable map is given by the classical \emph{Bowen--Series map} \cite{Bowen,bowen-series}. While such a map can be defined for arbitrary Fuchsian groups equipped with suitable fundamental domains, they are typically discontinuous. However, it turns out that for Fuchsian groups uniformizing spheres with punctures (possibly with one/two order two orbifold points), the Bowen--Series map is a covering map of the circle satisfying the defining properties of a mateable map \cite[\S 2]{mj-muk}. \\

\begin{figure}[ht!]
\begin{tikzpicture}
\node[anchor=south west,inner sep=0] at (0,0) {\includegraphics[width=0.46\linewidth]{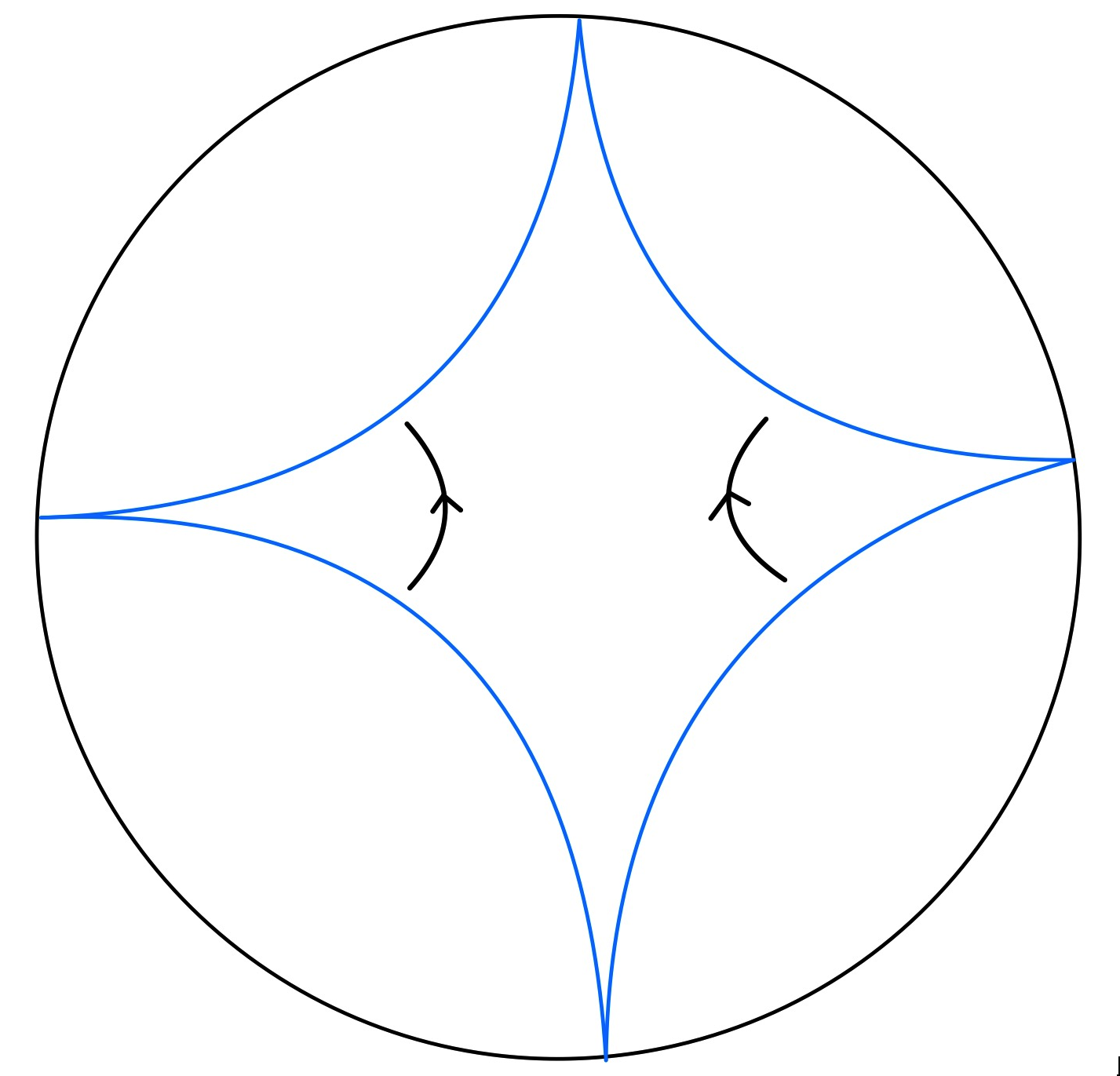}};
\node[anchor=south west,inner sep=0] at (6,0) {\includegraphics[width=0.49\linewidth]{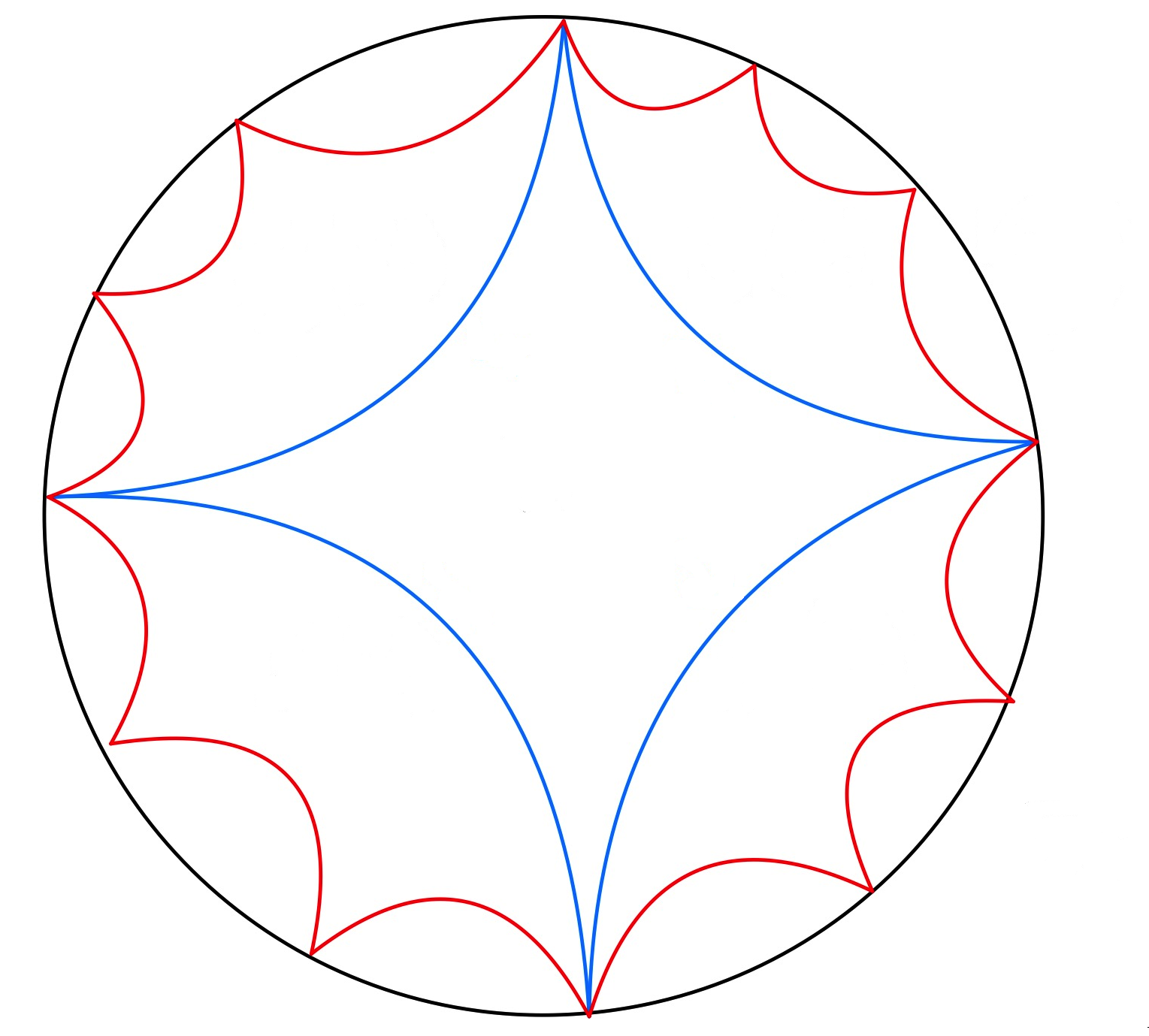}};
\node at (3,3.6) {\begin{Large}$R$\end{Large}};
\node at (2.58,2.7) {$g_1$};
\node at (3.56,2.7) {$g_2$};
\node at (4.4,4.2) {$g_2^{-1}$};
\node at (1.6,4.2) {$g_1^{-1}$};
\node at (4.4,1.6) {$g_2$};
\node at (1.6,1.6) {$g_1$};
\node at (9,3.2) {\begin{Large}$R$\end{Large}};
\node at (7.7,2) {$g_1^{-1}(R)$};
\node at (10.5,2.1) {$g_2^{-1}(R)$};
\node at (7.8,4) {$g_1(R)$};
\node at (10.3,4.1) {$g_2(R)$};
\node at (6.56,2.32) {\begin{tiny}$h_1$\end{tiny}};
\node at (7.2,1.2) {\begin{tiny}$h_2$\end{tiny}};
\node at (8.4,0.45) {\begin{tiny}$h_3$\end{tiny}};
\node at (10,0.64) {\begin{tiny}$h_3^{-1}$\end{tiny}};
\node at (11.72,1.2) {\begin{tiny}$h_3^{-1}h_1$\end{tiny}};
\node at (12.1,2.5) {\begin{tiny}$h_3^{-1}h_2$\end{tiny}};
\node at (12,3.84) {\begin{tiny}$h_2^{-1}h_3$\end{tiny}};
\node at (10.7,5.19) {\begin{tiny}$h_2^{-1}$\end{tiny}};
\node at (9.6,5.7) {\begin{tiny}$h_2^{-1}h_1$\end{tiny}};
\node at (7.7,5.6) {\begin{tiny}$h_1^{-1}h_2$\end{tiny}};
\node at (6.54,4.8) {\begin{tiny}$h_1^{-1}h_3$\end{tiny}};
\node at (6.1,3.6) {\begin{tiny}$h_1^{-1}$\end{tiny}};
\end{tikzpicture}
\caption{Left: $R$ is a fundamental domain of a Fuchsian group $\Gamma$ uniformizing $S_{0,3}$ with side pairing transformations $g_1^{\pm 1}, g_2^{\pm 1}$. The canonical extension of the Bowen--Series map of $\Gamma$ (associated with $R$) is defined on $\overline{\mathbb{D}}\setminus\textrm{int}~{R}$ in terms of $g_1^{\pm 1}, g_2^{\pm 1}$ as shown in the figure. Right: The second iterate of the Bowen--Series map of $\Gamma$ is a higher Bowen--Series map of the index two subgroup $\Gamma'=\langle g_1^2, g_1g_2, g_1g_2^{-1}\rangle \leq \Gamma$, which uniformizes $S_{0,4}$. Its canonical extension is defined on the region enclosed by $\mathbb{S}^1$ and the red (hyperbolic) geodesics in terms of $h_1=g_1^2, h_2=g_1g_2, h_3=g_1g_2^{-1}$ as shown in the figure. It maps the boundary of the red polygon onto the boundary of $R$. The degree of the Bowen--Series map of $\Gamma$ (as a circle covering) is $3$, so the degree of the higher Bowen--Series map of $\Gamma'$ is $9$.}
\label{bs_hbs_fig_1}
\end{figure}
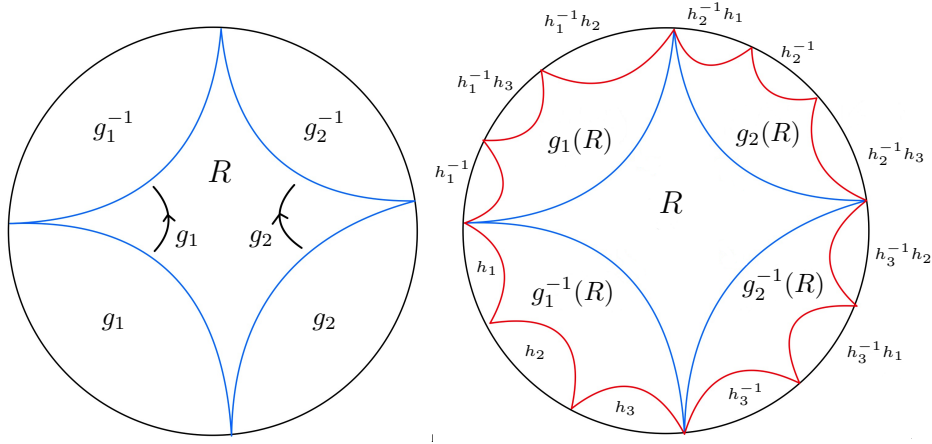

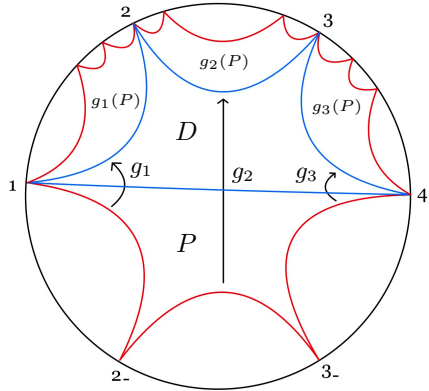
\begin{figure}[ht!]
\begin{tikzpicture}
\node[anchor=south west,inner sep=0] at (0,0) {\includegraphics[width=0.46\linewidth]{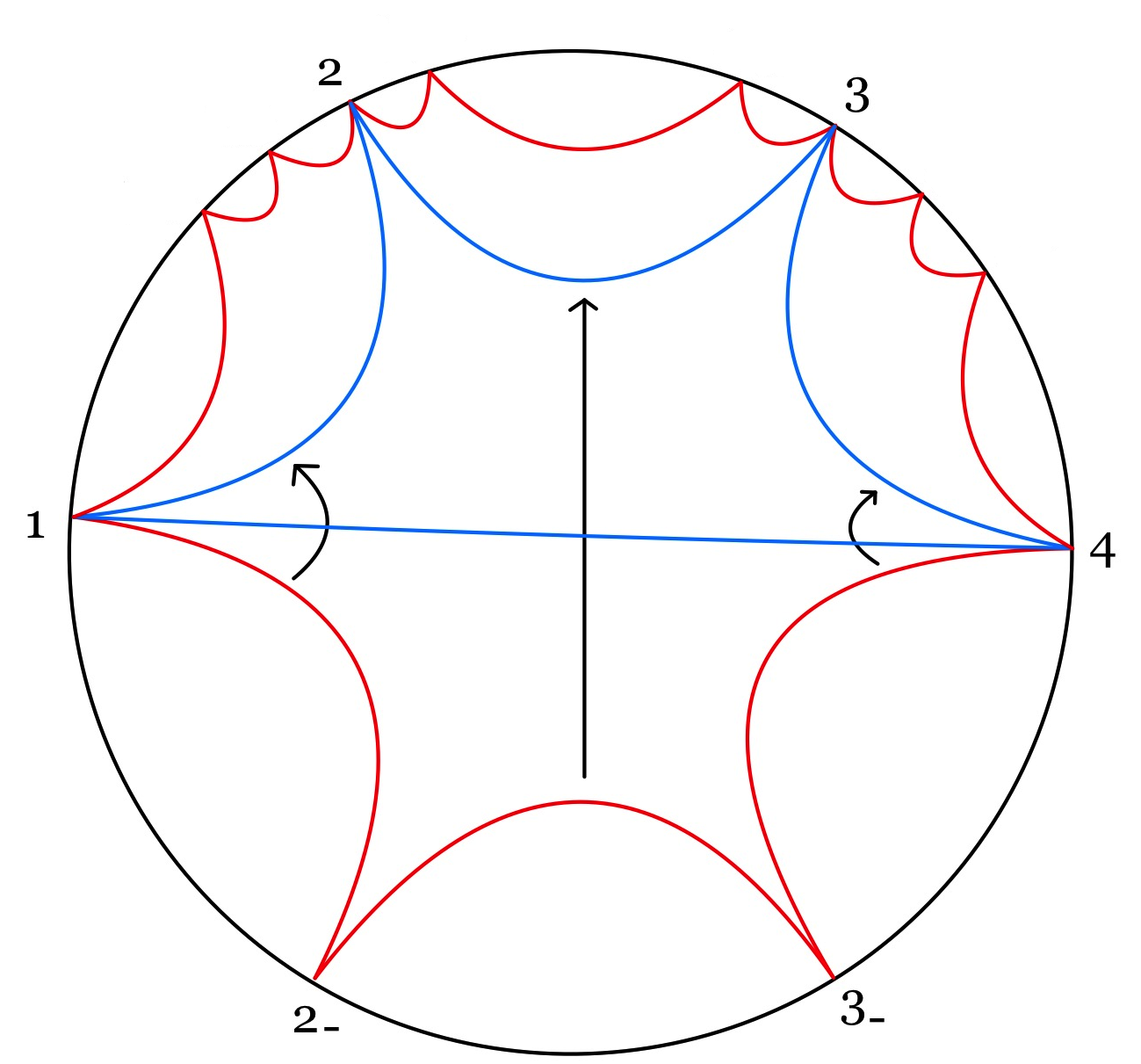}};
\node at (4.1,2.9) {\begin{small}$g_3$\end{small}};
\node at (3.25,2.9) {\begin{small}$g_2$\end{small}};
\node at (1.9,2.96) {\begin{small}$g_1$\end{small}};
\node at (2.5,2) {$P$};
\node at (2.5,3.48) {$D$};	
\node at (1.56,3.88) {\begin{tiny}$g_1(P)$\end{tiny}};	
\node at (3,4.4) {\begin{tiny}$g_2(P)$\end{tiny}};	
\node at (4.5,3.8) {\begin{tiny}$g_3(P)$\end{tiny}};		
\end{tikzpicture}\caption{The quadrilaterals $D$ and $P$ with ideal vertices at $1,2,3,4$ and $1,2_-,3_-,4$ (respectively) together form a fundamental domain of a Fuchsian group $\Gamma$ uniformizing $S_{0,4}$ with side pairing transformations $g_1^{\pm 1}, g_2^{\pm 1}, g_3^{\pm 1}$. The corresponding higher Bowen--Series map $A$ of $\Gamma$ acts on the anti-clockwise arc from $1$ to $4$ as the Bowen--Series map of $\Gamma$ associated with the fundamental domain $D\cup P$; while on the clockwise arc from $j$ to $j+1$, $A$ equals the Bowen--Series map of $\Gamma$ associated with the fundamental domain $D\cup g_j(P)$ ($j=1,2,3$). The degree of a Bowen--Series map of $\Gamma$ (associated with any fundamental domain) is $5$, while the degree of a higher Bowen--Series map of $\Gamma$ is $9$.}
\label{bs_hbs_fig_2}
\end{figure}

\noindent \textbf{Higher Bowen--Series maps.} More examples of mateable maps are given by \emph{higher Bowen--Series maps} of punctured sphere Fuchsian groups (see \cite[\S 4]{mj-muk} for their definition and basic properties). As suggested by the name, there are close connections between higher Bowen--Series maps and Bowen--Series maps. Indeed, the higher Bowen--Series map of a Fuchsian group uniformizing $S_{0,k}$ (a sphere with $k$ punctures) can be represented as the second iterate of the Bowen--Series map of a Fuchsian group uniformizing a sphere with roughly k/2 punctures and zero/one/two order two orbifold points \cite[Corollary~5.6]{mj-muk} (see Figure~\ref{bs_hbs_fig_1}). Alternatively, a higher Bowen--Series map of a Fuchsian group is obtained by `gluing together' several Bowen--Series maps of the same Fuchsian group with overlapping fundamental domains \cite[Proposition~4.5]{mj-muk} (see Figure~\ref{bs_hbs_fig_2}).

Every piecewise Fuchsian Markov map $A$ of the circle can be conformally extended to a canonically defined subset of $\overline{\mathbb{D}}$ (see \cite[\S 2.2]{mj-muk}). This extension is termed the \emph{canonical extension of $A$}. The following result, which is a conformal combination theorem for punctured sphere Fuchsian groups and hyperbolic polynomials with Jordan curve Julia sets, can be regarded as an analog of the Bers simultaneous uniformization theorem in the current setting.

\begin{theorem}\cite[Theorem~3.7, Theorem~4.8]{mj-muk}\label{holo_mating_thm}
	The canonical extensions of Bowen--Series maps and higher Bowen--Series maps of Fuchsian groups uniformizing punctured spheres (possibly with one/two orbifold points of order two) can be conformally mated with polynomials lying in principal hyperbolic components (of appropriate degree).
\end{theorem}

\begin{rmk}
As in the anti-holomorphic case, there is a key qualitative difference between the dynamics of Bowen--Series (respectively, higher Bowen--Series) maps on $\mathbb{S}^1$ and the dynamics of polynomials (lying in principal hyperbolic components) on their Julia set; namely, the former has parabolic fixed points on $\mathbb{S}^1$ while all fixed points of the latter on their Julia sets are repelling. Consequently, the topological conjugacy between a Bowen--Series (respectively, higher Bowen--Series) map and such a polynomial is not quasisymmetric. This forces one to abandon classical quasiconformal techniques (used in the proof of the Bers simultaneous uniformization theorem), and apply David homeomorphisms to prove Theorem~\ref{holo_mating_thm}.
\end{rmk}

\noindent \textbf{Moduli space of Fuchsian matings.} In the torsion-free case, the only topological surfaces that Theorem~\ref{holo_mating_thm} succeeds to combine with complex polynomials are punctured spheres (see \cite[6.35]{mj-muk} for the definition of moduli space of matings between a topological surface and hyperbolic complex polynomials with Jordan curve Julia sets). This naturally raises the following questions. 

\begin{enumerate}
\item Do mateable maps exist for higher genus surfaces (possibly with punctures)? More precisely, does there exist a mateable map $A$ with $\mathbb{D}/\Gamma_A\cong S_{g,k}$, for $g\geq 1$?

\item Are Bowen--Series and higher Bowen--Series maps the only mateable maps associated with punctured spheres?

\end{enumerate}

In this generality, the above questions remain open. However, \cite[Theorems~6.18, 6.33]{mj-muk} give a complete description of mateable maps satisfying some natural $2$-point  conditions over and above orbit equivalence. It turns out that under these additional hypotheses, punctured spheres are the only topological surfaces that can be combined with complex polynomials (see \cite[Theorem~6.36]{mj-muk} for a complete description of the interiors of such constrained moduli space of matings). A major part of the proofs of these theorems is to determine the topology of the surface $\mathbb{D}/\Gamma_A$ from the dynamical properties of a mateable map $A$, and this is accomplished by analyzing certain patterns and laminations associated with mateable maps.

\subsubsection{The case of Bers boundary groups}
We proceed to discuss the structure of the boundaries of the moduli spaces of Fuchsian matings arising from Theorem~\ref{holo_mating_thm}. 

For definiteness, let us fix a base Fuchsian group $\Gamma_0$ uniformizing $S_{0,k}$ ($k\geq 3$). For each $\Gamma$ lying on the boundary of the Bers slice $\mathcal{B}(\Gamma_0)$, there exists a continuous map $\phi_\Gamma:\mathbb{S}^1\to\Lambda_\Gamma$, called the \emph{Cannon--Thurston map} after \cite{CTpub}, that semi-conjugates the action of $\Gamma_0$ to that of $\Gamma$ \cite{mahan-split,mahan-elct,mahan-red,mahan-kl}. In fact, the data of the \emph{ending lamination} can be recovered from the Cannon--Thurston map. More precisely, the group $\Gamma$ can be obtained by `pinching a  lamination' on the surface $\mathbb{D}/\Gamma_0$ (while keeping the hyperbolic structure on the surface $(\widehat{\C}\setminus\overline{\mathbb{D}})/\Gamma_0$ unchanged), and the endpoints of the corresponding geodesic lamination on $\mathbb{D}$ generate a $\Gamma_0-$invariant equivalence relation on $\mathbb{S}^1$ which agrees with the one defined by the fibers of the Cannon--Thurston map $\phi_\Gamma$.  

Assume further that $A_{\Gamma_0}$ is a Bowen--Series (respectively, higher Bowen--Series) map of $\Gamma_0$. To extend the notion of mateability to a group $\Gamma\in\partial\mathcal{B}(\Gamma_0)$, one needs the limit set of $\Gamma$ to carry a continuous, piecewise complex-analytic self-map $A_\Gamma$ (that is orbit equivalent to $\Gamma$) defined by the following commutative diagram:

\[ \begin{tikzcd}
\mathbb{S}^1 \arrow{r}{A_{\Gamma_0}} \arrow[swap]{d}{\mathrm{\phi_\Gamma}} & \mathbb{S}^1 \arrow{d}{\mathrm{\phi_\Gamma}} \\
\Lambda_\Gamma \arrow[swap]{r}{A_{\Gamma}}& \Lambda_\Gamma 
\end{tikzcd}
\]
See \cite[\S 7.1]{mj-muk} for details and an alternative description of $A_\Gamma$ as a uniform limit of Bowen--Series (respectively, higher Bowen--Series) maps. The map $A_\Gamma$, if it exists, is called the \emph{Bowen--Series} (respectively, \emph{higher Bowen--Series}) map of $\Gamma$ and can be thought of as a mateable map associated with a Bers boundary group.

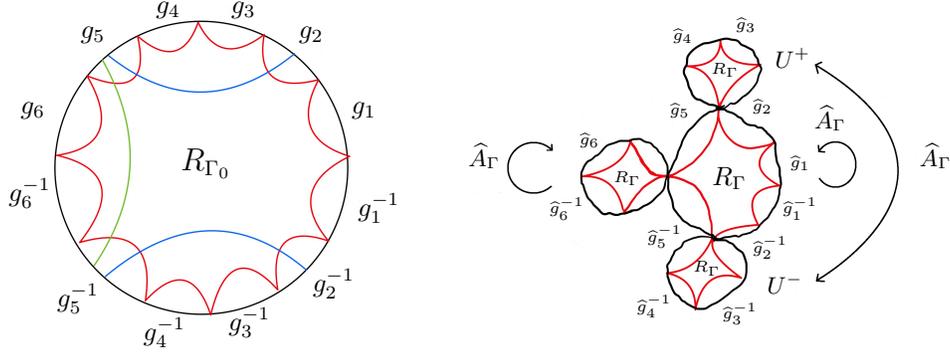
\begin{figure}[h!]
\begin{tikzpicture}
\node[anchor=south west,inner sep=0] at (0,0) {\includegraphics[width=0.42\textwidth]{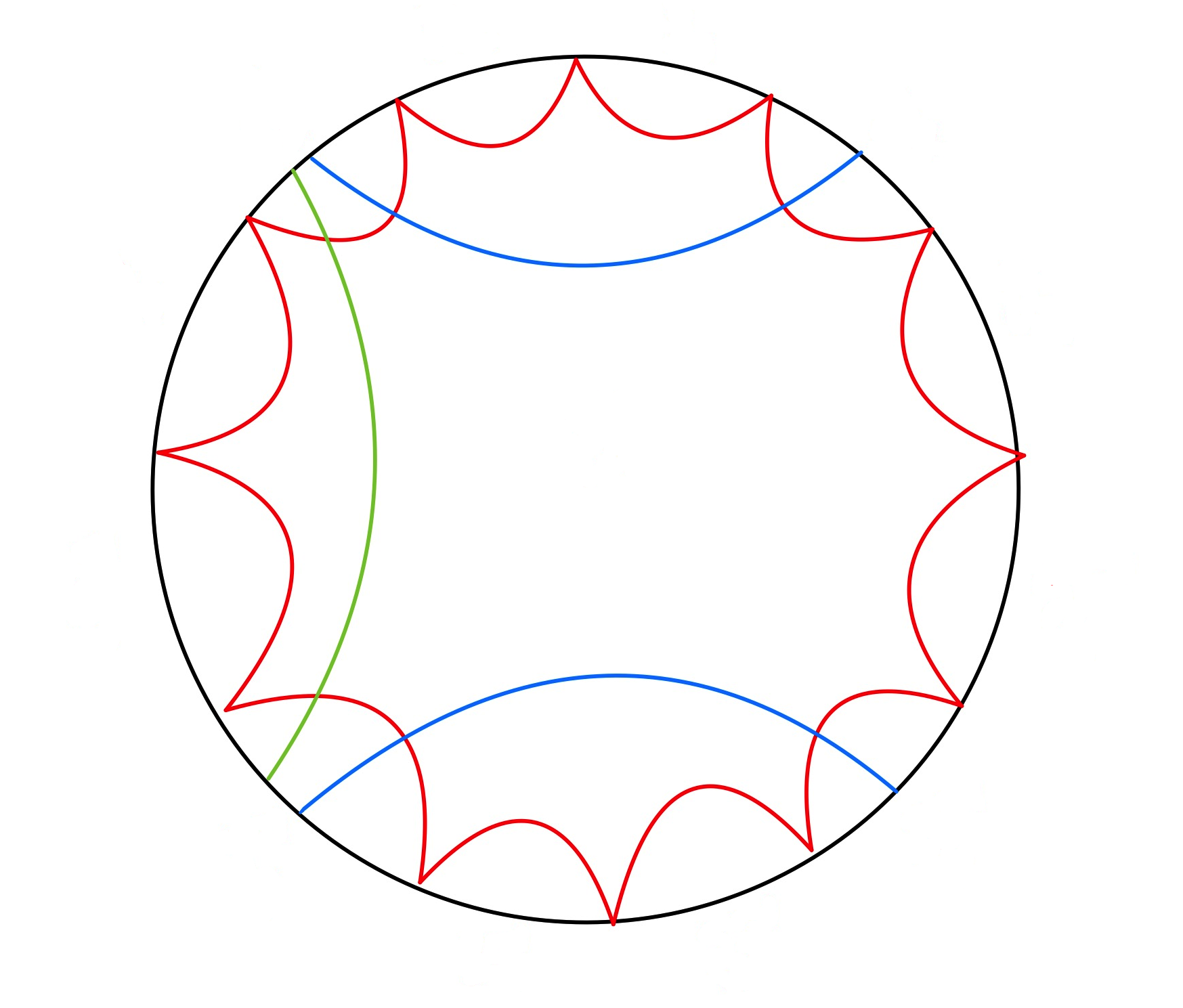}};
\node[anchor=south west,inner sep=0] at (5.4,0) {\includegraphics[width=0.64\textwidth]{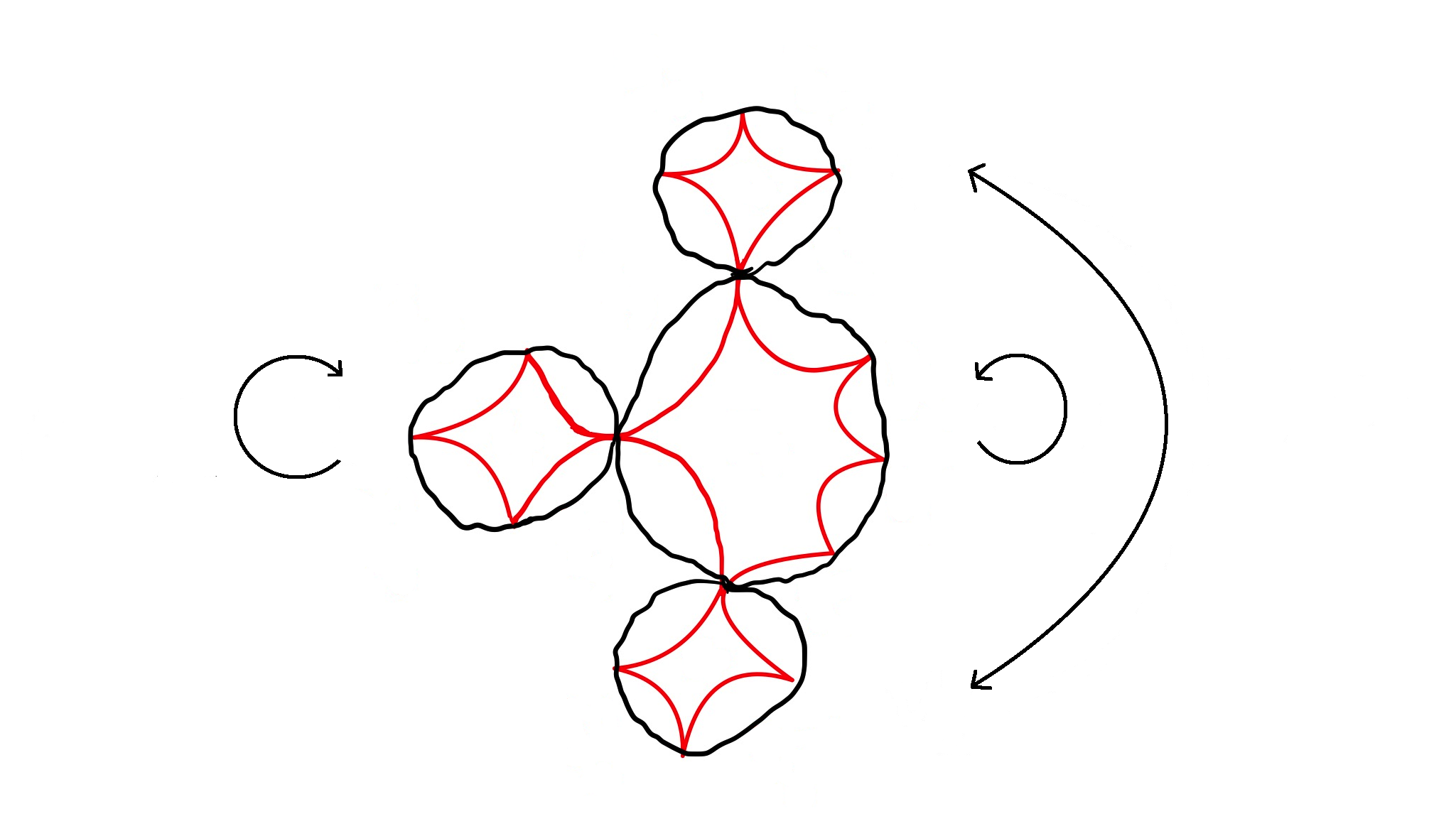}};
\node at (0.4,3.1) {$g_6$};
\node at (0.36,2) {$g_6^{-1}$};
\node at (1.2,4.1) {$g_5$};
\node at (1,0.6) {$g_5^{-1}$};
\node at (4.8,3.1) {$g_1$};
\node at (5,1.8) {$g_1^{-1}$};
\node at (4.1,4.1) {$g_2$};
\node at (4.4,0.8) {$g_2^{-1}$};
\node at (3.2,4.44) {$g_3$};
\node at (3.28,0.28) {$g_3^{-1}$};
\node at (2.2,4.44) {$g_4$};
\node at (2.15,0.15) {$g_4^{-1}$};
\node at (2.7,2.4) {\begin{Large}$R_{\Gamma_0}$\end{Large}};
\node at (10.5,3.8) {\begin{small}$U^+$\end{small}};
\node at (10.4,0.8) {\begin{small}$U^-$\end{small}};
\node at (9.66,2.2) {$R_\Gamma$};
\node at (8.3,2.2) {\begin{tiny}$R_\Gamma$\end{tiny}};
\node at (9.6,3.64) {\begin{tiny}$R_\Gamma$\end{tiny}};
\node at (9.36,1) {\begin{tiny}$R_\Gamma$\end{tiny}};
\node at (10.6,2.4) {\begin{tiny}$\widehat{g}_1$\end{tiny}};
\node at (10.6,1.8) {\begin{tiny}$\widehat{g}_1^{-1}$\end{tiny}};
\node at (10.1,3.16) {\begin{tiny}$\widehat{g}_2$\end{tiny}};
\node at (10.2,1.28) {\begin{tiny}$\widehat{g}_2^{-1}$\end{tiny}};
\node at (9.88,4.24) {\begin{tiny}$\widehat{g}_3$\end{tiny}};
\node at (9.8,0.4) {\begin{tiny}$\widehat{g}_3^{-1}$\end{tiny}};
\node at (9.04,4.1) {\begin{tiny}$\widehat{g}_4$\end{tiny}};
\node at (8.66,0.5) {\begin{tiny}$\widehat{g}_4^{-1}$\end{tiny}};
\node at (9,3.1) {\begin{tiny}$\widehat{g}_5$\end{tiny}};
\node at (8.8,1.42) {\begin{tiny}$\widehat{g}_5^{-1}$\end{tiny}};
\node at (7.8,2.7) {\begin{tiny}$\widehat{g}_6$\end{tiny}};
\node at (7.5,1.8) {\begin{tiny}$\widehat{g}_6^{-1}$\end{tiny}};
\node at (6.4,2.5) {\begin{small}$\widehat{A}_{\Gamma}$\end{small}};
\node at (12.4,2.5) {\begin{small}$\widehat{A}_{\Gamma}$\end{small}};
\node at (11,3) {\begin{small}$\widehat{A}_{\Gamma}$\end{small}};
\end{tikzpicture}
\caption{Left: $R_{\Gamma_0}$ is a fundamental domain of a Fuchsian group $\Gamma_0$ uniformizing $S_{0,7}$ with side pairing transformations $g_1^{\pm 1}, \cdots, g_6^{\pm 1}$. The geodesic lamination $\mathcal{L}^*$ on $\mathbb{D}/\Gamma_0$ consisting of two simple, closed curves corresponding to the elements $g_5, g_2g_5^{-1}\in \Gamma_0$ is admissible for the Bowen--Series map $A_{\Gamma_0}$. The blue and green geodesics are the connected components of the $\Gamma_0$-lift of $\mathcal{L}^*$ that intersect $R_{\Gamma_0}$.
Right: A cartoon of the limit set of a Bers boundary group $\Gamma$, which is obtained by pinching $\mathcal{L}^*$. The $\Gamma$-action on $\Omega(\Gamma)\setminus \Omega_\infty(\Gamma)$ admits a pinched fundamental domain $R_\Gamma$, and the canonical extension $\widehat{A}_\Gamma$ of the Bowen--Series map of $\Gamma$ is defined on $K(\Gamma)\setminus\textrm{int}~{R_\Gamma}$. Two of the components of $\Omega(\Gamma)$ intersecting $R_\Gamma$ are invariant under $\widehat{A}_\Gamma$, while the other two components $U^\pm$ form a $2$-cycle. The first return map of $\widehat{A}_{\Gamma}$ on $U^\pm$ is conformally conjugate to higher Bowen--Series maps of punctured sphere Fuchsian groups. The M{\"o}bius maps defining $\widehat{A}_{\Gamma}$ are also marked, where $\widehat{g}_i$ is the image of $g_i$ under the representation $\Gamma_0\rightarrow\Gamma$.}
\label{induced_bs_fig}
\end{figure}

It turns out that the existence of such a map $A_\Gamma$ imposes severe restrictions on the laminations that can be pinched. The next theorem says that only finitely many possibilities exist. We call such laminations \emph{admissible} (see Figure~\ref{induced_bs_fig} for an example of an admissible lamination in the Bowen--Series case). 

For $\Gamma\in\partial\mathcal{B}(\Gamma_0)$, we denote the unique $\Gamma$-invariant component of the domain of discontinuity $\Omega(\Gamma)$ by $\Omega_\infty(\Gamma)$, and set $K(\Gamma):=\widehat{\C}\setminus\Omega_\infty(\Gamma)$. If $\Gamma$ admits a Bowen--Series (respectively, higher Bowen--Series) map $A_\Gamma:\Lambda_\Gamma\to\Lambda_\Gamma$, then this map can be extended as a continuous, piecewise M{\"o}bius map to a canonical closed set $K(\Gamma)\setminus\textrm{int}~{R_\Gamma}$, where $R_\Gamma$ is a `pinched' fundamental domain for the $\Gamma-$action on $\Omega(\Gamma)\setminus \Omega_\infty(\Gamma)$ determined by $R_{\Gamma_0}$ (see \cite[\S 7.3]{mj-muk}). This canonical extension is denoted by $\widehat{A}_\Gamma$. The following theorem also demonstrates conformal mateability of groups $\Gamma\in\partial\mathcal{B}(\Gamma_0)$ admitting Bowen--Series/ higher Bowen--Series maps with polynomials lying in principal hyperbolic components (of suitable degree).

\begin{theorem}\cite[Lemma~7.3, Lemma~7.5, Theorem~7.19]{mj-muk}\label{mating_boundary_groups_thm} 
	Let $\Gamma_0$ be a Fuchsian group uniformizing $S_{0,k}$. Then, there are only finitely many quasiconformal conjugacy classes of groups $\Gamma\in\partial\mathcal{B}(\Gamma_0)$ for which the Cannon--Thurston map of $\Gamma$ semi-conjugates the Bowen--Series (respectively, higher Bowen--Series) map of $\Gamma_0$ to a self-map of $\Lambda(\Gamma)$ that is orbit equivalent to $\Gamma$. These Kleinian groups arise out of pinching finitely many disjoint, simple, closed curves (on the surface $\mathbb{D}/\Gamma_0$) out of an explicit finite list. In particular, all such groups $\Gamma$ are geometrically finite. 	
		
	Let $\Gamma\in\partial\mathcal{B}(\Gamma_0)$ be a group that admits a Bowen--Series (respectively, higher Bowen--Series) map $A_\Gamma$. Then the canonical extension $\widehat{A}_\Gamma:K(\Gamma)\setminus\textrm{int}~{R_\Gamma}\to K(\Gamma)$ can be conformally mated with polynomials lying in the principal hyperbolic component of degree $2k-3$ (respectively, $(k-1)^2$).
\end{theorem}

We refer the reader to \cite[Remark~7.20]{mj-muk} for a precise definition of conformal mateability of canonical extensions of the Bowen--Series/ higher Bowen--Series maps with polynomials lying in principal hyperbolic components (the definition is analogous to Definition~\ref{mating_def_2}).

The finiteness part of Theorem~\ref{mating_boundary_groups_thm} underscores the incompatibility between group invariant geodesic laminations and polynomial laminations (see \cite{kiwi1} for details on polynomial laminations) by establishing that the equivalence relation on $\mathbb{S}^1$ induced by a group invariant geodesic lamination on $\mathbb{D}$ is seldom invariant under $A_{\Gamma_0}$ (since $A_{\Gamma_0}\vert_{\mathbb{S}^1}$ is topologically conjugate to $z^d\vert_{\mathbb{S}^1}$ for some $d\geq 2$, invariance under $A_{\Gamma_0}$ should be thought of as $z^d-$invariance).

The proof of existence of a conformal mating between $\widehat{A}_\Gamma:K(\Gamma)\setminus\textrm{int}~{R_\Gamma}\to K(\Gamma)$ and polynomials lying in principal hyperbolic components has two main steps. The first one is to topologically realize the action of $A_\Gamma\vert_{\Lambda_\Gamma}$ by the dynamics of a postcritically finite polynomial $P_\Gamma$ on its Julia set, which is the content of \cite[Theorem~7.16]{mj-muk}. Once this is achieved, one needs to replace the dynamics of $P_\Gamma$ on periodic Fatou components by the action of Bowen--Series/ higher Bowen--Series maps of suitable punctured sphere Fuchsian groups. This involves a rather delicate surgery technique using David homeomorphisms.

\newcommand{\etalchar}[1]{$^{#1}$}

\end{document}